\documentclass{article}
\usepackage{amsmath}
\usepackage{amssymb}
\begin{document}
\def\e#1\e{\begin{equation}#1\end{equation}}
\def\ea#1\ea{\begin{align}#1\end{align}}
\def\eq#1{{\rm(\ref{#1})}}
\newtheorem{thm}{Theorem}[section]
\newtheorem{lem}[thm]{Lemma}
\newtheorem{prop}[thm]{Proposition}
\newtheorem{conj}[thm]{Conjecture}
\newtheorem{cor}[thm]{Corollary}
\newenvironment{dfn}{\medskip\refstepcounter{thm}
\noindent{\bf Definition \thesection.\arabic{thm}\ }}{\medskip}
\newenvironment{ex}{\medskip\refstepcounter{thm}
\noindent{\bf Example \thesection.\arabic{thm}\ }}{\medskip}
\newenvironment{proof}[1][,]{\medskip\ifcat,#1
\noindent{\it Proof.\ }\else\noindent{\it Proof of #1.\ }\fi}
{\hfill$\square$\medskip}
\def\dim{\mathop{\rm dim}}
\def\Re{\mathop{\rm Re}}
\def\Im{\mathop{\rm Im}}
\def\ind{\mathop{\rm ind}}
\def\vol{\mathop{\rm vol}}
\def\sind{{\ts\mathop{\text{\rm s-ind}}}}
\def\id{\mathop{\rm id}}
\def\U{\mathbin{\rm U}}
\def\SU{\mathop{\rm SU}}
\def\g{\mathfrak{g}} 
\def\su{\mathfrak{su}} 
\def\ge{\geqslant} 
\def\le{\leqslant} 
\def\R{\mathbin{\mathbb R}}
\def\N{\mathbin{\mathbb N}}
\def\Z{\mathbin{\mathbb Z}}
\def\C{\mathbin{\mathbb C}}
\def\al{\alpha}
\def\be{\beta}
\def\ga{\gamma}
\def\de{\delta}
\def\ep{\epsilon}
\def\la{\lambda}
\def\ka{\kappa}
\def\th{\theta}
\def\ze{\zeta}
\def\up{\upsilon}
\def\vp{\varphi}
\def\si{\sigma}
\def\om{\omega}
\def\De{\Delta}
\def\La{\Lambda}
\def\Om{\Omega}
\def\Ga{\Gamma}
\def\Si{\Sigma}
\def\Th{\Theta}
\def\Up{\Upsilon}
\def\d{{\rm d}}
\def\pd{\partial}
\def\ts{\textstyle}
\def\sst{\scriptscriptstyle}
\def\w{\wedge}
\def\lt{\ltimes}
\def\sm{\setminus}
\def\op{\oplus}
\def\ot{\otimes}
\def\iy{\infty}
\def\ra{\rightarrow}
\def\hookra{\hookrightarrow}
\def\t{\times}
\def\na{\nabla}
\def\ha{{\textstyle\frac{1}{2}}}
\def\ti{\tilde}
\def\bs{\boldsymbol}
\def\ov{\overline}
\def\sE{{\smash{\sst\cal E}}}
\def\sF{{\smash{\sst\cal F}}}
\def\sFp{{\smash{\sst{\cal F}'}}}
\def\sN{{\smash{\sst X}}}
\def\sX{{\smash{\sst X}}}
\def\sXp{{\smash{\sst X'}}}
\def\D{{\cal D}}
\def\E{{\cal E}}
\def\F{{\cal F}}
\def\H{{\cal H}}
\def\I{{\cal I}}
\def\K{{\cal K}}
\def\M{{\cal M}}
\def\O{{\cal O}}
\def\V{{\cal V}}
\def\ms#1{\vert#1\vert^2}
\def\bms#1{\bigl\vert#1\bigr\vert^2}
\def\md#1{\vert #1 \vert}
\def\bmd#1{\big\vert #1 \big\vert}
\def\nm#1{\Vert #1 \Vert}
\def\cnm#1#2{\Vert #1 \Vert_{C^{#2}}} 
\def\lnm#1#2{\Vert #1 \Vert_{L^{#2}}} 
\def\snm#1#2#3{\Vert #1 \Vert_{L^{#2}_{#3}}} 
\def\bnm#1{\bigl\Vert #1 \bigr\Vert}
\def\an#1{\langle#1\rangle}
\def\ban#1{\bigl\langle#1\bigr\rangle}
\title{Special Lagrangian submanifolds with isolated \\
conical singularities. II. Moduli spaces}
\author{Dominic Joyce \\ Lincoln College, Oxford}
\date{}
\maketitle

\section{Introduction}
\label{cm1}

{\it Special Lagrangian $m$-folds (SL\/ $m$-folds)} are a
distinguished class of real $m$-dimensional minimal submanifolds
which may be defined in $\C^m$, or in {\it Calabi--Yau $m$-folds},
or more generally in {\it almost Calabi--Yau $m$-folds} (compact
K\"ahler $m$-folds with trivial canonical bundle).

This is the second in a series of five papers
\cite{Joyc7,Joyc8,Joyc9,Joyc10} studying SL $m$-folds with
{\it isolated conical singularities}. That is, we consider an
SL $m$-fold $X$ in $M$ with singularities at $x_1,\ldots,x_n$
in $M$, such that for some SL cones $C_i$ in $T_{\smash{x_i}}M
\cong\C^m$ with $C_i\sm\{0\}$ nonsingular, $X$ approaches $C_i$
near $x_i$ in an asymptotic $C^1$ sense. Readers are advised to
begin with the final paper \cite{Joyc10}, which surveys the
series, and applies the results to prove some conjectures.

Having a good understanding of the singularities of special
Lagrangian submanifolds will be essential in clarifying the
Strominger--Yau--Zaslow conjecture on the Mirror Symmetry
of Calabi--Yau 3-folds \cite{SYZ}, and also in resolving
conjectures made by the author \cite{Joyc1} on defining
new invariants of Calabi--Yau 3-folds by counting special
Lagrangian homology 3-spheres with weights. The series aims
to develop such an understanding for simple singularities
of SL $m$-folds.

In this paper we study the {\it deformation theory} of
compact SL $m$-folds $X$ with conical singularities
$x_1,\ldots,x_n$ with cones $C_1,\ldots,C_n$ in an almost
Calabi--Yau $m$-fold $M$, extending results of McLean
\cite{McLe} for nonsingular compact SL $m$-folds. We
define the {\it moduli space} $\M_\sX$ of deformations
of $X$ as an SL $m$-fold with conical singularities in
$M$, and construct a natural {\it topology} on~$\M_\sX$.

We prove that $\M_\sX$ is locally homeomorphic to the
zeroes of a smooth map $\Phi:\I_\sXp\ra\O_\sXp$, where the
{\it infinitesimal deformation space} $\I_\sXp$ and the
{\it obstruction space} $\O_\sXp$ are finite-dimensional
vector spaces. Here $\I_\sXp$ depends only on the topology
of $X$, and $\O_\sXp$ only on the singular cones
$C_1,\ldots,C_n$. If $\O_\sXp$ is zero then $\M_\sX$ is a
{\it smooth manifold}. We also consider deformations of
$X$ in a {\it smooth family} of almost Calabi--Yau
$m$-folds~$\bigl\{(M,J^s,\om^s,\Om^s):s\in\F\bigr\}$.

The first paper \cite{Joyc7} laid the foundations for the
series, and studied the {\it regularity} of SL $m$-folds
with conical singularities near their singular points. The
sequels \cite{Joyc8,Joyc9} will consider
{\it desingularizations} of a compact SL $m$-fold $X$ with
conical singularities $x_1,\ldots,x_n$ with cones
$C_1,\ldots,C_n$ in $M$. We will take nonsingular SL
$m$-folds $L_1,\ldots,L_n$ in $\C^m$ asymptotic to
$C_1,\ldots,C_n$ at infinity, and glue them in to $X$ at
$x_1,\ldots,x_n$ to get a smooth family of compact,
{\it nonsingular} SL $m$-folds $\ti N$ in $M$ which
converge to~$X$.

We begin in \S\ref{cm2} with an introduction to special
Lagrangian geometry, and the deformation theory of nonsingular
compact SL $m$-folds. Section \ref{cm3} discusses {\it special
Lagrangian cones} and {\it conical singularities} of SL $m$-folds.
The previous paper \cite{Joyc7} is reviewed in \S\ref{cm4}. To
keep this paper and \cite{Joyc8,Joyc9} to a manageable length
we have done quite a lot of work on symplectic geometry and
asymptotic analysis in advance in \cite{Joyc7}, and we just
quote the results.

Section \ref{cm5} defines the moduli space $\M_\sX$ of SL
$m$-folds and its topology, and explains why this definition
of topology is a good one. In \S\ref{cm6} we define the
{\it infinitesimal deformation space} $\I_\sXp$ and the
{\it obstruction space} $\O_\sXp$, and prove our first main
result, Theorem \ref{cm6thm2}, which shows that the moduli
space $\M_\sX$ is locally homeomorphic to the zeroes of a
smooth map $\Phi:\I_\sXp\ra\O_\sXp$. Thus, if $\O_\sXp$ is zero
then $\M_\sX$ is a manifold. More generally, if $\d\Phi\vert_0$
is surjective then $\M_\sX$ is a manifold near~$X$.

Section \ref{cm7} extends \S\ref{cm5}--\S\ref{cm6} to {\it
families} $\bigl\{(M,J^s,\om^s,\Om^s):s\in\F\bigr\}$ of almost
Calabi--Yau $m$-folds. We define a {\it joint moduli space}
$\M_\sX^\sF$ with {\it projection} $\pi^\sF:\M_\sX^\sF\ra\F$
such that $\M_\sX^s=(\pi^\sF)^{-1}(s)$ is the moduli space of
deformations of $X$ in $(M,J^s,\om^s,\Om^s)$ for $s\in\F$.
Then we show that $\M_\sX^\sF$ is locally homeomorphic to the
zeroes of a smooth map $\Phi^\sF:\F\t\I_\sXp\ra\O_\sXp$, where
$\I_\sXp,\O_\sXp$ are as before.

Section \ref{cm8} briefly describes various other extensions
of the results to immersions, families of SL cones in $\C^m$,
and so on. Finally, \S\ref{cm9} considers {\it genericity} and
{\it transversality} results. We show that for any compact SL
$m$-fold $X$ with conical singularities in $(M,J,\om,\Om)$, we
can choose a family of deformations $\bigl\{(M,J,\om^s,\Om):s\in
\F\bigr\}$ such that $\M_\sX^\sF$ is a manifold near $(0,X)$, and
for small generic $s\in\F$ the deformed moduli space $\M_\sX^s=
(\pi^\sF)^{-1}(s)$ is smooth near $(0,X)$. We conjecture that
if the K\"ahler form $\om$ is chosen {\it generically in its
K\"ahler class}, then $\M_\sX$ is smooth.
\medskip

\noindent{\it Acknowledgements.} I would like to thank Stephen
Marshall, Mark Haskins and Nigel Hitchin for helpful conversations.
I was supported by an EPSRC Advanced Research Fellowship whilst
writing this paper.

\section{Special Lagrangian geometry}
\label{cm2}

We now introduce special Lagrangian submanifolds (SL $m$-folds)
in two different geometric contexts. First, in \S\ref{cm21}, we
define SL $m$-folds in $\C^m$. Then \S\ref{cm22} discusses SL
$m$-folds in {\it almost Calabi--Yau $m$-folds}, compact K\"ahler
manifolds with a holomorphic volume form, which generalize
Calabi--Yau manifolds. Section \ref{cm23} describes the
{\it deformation theory} of compact SL $m$-folds. Some
references for this section are Harvey and Lawson \cite{HaLa},
McLean \cite{McLe}, and the author~\cite{Joyc6}.

\subsection{Special Lagrangian submanifolds in $\C^m$}
\label{cm21}

We begin by defining {\it calibrations} and {\it calibrated 
submanifolds}, following Harvey and Lawson~\cite{HaLa}.

\begin{dfn} Let $(M,g)$ be a Riemannian manifold. An {\it oriented
tangent $k$-plane} $V$ on $M$ is a vector subspace $V$ of
some tangent space $T_xM$ to $M$ with $\dim V=k$, equipped
with an orientation. If $V$ is an oriented tangent $k$-plane
on $M$ then $g\vert_V$ is a Euclidean metric on $V$, so 
combining $g\vert_V$ with the orientation on $V$ gives a 
natural {\it volume form} $\vol_V$ on $V$, which is a 
$k$-form on~$V$.

Now let $\vp$ be a closed $k$-form on $M$. We say that
$\vp$ is a {\it calibration} on $M$ if for every oriented
$k$-plane $V$ on $M$ we have $\vp\vert_V\le \vol_V$. Here
$\vp\vert_V=\al\cdot\vol_V$ for some $\al\in\R$, and 
$\vp\vert_V\le\vol_V$ if $\al\le 1$. Let $N$ be an 
oriented submanifold of $M$ with dimension $k$. Then 
each tangent space $T_xN$ for $x\in N$ is an oriented
tangent $k$-plane. We say that $N$ is a {\it calibrated 
submanifold\/} if $\vp\vert_{T_xN}=\vol_{T_xN}$ for all~$x\in N$.
\label{cm2def1}
\end{dfn}

It is easy to show that calibrated submanifolds are automatically
{\it minimal submanifolds} \cite[Th.~II.4.2]{HaLa}. Here is the 
definition of special Lagrangian submanifolds in $\C^m$, taken
from~\cite[\S III]{HaLa}.

\begin{dfn} Let $\C^m$ have complex coordinates $(z_1,\dots,z_m)$, 
and define a metric $g'$, a real 2-form $\om'$ and a complex $m$-form 
$\Om'$ on $\C^m$ by
\e
\begin{split}
g'=\ms{\d z_1}+\cdots+\ms{\d z_m},\quad
\om'&=\ts\frac{i}{2}(\d z_1\w\d\bar z_1+\cdots+\d z_m\w\d\bar z_m),\\
\text{and}\quad\Om'&=\d z_1\w\cdots\w\d z_m.
\end{split}
\label{cm2eq1}
\e
Then $\Re\Om'$ and $\Im\Om'$ are real $m$-forms on $\C^m$. Let $L$
be an oriented real submanifold of $\C^m$ of real dimension $m$. We
say that $L$ is a {\it special Lagrangian submanifold\/} of $\C^m$,
or {\it SL\/ $m$-fold}\/ for short, if $L$ is calibrated with respect
to $\Re\Om'$, in the sense of Definition~\ref{cm2def1}.
\label{cm2def2}
\end{dfn}

Harvey and Lawson \cite[Cor.~III.1.11]{HaLa} give the following
alternative characterization of special Lagrangian submanifolds:

\begin{prop} Let\/ $L$ be a real $m$-dimensional submanifold 
of\/ $\C^m$. Then $L$ admits an orientation making it into an
SL submanifold of\/ $\C^m$ if and only if\/ $\om'\vert_L\equiv 0$ 
and\/~$\Im\Om'\vert_L\equiv 0$.
\label{cm2prop1}
\end{prop}

Thus SL $m$-folds are {\it Lagrangian submanifolds} in
$\R^{2m}\cong\C^m$ satisfying the extra condition that
$\Im\Om'\vert_L\equiv 0$, which is how they get their name.

\subsection{Almost Calabi--Yau $m$-folds and SL $m$-folds} 
\label{cm22}

We shall define special Lagrangian submanifolds not just in
Calabi--Yau manifolds, as usual, but in the much larger
class of {\it almost Calabi--Yau manifolds}.

\begin{dfn} Let $m\ge 2$. An {\it almost Calabi--Yau $m$-fold\/}
is a quadruple $(M,J,\om,\Om)$ such that $(M,J)$ is a compact
$m$-dimensional complex manifold, $\om$ is the K\"ahler form
of a K\"ahler metric $g$ on $M$, and $\Om$ is a non-vanishing
holomorphic $(m,0)$-form on~$M$.

We call $(M,J,\om,\Om)$ a {\it Calabi--Yau $m$-fold\/} if in
addition $\om$ and $\Om$ satisfy
\e
\om^m/m!=(-1)^{m(m-1)/2}(i/2)^m\Om\w\bar\Om.
\label{cm2eq2}
\e
Then for each $x\in M$ there exists an isomorphism $T_xM\cong\C^m$
that identifies $g_x,\om_x$ and $\Om_x$ with the flat versions
$g',\om',\Om'$ on $\C^m$ in \eq{cm2eq1}. Furthermore, $g$ is
Ricci-flat and its holonomy group is a subgroup of~$\SU(m)$.
\label{cm2def3}
\end{dfn}

This is not the usual definition of a Calabi--Yau manifold, but
is essentially equivalent to it.

\begin{dfn} Let $(M,J,\om,\Om)$ be an almost Calabi--Yau $m$-fold,
and $N$ a real $m$-dimensional submanifold of $M$. We call $N$ a
{\it special Lagrangian submanifold}, or {\it SL $m$-fold\/} for
short, if $\om\vert_N\equiv\Im\Om\vert_N\equiv 0$. It easily
follows that $\Re\Om\vert_N$ is a nonvanishing $m$-form on $N$.
Thus $N$ is orientable, with a unique orientation in which
$\Re\Om\vert_N$ is positive.
\label{cm2def4}
\end{dfn}

Again, this is not the usual definition of SL $m$-fold, but is
essentially equivalent to it. Suppose $(M,J,\om,\Om)$ is an
almost Calabi--Yau $m$-fold, with metric $g$. Let
$\psi:M\ra(0,\iy)$ be the unique smooth function such that
\e
\psi^{2m}\om^m/m!=(-1)^{m(m-1)/2}(i/2)^m\Om\w\bar\Om,
\label{cm2eq3}
\e
and define $\ti g$ to be the conformally equivalent metric $\psi^2g$
on $M$. Then $\Re\Om$ is a {\it calibration} on the Riemannian manifold
$(M,\ti g)$, and SL $m$-folds $N$ in $(M,J,\om,\Om)$ are calibrated
with respect to it, so that they are minimal with respect to~$\ti g$.

If $M$ is a Calabi--Yau $m$-fold then $\psi\equiv 1$ by \eq{cm2eq2},
so $\ti g=g$, and an $m$-submanifold $N$ in $M$ is special Lagrangian
if and only if it is calibrated w.r.t.\ $\Re\Om$ on $(M,g)$, as in
Definition \ref{cm2def2}. This recovers the usual definition of
special Lagrangian $m$-folds in Calabi--Yau $m$-folds.

\subsection{Deformations of compact SL $m$-folds} 
\label{cm23}

The {\it deformation theory} of special Lagrangian submanifolds
was studied by McLean \cite[\S 3]{McLe}, who proved the following
result in the Calabi--Yau case. The extension to the almost
Calabi--Yau case is described in~\cite[\S 9.5]{Joyc6}.

\begin{thm} Let\/ $N$ be a compact SL\/ $m$-fold in an almost
Calabi--Yau $m$-fold\/ $(M,J,\om,\Om)$. Then the moduli space
$\M_\sN$ of special Lagrangian deformations of\/ $N$ is a
smooth manifold of dimension $b^1(N)$, the first Betti number of\/~$N$.
\label{cm2thm1}
\end{thm}

We now give a partial proof of Theorem \ref{cm2thm1}, glossing over
the analytic details, and concentrating on the parts we will use
later. We start by recalling some symplectic geometry, which can
be found in McDuff and Salamon~\cite{McSa2}.

Let $N$ be a real $m$-manifold. Then its tangent bundle $T^*N$ has
a {\it canonical symplectic form} $\hat\om$, defined as follows.
Let $(x_1,\ldots,x_m)$ be local coordinates on $N$. Extend them to
local coordinates $(x_1,\ldots,x_m,y_1,\ldots,y_m)$ on $T^*N$ such
that $(x_1,\ldots,y_m)$ represents the 1-form $y_1\d x_1+\cdots+y_m
\d x_m$ in $T_{(x_1,\ldots,x_m)}^*N$. Then~$\hat\om=\d x_1\w\d y_1+
\cdots+\d x_m\w\d y_m$.

Identify $N$ with the zero section in $T^*N$. Then $N$ is a
{\it Lagrangian submanifold\/} of $T^*N$. The {\it Lagrangian
Neighbourhood Theorem} \cite[Th.~3.33]{McSa2} shows that any
compact Lagrangian submanifold $N$ in a symplectic manifold
looks locally like the zero section in~$T^*N$.

\begin{thm} Let\/ $(M,\om)$ be a symplectic manifold and\/
$N\subset M$ a compact Lagrangian submanifold. Then there
exists an open tubular neighbourhood\/ $U$ of the zero
section $N$ in $T^*N$, and an embedding $\Phi:U\ra M$ with\/
$\Phi\vert_N=\id:N\ra N$ and\/ $\Phi^*(\om)=\hat\om$, where
$\hat\om$ is the canonical symplectic structure on~$T^*N$.
\label{cm2thm2}
\end{thm}

In the situation of Theorem \ref{cm2thm1}, let $g$ be the K\"ahler
metric on $M$, and define $\psi:M\ra(0,\iy)$ by \eq{cm2eq3}.
Applying Theorem \ref{cm2thm2} gives an open neighbourhood $U$ of
$N$ in $T^*N$ and an embedding $\Phi:U\ra M$. Let $\pi:U\ra N$ be
the natural projection. Define an $m$-form $\be$ on $U$ by $\be=
\Phi^*(\Im\Om)$. If $\al$ is a 1-form on $N$ let $\Ga(\al)$ be
the graph of $\al$ in $T^*N$, and write $C^\iy(U)\subset
C^\iy(T^*N)$ for the subset of 1-forms whose graphs lie in~$U$.

Then each submanifold $\ti N$ of $M$ which is $C^1$-close to $N$
is $\Phi(\Ga(\al))$ for some small $\al\in C^\iy(U)$. Here is the
condition for $\ti N$ to be special Lagrangian.

\begin{lem} In the situation above, if\/ $\al\in C^\iy(U)$
then $\ti N=\Phi\bigl(\Ga(\al)\bigr)$ is a special
Lagrangian $m$-fold in $M$ if and only if\/ $\d\al=0$
and\/~$\pi_*\bigl(\be\vert_{\Ga(\al)}\bigr)=0$.
\label{cm2lem1}
\end{lem}

\begin{proof} By Definition \ref{cm2def4}, $\ti N$ is an SL $m$-fold
in $M$ if and only if $\om\vert_{\ti N}\equiv\Im\Om\vert_{\ti N}
\equiv 0$. Pulling back by $\Phi$ and pushing forward by $\pi:
\Ga(\al)\ra N$, we see that $\ti N$ is special Lagrangian if
and only if $\pi_*\bigl(\hat\om\vert_{\Ga(\al)}\bigr)\equiv
\pi_*\bigl(\be\vert_{\Ga(\al)}\bigr)\equiv 0$, since
$\Phi^*(\om)=\hat\om$ and $\Phi^*(\Im\Om)=\be$. But as
$\hat\om$ is the natural symplectic structure on $U\subset
T^*N$ we have $\pi_*\bigl(\hat\om\vert_{\Ga(\al)}\bigr)=-\d\al$,
and the lemma follows.
\end{proof}

We rewrite the condition $\pi_*\bigl(\be\vert_{\Ga(\al)}\bigr)=0$
in terms of a function~$F$.

\begin{dfn} Define $F:C^\iy(U)\ra C^\iy(N)$ by $\pi_*\bigl(
\be\vert_{\Ga(\al)}\bigr)=F(\al)\,\d V_g$, where $\d V_g$ is the
volume form of $g\vert_N$ on $N$. Then Lemma \ref{cm2lem1}
shows that if $\al\in C^\iy(U)$ then $\Phi\bigl(\Ga(\al)\bigr)$
is special Lagrangian if and only if~$\d\al=F(\al)=0$.
\label{cm2def5}
\end{dfn}

We compute the expansion of $F$ up to first order in~$\al$.

\begin{prop} This function $F$ may be written
\e
F(\al)[x]=-\d^*\bigl(\psi^m\al\bigr)+Q\bigl(x,\al(x),\na\al(x)\bigr)
\quad\text{for $x\in N$,}
\label{cm2eq4}
\e
where $Q:\bigl\{(x,y,z):x\in N$, $y\in T^*_xN\cap U$,
$z\in\ot^2T^*_xN\bigr\}\ra\R$ is smooth and\/ $Q(x,y,z)
=O(\ms{y}+\ms{z})$ for small\/~$y,z$.
\label{cm2prop2}
\end{prop}

\begin{proof} The value of $F(\al)$ at $x\in N$ depends on the
tangent space $T_{x'}\Ga(\al)$, where $x'\in\Ga(\al)$ with $\pi(x')=x$.
But $T_{x'}\Ga(\al)$ depends on both $\al\vert_x$ and $\na\al\vert_x$.
Hence $F(\al)$ depends pointwise on both $\al$ and $\na\al$, rather
than just $\al$. So we may take \eq{cm2eq4} as a {\it definition}
of $Q$, and $Q$ is then well-defined on the set of all $(x,y,z)$
realized by $\bigl(x,\al(x),\na\al(x)\bigr)$ for $\al\in C^\iy(U)$,
which is the domain given for~$Q$.

As $F$ depends smoothly on $\al$ we see that $Q$ is a smooth
function of its arguments. Therefore Taylor's theorem yields
\begin{equation*}
Q(x,y,z)=Q(x,0,0)+y\cdot(\pd_yQ)(x,0,0)+z\cdot
(\pd_zQ)(x,0,0)+O(\ms{y}+\ms{z})
\end{equation*}
for small $y,z$. So to prove that $Q(x,y,z)=
O(\ms{y}+\ms{z})$ we just need to show that $Q(x,0,0)=
\pd_yQ(x,0,0)=\pd_zQ(x,0,0)=0$. Now $N=\Phi(\Ga(0))$
is special Lagrangian, so $\al=0$ satisfies $F(\al)=0$ by
Definition \ref{cm2def5}. Thus \eq{cm2eq4} gives~$Q(x,0,0)\equiv 0$.

To compute $\pd_yQ(x,0,0)$ and $\pd_zQ(x,0,0)$, let
$\al\in C^\iy(U)$ be small, and let $v$ be the vector field on
$T^*N$ with $v\cdot\hat\om=-\pi^*(\al)$. Then $v$ is tangent to
the fibres of $\pi:T^*N\ra N$, and $\exp(v)$ maps $T^*N\ra T^*N$
taking $\ga\mapsto\al+\ga$ for 1-forms $\ga$ on $N$. Identifying
$N$ with the zero section of $T^*N$, the image $\exp(sv)[N]$ of
$N$ under $\exp(sv)$ is $\Ga(s\al)$ for~$s\in[0,1]$.

Therefore $F(s\al)\,\d V_g=\exp(sv)^*(\be)$ for $s\in[0,1]$.
Differentiating gives
\e
\begin{split}
\d F\vert_0(\al)\,\d V_g&=\frac{\d}{\d s}\bigl(F(s\al)\bigr)\Big\vert_{s=0}
\d V_g=\frac{\d}{\d s}\bigl(\exp(sv)^*(\be)\bigr)\Big\vert_{s=0}\\
&=\bigl({\cal L}_v\be\bigr)\Big\vert_N
=\bigl(\d(v\cdot\be)+v\cdot(\d\be)\bigr)\Big\vert_N
=\d\bigl((v\cdot\be)\vert_N\bigr),
\end{split}
\label{cm2eq5}
\e
where ${\cal L}_v$ is the Lie derivative, `$\,\cdot\,$'
contracts together vector fields and forms, and $\d\be=0$
as $\Om$ is closed and~$\be=\Phi^*(\Im\Om)$.

Calculation at a point $x\in N$ shows that $(v\cdot\be)\vert_N
=\psi^m*\al$, where $*$ is the Hodge star of $g$ on $N$. As
$*\d V_g=1$ and $*\d *=-\d^*$ on 1-forms, \eq{cm2eq5} gives
\begin{equation*}
\d F\vert_0(\al)\,\d V_g=\d(\psi^m*\al)=\bigl(*\d*(\psi^m\al)\bigr)\,\d V_g
=\bigl(-\d^*(\psi^m\al)\bigr)\,\d V_g.
\end{equation*}
Comparing this with \eq{cm2eq4} shows that $\pd_yQ(x,0,0)=
\pd_zQ(x,0,0)=0$, which completes the proof.
\end{proof}

We briefly sketch the remainder of the proof of Theorem
\ref{cm2thm1}. From Definition \ref{cm2def5} and Proposition
\ref{cm2prop2} we see that $\M_\sN$ is locally
approximately isomorphic to the vector space of 1-forms $\al$
with $\d\al=\d^*(\psi^m\al)=0$. But by Hodge theory, this is
isomorphic to the de Rham cohomology group $H^1(N,\R)$, and
is a manifold with dimension~$b^1(N)$.

To carry out this last step rigorously requires some technical
machinery: one must work with certain {\it Banach spaces} of 
sections of $\La^kT^*N$ for $k=0,1,2$, use {\it elliptic
regularity results} to prove that the map $\al\mapsto\bigl(\d\al,
\d F\vert_0(\al)\bigr)$ is {\it surjective} upon the appropriate
Banach spaces, and then use the {\it Implicit Mapping Theorem
for Banach spaces} to show that the kernel of the map is what
we expect.

Finally we extend of Theorem \ref{cm2thm1} to {\it families}
of almost Calabi--Yau $m$-folds.

\begin{dfn} Let $(M,J,\om,\Om)$ be an almost Calabi--Yau
$m$-fold. A {\it smooth family of deformations of\/}
$(M,J,\om,\Om)$ is a connected open set $\F\subset\R^d$
for $d\ge 0$ with $0\in\F$ called the {\it base space},
and a smooth family $\bigl\{(M,J^s,\om^s,\Om^s):s\in{\cal
F}\bigr\}$ of almost Calabi--Yau structures on $M$
with~$(J^0,\om^0,\Om^0)=(J,\om,\Om)$.
\label{cm2def6}
\end{dfn}

If $N$ is a compact SL $m$-fold in $(M,J,\om,\Om)$, the moduli
of deformations of $N$ in each $(M,J^s,\om^s,\Om^s)$ for $s\in\F$
make up a big moduli space~$\M_\sN^\sF$.

\begin{dfn} Let $\bigl\{(M,J^s,\om^s,\Om^s):s\in\F\bigr\}$
be a smooth family of deformations of an almost Calabi--Yau
$m$-fold $(M,J,\om,\Om)$, and $N$ be a compact SL $m$-fold
in $(M,J,\om,\Om)$. Define the {\it moduli space $\M_\sN^\sF$
of deformations of\/ $N$ in the family} $\F$ to be the set
of pairs $(s,\hat N)$ for which $s\in\F$ and $\hat N$ is a
compact SL $m$-fold in $(M,J^s,\om^s,\Om^s)$ which is
diffeomorphic to $N$ and isotopic to $N$ in $M$. Define a
{\it projection} $\pi^\sF:\M_\sN^\sF\ra\F$ by
$\pi^\sF(s,\hat N)=s$. Then $\M_\sN^\sF$ has a
natural topology, and $\pi^\sF$ is continuous.
\label{cm2def7}
\end{dfn}

The following result is proved by Marshall \cite[Th.~3.2.9]{Mars},
using similar methods to Theorem~\ref{cm2thm1}.

\begin{thm} Let\/ $\bigl\{(M,J^s,\om^s,\Om^s):s\in\F\bigr\}$ be
a smooth family of deformations of an almost Calabi--Yau $m$-fold\/
$(M,J,\om,\Om)$, with base space $\F\subset\R^d$. Suppose $N$ is a
compact SL\/ $m$-fold in $(M,J,\om,\Om)$ with\/ $[\om^s\vert_N]=0$
in $H^2(N,\R)$ and\/ $[\Im\Om^s\vert_N]=0$ in $H^m(N,\R)$ for all\/
$s\in\F$. Let\/ $\M_\sN^\sF$ be the moduli space of deformations
of\/ $N$ in $\F$, and\/ $\pi^\sF:\M_\sN^\sF\ra\F$ the natural
projection.

Then $\M_\sN^\sF$ is a smooth manifold of dimension\/ $d+b^1(N)$,
and\/ $\pi^\sF:\M_\sN^\sF\ra\F$ a smooth submersion. For small
$s\in\F$ the moduli space $\M_\sN^s=(\pi^\sF)^{-1}(s)$ of
deformations of\/ $N$ in $(M,J^s,\om^s,\Om^s)$ is a nonempty
smooth manifold of dimension $b^1(N)$, with\/~$\M_\sN^0=\M_\sN$.
\label{cm2thm3}
\end{thm}

Here a necessary condition for the existence of an SL $m$-fold
$\hat N$ isotopic to $N$ in $(M,J^s,\om^s,\Om^s)$ is that
$[\om^s\vert_N]=[\Im\Om^s\vert_N]=0$ in $H^*(N,\R)$, since
$[\om^s\vert_N]$ and $[\om^s\vert_{\smash{\hat N}}]$ are
identified under the natural isomorphism between $H^2(N,\R)$
and $H^2(\hat N,\R)$, and similarly for~$\Im\Om^s$.

The point of the theorem is that these conditions $[\om^s\vert_N]
=[\Im\Om^s\vert_N]=0$ are also {\it sufficient\/} for the
existence of $\hat N$ when $s$ is close to 0 in $\F$. That is,
the only {\it obstructions} to existence of compact SL $m$-folds
when we deform the underlying almost Calabi--Yau $m$-fold are
the obvious cohomological ones.

\section{SL cones and conical singularities}
\label{cm3}

After some preliminary work in \S\ref{cm31} on {\it special
Lagrangian cones}, and some examples in \S\ref{cm32}, section
\ref{cm33} defines {\it special Lagrangian $m$-folds with
conical singularities} in almost Calabi--Yau manifolds,
which are the subject of the paper.

\subsection{Preliminaries on special Lagrangian cones}
\label{cm31}

We now give some definitions and results on {\it special
Lagrangian cones}. Some are quoted from \cite{Joyc7}, and
some are new.

\begin{dfn} A (singular) SL $m$-fold $C$ in $\C^m$ is called a
{\it cone} if $C=tC$ for all $t>0$, where $tC=\{t{\bf x}:{\bf x}
\in C\}$. Let $C$ be a closed SL cone in $\C^m$ with an isolated
singularity at 0. Then $\Si=C\cap{\cal S}^{2m-1}$ is a compact,
nonsingular $(m\!-\!1)$-submanifold of ${\cal S}^{2m-1}$, not
necessarily connected. Let $g_{\sst\Si}$ be the restriction
of $g'$ to $\Si$, where $g'$ is as in~\eq{cm2eq1}.

Set $C'=C\sm\{0\}$. Define $\iota:\Si\t(0,\iy)\ra\C^m$ by
$\iota(\si,r)=r\si$. Then $\iota$ has image $C'$. By an abuse
of notation, {\it identify} $C'$ with $\Si\t(0,\iy)$ using
$\iota$. The {\it cone metric} on $C'\cong\Si\t(0,\iy)$
is~$g'=\iota^*(g')=\d r^2+r^2g_{\sst\Si}$.

For $\al\in\R$, we say that a function $u:C'\ra\R$ is
{\it homogeneous of order} $\al$ if $u\circ t\equiv t^\al u$ for
all $t>0$. Equivalently, $u$ is homogeneous of order $\al$ if
$u(\si,r)\equiv r^\al v(\si)$ for some function~$v:\Si\ra\R$.
\label{cm3def1}
\end{dfn}

In \cite[Lem.~2.3]{Joyc7} we study {\it homogeneous harmonic
functions} on~$C'$.

\begin{lem} In the situation of Definition \ref{cm3def1},
let\/ $u(\si,r)\equiv r^\al v(\si)$ be a homogeneous function
of order $\al$ on $C'=\Si\t(0,\iy)$, for $v\in C^2(\Si)$. Then
\begin{equation*}
\De u(\si,r)=r^{\al-2}\bigl(\De_{\sst\Si} v-\al(\al+m-2)v\bigr),
\end{equation*}
where $\De$, $\De_{\sst\Si}$ are the Laplacians on $(C',g')$
and\/ $(\Si,g_{\sst\Si})$. Hence, $u$ is harmonic on $C'$ if
and only if\/ $v$ is an eigenfunction of\/ $\De_{\sst\Si}$
with eigenvalue~$\al(\al+m-2)$.
\label{cm3lem1}
\end{lem}

Following \cite[Def.~2.5]{Joyc7}, we define:

\begin{dfn} In the situation of Definition \ref{cm3def1},
suppose $m>2$ and define
\e
\D_{\sst\Si}=\bigl\{\al\in\R:\text{$\al(\al+m-2)$ is
an eigenvalue of $\De_{\sst\Si}$}\bigr\}.
\label{cm3eq2}
\e
Then $\D_{\sst\Si}$ is a countable, discrete subset of
$\R$. By Lemma \ref{cm3lem1}, an equivalent definition is that
$\D_{\sst\Si}$ is the set of $\al\in\R$ for which there
exists a nonzero homogeneous harmonic function $u$ of order
$\al$ on~$C'$.

Define $m_{\sst\Si}:\D_{\sst\Si}\ra\N$ by taking
$m_{\sst\Si}(\al)$ to be the multiplicity of the eigenvalue
$\al(\al+m-2)$ of $\De_{\sst\Si}$, or equivalently the
dimension of the vector space of homogeneous harmonic
functions $u$ of order $\al$ on $C'$. Define
$N_{\sst\Si}:\R\ra\Z$ by
\e
N_{\sst\Si}(\de)=
-\sum_{\!\!\!\!\al\in\D_{\sst\Si}\cap(\de,0)\!\!\!\!}m_{\sst\Si}(\al)
\;\>\text{if $\de<0$, and}\;\>
N_{\sst\Si}(\de)=
\sum_{\!\!\!\!\al\in\D_{\sst\Si}\cap[0,\de]\!\!\!\!}m_{\sst\Si}(\al)
\;\>\text{if $\de\ge 0$.}
\label{cm3eq3}
\e
Then $N_{\sst\Si}$ is monotone increasing and upper semicontinuous,
and is discontinuous exactly on $\D_{\sst\Si}$, increasing by
$m_{\sst\Si}(\al)$ at each $\al\in\D_{\sst\Si}$. As the
eigenvalues of $\De_{\sst\Si}$ are nonnegative, we see that
$\D_{\sst\Si}\cap(2-m,0)=\emptyset$ and $N_{\sst\Si}\equiv 0$
on~$(2-m,0)$.
\label{cm3def2}
\end{dfn}

We shall show that there automatically exist homogeneous
harmonic functions on $C'$ of orders 1 and 2, using the
idea of {\it moment map}. The group of automorphisms of
$\C^m$ preserving $g',\om'$ and $\Om'$ is $\SU(m)\lt\C^m$,
where $\C^m$ acts by translations. Its Lie algebra
$\su(m)\lt\C^m$ acts on $\C^m$ by vector fields.

Let $v$ be such a vector field in $\su(m)\lt\C^m$. Then
$v\cdot\om'$ is a closed 1-form on $\C^m$, and we may write
$v\cdot\om'=\d\mu$ for some function $\mu:\C^m\ra\R$, which
is unique up to addition of constants, and is in fact a real
quadratic polynomial. We call $\mu$ a {\it moment map} for~$v$.

\begin{lem} Let\/ $L$ be an SL $m$-fold in $\C^m$, and
let\/ $\mu:\C^m\ra\R$ be a moment map for a vector field\/
$v$ in $\su(m)\lt\C^m$. Then $\mu\vert_L$ is a harmonic
function on $L$, using the obvious metric~$g'\vert_L$.
\label{cm3lem2} 
\end{lem}

\begin{proof} In the proof of Theorem \ref{cm2thm1} we
saw that infinitesimal deformations of an SL $m$-fold
$L$ as a submanifold correspond to 1-forms $\al$ on $L$,
and infinitesimal deformations as an SL $m$-fold to
closed and coclosed 1-forms $\al$ on~$L$.

Now as $\SU(m)\lt\C^m$ takes SL $m$-folds in $\C^m$ to SL
$m$-folds in $\C^m$, the vector field $v$ in $\su(m)\lt\C^m$
gives an infinitesimal deformation of $L$ as an SL $m$-fold
in $\C^m$. It is easy to see that the corresponding 1-form
on $L$ is $(v\cdot\om)\vert_L$. Therefore $(v\cdot\om)\vert_L
=\d\mu\vert_L$ is a closed and coclosed 1-form on $L$, and
thus $\d^*(\d\mu\vert_L)=0$, so $\mu\vert_L$ is harmonic.
\end{proof}

\begin{prop} Let\/ $C$ be an SL cone in $\C^m$ with isolated
singularity at\/ $0$, and\/ $G$ the Lie subgroup of\/ $\SU(m)$
preserving $C$. Set\/ $C'=C\sm\{0\}$ and\/ $\Si=C\cap{\cal S}^{2m-1}$,
and let\/ $m_{\sst\Si}$ be as in Definition \ref{cm3def2}. Then
\begin{itemize}
\item[{\rm(a)}] The restriction of real linear functions on $\C^m$
to $C'$ form a vector space of order $1$ homogeneous harmonic
functions on $C'$, with dimension $2m$. Hence~$m_{\sst\Si}(1)\ge 2m$.
\item[{\rm(b)}] The restriction of\/ $\su(m)$ moment maps
$\mu:\C^m\ra\R$ with\/ $\mu(0)=0$ to $C'$ form a vector space
of order $2$ homogeneous harmonic functions on $C'$, with
dimension $m^2\!-\!1\!-\!\dim G$. Hence~$m_{\sst\Si}(2)\ge
m^2\!-\!1\!-\!\dim G$.
\end{itemize}
\label{cm3prop1}
\end{prop}

\begin{proof} Real linear functions on $\C^m$ are moment maps
of {\it translations} on $\C^m$, and so restrict to harmonic
maps on SL $m$-folds $L$ in $\C^m$ by Lemma \ref{cm3lem2}.
Thus the vector space in (a) is of harmonic functions on $C'$,
which are clearly homogeneous of order 1. Now $C$ has a unique
singular point at 0, so it cannot be invariant under nontrivial
translations. Therefore the moment map of a nontrivial translation
cannot vanish on $C'$, and the restriction in (a) is injective.
It follows that the vector space has dimension $2m$, proving
part~(a).

For (b), each $\su(m)$ vector field has a unique moment map
$\mu:\C^m\ra\R$ with $\mu(0)=0$, which is a homogeneous real
quadratic polynomial. It follows as for (a) that the vector
space in (b) consists of order 2 homogeneous harmonic functions
on $C'$. This vector space is the image of a linear map from
$\su(m)$, and it is easy to show that the kernel of this map
is $\g$, the Lie algebra of $G$. Hence the dimension of the
vector space is $\dim\su(m)-\dim\g$ by rank-nullity, and the
proposition follows.
\end{proof}

We define the {\it stability index} of $C$, and {\it stable}
and {\it rigid} cones.

\begin{dfn} Let $C$ be an SL cone in $\C^m$ for $m>2$ with an
isolated singularity at 0, let $G$ be the Lie subgroup of\/
$\SU(m)$ preserving $C$, and use the notation of Definitions
\ref{cm3def1} and \ref{cm3def2}. Then
\e
m_{\sst\Si}(0)=b^0(\Si),\quad
m_{\sst\Si}(1)\ge 2m \quad\text{and}\quad
m_{\sst\Si}(2)\ge m^2-1-\dim G,
\label{cm3eq4}
\e
where the first equation follows as $m_{\sst\Si}(0)$ is the
multiplicity of the eigenvalue 0 of $\De_{\sst\Si}$, and the
others from Proposition \ref{cm3prop1}.

Define the {\it stability index} $\sind(C)$ to be
\e
\sind(C)=N_{\sst\Si}(2)-b^0(\Si)-m^2-2m+1+\dim G.
\label{cm3eq5}
\e
Then $\sind(C)\ge 0$ by \eq{cm3eq4}, as $N_{\sst\Si}(2)\ge
m_{\sst\Si}(0)+m_{\sst\Si}(1)+m_{\sst\Si}(2)$ by \eq{cm3eq3}.
We call $C$ {\it stable} if~$\sind(C)=0$.

Following \cite[Def.~6.7]{Joyc7}, we call $C$ {\it rigid} if
$m_{\sst\Si}(2)=m^2-1-\dim G$. As
\begin{equation*}
\sind(C)\ge m_{\sst\Si}(2)-(m^2-1-\dim G)\ge 0,
\end{equation*}
we see that {\it if\/ $C$ is stable, then $C$ is rigid}.
\label{cm3def3}
\end{dfn}

Here is the point of this definition. In deforming SL
$m$-folds $X$ in an almost Calabi--Yau $m$-fold $M$ with
a conical singularity $x$ modelled on $C$, it will turn out
in \S\ref{cm6} that $x$ contributes an {\it obstruction space}
of dimension $N_{\sst\Si}(2)$ to deforming $X$. However, we
will be able to {\it overcome} a subspace of these obstructions
with dimension $b^0(\Si)+m^2+2m-1-\dim G$ automatically, by
moving $x$ around in $M$, and changing the identification
$\C^m\cong T_xM$. Thus $\sind(C)$ is the dimension of the
{\it residual obstruction space}, which we cannot get rid of.

If $C$ is {\it stable} then the deformation problem is {\it
unobstructed}. Rigid (and more generally {\it Jacobi integrable})
SL cones were discussed in \cite[\S 6]{Joyc7}. An SL cone $C$
is {\it rigid\/} if all infinitesimal deformations of $C$ as
an SL cone come from $\su(m)$ rotations.

\subsection{Examples of special Lagrangian cones}
\label{cm32}

Examples of SL cones are constructed by Harvey and Lawson
\cite[\S III.3]{HaLa}, the author \cite{Joyc2,Joyc3}, and
others. We will study a family of special Lagrangian cones
in $\C^m$ constructed by Harvey and Lawson \cite[\S III.3.A]{HaLa}.
For $m\ge 3$, define
\e
C_{\sst\rm HL}^m=\bigl\{(z_1,\ldots,z_m)\in\C^m:\md{z_1}=
\cdots=\md{z_m}, \quad i^{m+1}z_1\cdots z_m\in[0,\iy)\bigr\}.
\label{cm3eq6}
\e

Then $C_{\sst\rm HL}^m$ is a special Lagrangian cone in $\C^m$
with an isolated singularity at 0, and $\Si_{\sst\rm HL}^m=
C_{\sst\rm HL}^m\cap{\cal S}^{2m-1}$ is an $(m\!-\!1)$-torus
$T^{m-1}$ with a flat metric. Also $C_{\sst\rm HL}^m$ and
$\Si_{\sst\rm HL}^m$ are invariant under the $\U(1)^{m-1}$
subgroup of $\SU(m)$ acting by
\e
(z_1,\ldots,z_m)\mapsto({\rm e}^{i\th_1}z_1,\ldots,{\rm e}^{i\th_m}
z_m) \quad\text{for $\th_j\in\R$ with $\th_1+\cdots+\th_m=0$.}
\label{cm3eq7}
\e
In fact $\pm\,C_{\sst\rm HL}^m$ for $m$ odd, and $C_{\sst\rm
HL}^m,iC_{\sst\rm HL}^m$ for $m$ even, are the unique SL cones
in $\C^m$ invariant under \eq{cm3eq7}, which is how Harvey and
Lawson constructed them.

We shall find the {\it stability index} $\sind(C_{\sst\rm HL}^m)$ of
these cones, and test whether they are {\it stable} or {\it rigid}.
This was first done by the author \cite[\S 3.2]{Joyc1} for $m=3$
and Marshall \cite[\S 6.3.4]{Mars} for $3\le m\le 8$. The metric on
$\Si_{\sst\rm HL}^m\cong T^{m-1}$ is flat, so it is not difficult to
compute the eigenvalues of $\De_{\smash{\sst\Si_{\rm HL}^m}}$. There
is a 1-1 correspondence between $(n_1,\ldots,n_{m-1})\in\Z^{m-1}$ and
eigenvectors of $\De_{\smash{\sst\Si_{\rm HL}^m}}$ with eigenvalue
\e
m\sum_{i=1}^{m-1}n_i^2-\sum_{i,j=1}^{m-1}n_in_j.
\label{cm3eq8}
\e

Using \eq{cm3eq8} and a computer we can find the eigenvalues of
$\De_{\smash{\sst\Si_{\rm HL}^m}}$, and their multiplicities.
Thus we can calculate $N_{\smash{\sst\Si_{\rm HL}^m}}(2)$, which
is the sum of multiplicities of eigenvalues in $[0,2m]$, and
$m_{\smash{\sst\Si_{\rm HL}^m}}(2)$, which is the multiplicity of
the eigenvalue $2m$. A table of eigenvalues and multiplicities
for $3\le m\le 8$ is given in Marshall \cite[Table 6.1]{Mars}.
Now the subgroup $G_m$ of $\SU(m)$ preserving $C_{\sst\rm HL}^m$
is $\U(1)^{m-1}$, with dimension $m-1$. Thus \eq{cm3eq5} gives
$\sind(C_{\sst\rm HL}^m)=N_{\smash{\sst\Si_{\rm HL}^m}}(2)-
m^2-m-1$. Table \ref{cm3table} gives the data
$m,N_{\smash{\sst\Si_{\rm HL}^m}}(2),
m_{\smash{\sst\Si_{\rm HL}^m}}(2)$ and
$\sind(C_{\sst\rm HL}^m)$ for~$3\le m\le 12$.

\begin{table}[htb]
\center{
\begin{tabular}{|l|r|r|r|r|r|r|r|r|r|r|}\hline
$m\vphantom{\bigr(^l_j}$
&  3 &  4 &  5 &  6 &  7 &  8 &  9 &  10 & 11 & 12 \\
\hline
$N_{\smash{\sst\Si_{\rm HL}^m}}(2)\vphantom{\bigr(^l_j}$
& 13 & 27 & 51 & 93 & 169& 311& 331& 201 & 243& 289\\
\hline
$m_{\smash{\sst\Si_{\rm HL}^m}}(2)\vphantom{\bigr(^l_j}$
&  6 & 12 & 20 & 30 & 42 & 126& 240&  90 & 110& 132\\
\hline
$\sind(C_{\sst\rm HL}^m)\vphantom{\bigr(^l_j}$
&  0 &  6 & 20 & 50 & 112& 238& 240&  90 & 110& 132\\
\hline
\end{tabular}
}
\caption{Data for $\U(1)^{m-1}$-invariant SL cones
$C_{\sst\rm HL}^m$ in $\C^m$}
\label{cm3table}
\end{table}

Motivated by Table \ref{cm3table}, with some more work one can
prove that
\e
N_{\smash{\sst\Si_{\rm HL}^m}}(2)=2m^2+1\;\>\text{and}\;\>
m_{\smash{\sst\Si_{\rm HL}^m}}(2)=\sind(C_{\sst\rm HL}^m)=m^2-m
\;\>\text{for $m\ge 10$.}
\label{cm3eq9}
\e
As $C_{\sst\rm HL}^m$ is {\it stable} when $\sind(C_{\sst\rm
HL}^m)=0$ we see from Table \ref{cm3table} and \eq{cm3eq9}
that $C_{\sst\rm HL}^3$ is a {\it stable} cone in $\C^3$, but
$C_{\sst\rm HL}^m$ is {\it unstable} for $m\ge 4$.

Also $C_{\sst\rm HL}^m$ is {\it rigid\/} when $m_{\smash{\sst\Si_{\rm
HL}^m}}(2)=m^2-m$. Thus $C_{\sst\rm HL}^m$ is {\it rigid\/} if and
only if $m\ne 8,9$, by Table \ref{cm3table} and \eq{cm3eq9}. It
would be interesting to know whether the SL cones $C_{\sst\rm HL}^8$
and $C_{\sst\rm HL}^9$ are {\it Jacobi integrable} in the sense of
\cite[\S 6]{Joyc7}, as rigid implies Jacobi integrable but not vice
versa. The author guesses that $C_{\sst\rm HL}^8,C_{\sst\rm HL}^9$
are not Jacobi integrable.

\subsection{Special Lagrangian $m$-folds with conical singularities}
\label{cm33}

Now we can define {\it conical singularities} of SL $m$-folds,
following~\cite[Def.~3.6]{Joyc7}.

\begin{dfn} Let $(M,J,\om,\Om)$ be an almost Calabi--Yau $m$-fold
for $m>2$, and define $\psi:M\ra(0,\iy)$ as in \eq{cm2eq3}. Suppose
$X$ is a compact singular SL $m$-fold in $M$ with singularities at
distinct points $x_1,\ldots,x_n\in X$, and no other singularities.

Fix isomorphisms $\up_i:\C^m\ra T_{x_i}M$ for $i=1,\ldots,n$
such that $\up_i^*(\om)=\om'$ and $\up_i^*(\Om)=\psi(x_i)^m\Om'$,
where $\om',\Om'$ are as in \eq{cm2eq1}. Let $C_1,\ldots,C_n$ be SL
cones in $\C^m$ with isolated singularities at 0. For $i=1,\ldots,n$
let $\Si_i=C_i\cap{\cal S}^{2m-1}$, and let $\mu_i\in(2,3)$ with
\e
(2,\mu_i]\cap\D_{\sst\Si_i}=\emptyset,
\quad\text{where $\D_{\sst\Si_i}$ is defined in \eq{cm3eq2}.}
\label{cm3eq10}
\e
Then we say that $X$ has a {\it conical singularity} at $x_i$,
with {\it rate} $\mu_i$ and {\it cone} $C_i$ for $i=1,\ldots,n$,
if the following holds.

By Darboux' Theorem \cite[Th.~3.15]{McSa2} there exist embeddings
$\Up_i:B_R\ra M$ for $i=1,\ldots,n$ satisfying $\Up_i(0)=x_i$,
$\d\Up_i\vert_0=\up_i$ and $\Up_i^*(\om)=\om'$, where $B_R$ is
the open ball of radius $R$ about 0 in $\C^m$ for some small $R>0$.
Define $\iota_i:\Si_i\t(0,R)\ra B_R$ by $\iota_i(\si,r)=r\si$
for~$i=1,\ldots,n$.

Define $X'=X\sm\{x_1,\ldots,x_n\}$. Then there should exist a
compact subset $K\subset X'$ such that $X'\sm K$ is a union of
open sets $S_1,\ldots,S_n$ with $S_i\subset\Up_i(B_R)$, whose
closures $\bar S_1,\ldots,\bar S_n$ are disjoint in $X$. For
$i=1,\ldots,n$ and some $R'\in(0,R]$ there should exist a smooth
$\phi_i:\Si_i\t(0,R')\ra B_R$ such that $\Up_i\circ\phi_i:\Si_i
\t(0,R')\ra M$ is a diffeomorphism $\Si_i\t(0,R')\ra S_i$, and
\e
\bmd{\na^k(\phi_i-\iota_i)}=O(r^{\mu_i-1-k})
\quad\text{as $r\ra 0$ for $k=0,1$.}
\label{cm3eq11}
\e
Here $\na,\md{\,.\,}$ are computed using the cone metric
$\iota_i^*(g')$ on~$\Si_i\t(0,R')$.

If the cones $C_1,\ldots,C_n$ are {\it stable} in the sense of
Definition \ref{cm3def3}, then we say that $X$ has {\it stable
conical singularities}.
\label{cm3def4}
\end{dfn}

The reasoning behind this definition was discussed in
\cite[\S 3.3]{Joyc7}. Here we just make two remarks:
\begin{itemize}
\item We suppose $m>2$ for two reasons. Firstly, the only
SL cones in $\C^2$ are finite unions of SL planes $\R^2$
in $\C^2$ intersecting only at 0. Thus any SL 2-fold with
conical singularities is actually {\it nonsingular} as an
immersed 2-fold, so there is really no point in studying
them. Secondly, $m=2$ is a special case in the analysis
of \cite[\S 2]{Joyc7}, and it is simpler to exclude it.

In the rest of the paper we shall assume~$m>2$.
\item The purpose of \eq{cm3eq10} is to reduce to a minimum
the obstructions to deforming $X$ as an SL $m$-fold with
conical singularities. If we omitted condition \eq{cm3eq10} then
each $\al\in(2,\mu_i]\cap\D_{\sst\Si_i}$ would contribute
additional obstructions to deforming $X$ in~\S\ref{cm6}.
\end{itemize}

\section{Review of material from \cite{Joyc7}}
\label{cm4}

We now review the definitions and results from the preceding
paper \cite{Joyc7} which we will need later. Throughout we
suppose~$m>2$.

\subsection{Analysis on SL $m$-folds with conical singularities}
\label{cm41}

We will need the following tool \cite[Def.~2.6]{Joyc7}, a
smoothed out version of the distance from the singular set
$\{x_1,\ldots,x_n\}$ in~$X$.

\begin{dfn} Let $(M,J,\om,\Om)$ be an almost Calabi--Yau
$m$-fold and $X$ a compact SL $m$-fold in $M$ with conical
singularities at $x_1,\ldots,x_n$, and use the notation of
Definition \ref{cm3def4}. Define a {\it radius function}
$\rho$ on $X'$ to be a smooth function $\rho:X'\ra(0,1]$
such that $\rho\equiv 1$ on $K$ and $\rho(y)=d(x_i,y)$ for
$y\in S_i$ close to $x_i$, where $d$ is the metric on $X$.
Radius functions always exist.

For $\bs{\be}=(\be_1,\ldots,\be_n)\in\R^n$, define a function
$\rho^{\bs\be}$ on $X'$ by $\rho^{\bs\be}(y)=\rho(y)^{\be_i}$
on $S_i$ for $i=1,\ldots,n$ and $\rho^{\bs\be}(y)=1$ on $K$.
Then $\rho^{\bs\be}$ is well-defined and smooth on $X'$, and
equals $\rho^{\be_i}$ near $x_i$ in $X'$. If ${\bs\be},{\bs\ga}
\in\R^n$, write ${\bs\be}\ge{\bs\ga}$ if $\be_i\ge\ga_i$ for
$i=1,\ldots,n$. If ${\bs\be}\in\R^n$ and $a\in\R$, write
${\bs\be}+a=(\be_1+a,\ldots,\be_n+a)$ in~$\R^n$.
\label{cm4def1}
\end{dfn}

Now we define some Banach spaces of functions on $X'$,
\cite[Def.~2.7]{Joyc7}.

\begin{dfn} Let $(M,J,\om,\Om)$ be an almost Calabi--Yau
$m$-fold with metric $g$, and $X$ a compact SL $m$-fold in $M$
with conical singularities at $x_1,\ldots,x_n$, and use the
notation of Definitions \ref{cm3def4} and \ref{cm4def1}. Let
$\rho$ be a radius function on $X'$. Regard $X'$ as a Riemannian
manifold, with metric $g$ restricted from~$M$.

For ${\bs\be}\in\R^n$ and $k\ge 0$ define $C^k_{\smash{\bs\be}}(X')$
to be the space of continuous functions $f$ on $X'$ with $k$
continuous derivatives, such that $\bmd{\rho^{-{\bs\be}+j}\na^jf}$
is bounded on $X'$ for $j=0,\dots,k$. Define the norm $\nm{\,.\,}_{
\smash{C^k_{\bs\be}}}$ on $C^k_{\smash{\bs\be}}(X')$ by
\e
\nm{f}_{\smash{C^k_{\bs\be}}}=\sum_{j=0}^k\sup_{X'}
\bmd{\rho^{-{\bs\be}+j}\na^jf}.
\label{cm4eq1}
\e
Then $C^k_{\smash{\bs\be}}(X')$ is a Banach space.
Define~$C^\iy_{\smash{\bs\be}}(X')=\bigcap_{k\ge 0}
C^k_{\smash{\bs\be}}(X')$.

For $p\ge 1$, ${\bs\be}\in\R^n$ and $k\ge 0$ define the
{\it weighted Sobolev space} $L^p_{k,{\bs\be}}(X')$ to be
the set of functions $f$ on $X'$ that are locally integrable
and $k$ times weakly differentiable, and for which the norm
\e
\snm{f}p{k,{\bs\be}}=
\left(\sum_{j=0}^k\int_{X'}\bmd{\rho^{-{\bs\be}+j}\na^jf}^p
\rho^{-m}\d V_g\right)^{1/p}
\label{cm4eq2}
\e
is finite. Then $L^p_{k,{\bs\be}}(X')$ is a Banach space, and
$L^2_{k,{\bs\be}}(X')$ a Hilbert space. 
\label{cm4def2}
\end{dfn}

We call these {\it weighted Banach spaces} since the norms are
locally weighted by a power of $\rho$. Roughly speaking, if $f$
lies in $L^p_{k,{\bs\be}}(X')$ or $C^k_{\smash{\bs\be}}(X')$ then
$f$ grows at most like $\rho^{\be_i}$ near $x_i$ as $\rho\ra 0$,
and so the multi-index ${\bs\be}=(\be_1,\ldots,\be_n)$ should be
interpreted as an {\it order of growth}.

Here is a weighted version of the {\it Sobolev Embedding
Theorem},~\cite[Th.~2.9]{Joyc7}.

\begin{thm} In the situation above, suppose $k>l\ge 0$ are
integers and\/ $p>1$ with\/ $\frac{1}{p}<\frac{k-l}{m}$, and\/
${\bs\be},{\bs\ga}\in\R^n$ with\/ ${\bs\be}\ge{\bs\ga}$.
Then $L^p_{k,{\bs\be}}(X')\hookra C^l_{\bs\ga}(X')$ is a
continuous inclusion.
\label{cm4thm1}
\end{thm}

Here is a Fredholm result for the operator $P:f\mapsto\d^*(\psi^m\d f)$
on weighted Sobolev spaces, \cite[Th.~5.3]{Joyc7}. Putting $\al=\d f$
in \eq{cm2eq4}, we see that $P$ appears in the linearization of the
deformation problem for SL $m$-folds.

\begin{thm} Let\/ $(M,J,\om,\Om)$ be an almost Calabi--Yau
$m$-fold, and define $\psi:M\ra(0,\iy)$ as in \eq{cm2eq3}.
Suppose $X$ is a compact SL\/ $m$-fold in $M$ with conical
singularities at $x_1,\ldots,x_n$ with cones $C_i$. Define
$\D_{\sst\Si_i},N_{\sst\Si_i}$ and\/ $L^p_{k,{\bs\be}}(X')$
as in Definitions \ref{cm3def2} and \ref{cm4def2}. Fix $p>1$
and\/ $k\ge 2$, and for ${\bs\be}\in\R^n$ define $P_{\bs\be}:
L^p_{k,{\bs\be}}(X')\ra L^p_{k-2,{\bs\be}-2}(X')$ by
$P_{\bs\be}(f)=\d^*(\psi^m\d f)$. Then
\begin{itemize}
\item[{\rm(a)}] $P_{\bs\be}$ is Fredholm if and only if\/
${\bs\be}\in\bigl(\R\sm\D_{\sst\Si_1}\bigr)\t\cdots\t
\bigl(\R\sm\D_{\sst\Si_n}\bigr)$, and then
\e
\ind(P_{\bs\be})=-\sum_{i=1}^nN_{\sst\Si_i}(\be_i).
\label{cm4eq3}
\e
\item[{\rm(b)}] If\/ $\be_i>0$ for all\/ $i$ then $P_{\bs\be}$
is injective.
\end{itemize}
\label{cm4thm2}
\end{thm}

\subsection{Homology, cohomology and Hodge theory}
\label{cm42}

Next we discuss {\it homology} and {\it cohomology} of SL
$m$-folds with conical singularities, following \cite[\S 2.4]{Joyc7}.
For a general reference, see for instance Bredon \cite{Bred}.
When $Y$ is a manifold, write $H^k(Y,\R)$ for the $k^{\rm th}$
{\it de Rham cohomology group} and $H^k_{\rm cs}(Y,\R)$ for the
$k^{\rm th}$ {\it compactly-supported de Rham cohomology group}
of $Y$. If $Y$ is compact then~$H^k(Y,\R)=H^k_{\rm cs}(Y,\R)$.

Let $Y$ be a topological space, and $Z\subset Y$ a subspace.
Write $H_k(Y,\R)$ for the $k^{\rm th}$ {\it real singular
homology group} of $Y$, and $H_k(Y;Z,\R)$ for the $k^{\rm th}$
{\it real singular relative homology group} of $(Y;Z)$. When
$Y$ is a manifold and $Z$ a submanifold we define $H_k(Y,\R)$
and $H_k(Y;Z,\R)$ using {\it smooth\/} simplices, as in
\cite[\S V.5]{Bred}. Then the pairing between (singular)
homology and (de Rham) cohomology is defined at the chain
level by integrating $k$-forms over $k$-simplices.

Suppose $X$ is a compact SL $m$-fold in $M$ with conical
singularities $x_1,\ldots,x_n$ and cones $C_1,\ldots,C_n$, and set
$X'=X\sm\{x_1,\ldots,x_n\}$ and $\Si_i=C_i\cap{\cal S}^{2m-1}$ as
above. Then by \cite[\S 2.4]{Joyc7} there is a natural long exact
sequence
\e
\cdots\ra
H^k_{\rm cs}(X',\R)\ra H^k(X',\R)\ra\bigoplus_{i=1}^n
H^k(\Si_i,\R)\ra H^{k+1}_{\rm cs}(X',\R)\ra\cdots,
\label{cm4eq4}
\e
and natural isomorphisms
\begin{gather}
H_k\bigl(X;\{x_1,\ldots,x_n\},\R\bigr)^*\!\cong\!
H^k_{\rm cs}(X',\R)\!\cong\!H_{m-k}(X',\R)\!\cong\!H^{m-k}(X',\R)^*
\label{cm4eq5}\\
\text{and}\quad
H^k_{\rm cs}(X',\R)\cong H_k(X,\R)^*
\quad\text{for all $k>1$.}
\label{cm4eq6}
\end{gather}
The inclusion $\iota:X\ra M$ induces homomorphisms $\iota_*:
H_k(X,\R)\ra H_k(M,\R)$ and~$\iota^*:H^k(M,\R)\ra H^k(X',\R)$.

If $(Y,g)$ is a compact Riemannian manifold, then {\it Hodge theory}
shows that each class in $H^k(Y,\R)$ is represented by a unique
$k$-form $\al$ with $\d\al=\d^*\al=0$. Here is an analogue of this
on $X'$ when $k=1$, part of~\cite[Th.~5.4]{Joyc7}.

\begin{thm} Let\/ $(M,J,\om,\Om)$ be an almost Calabi--Yau
$m$-fold, and define $\psi:M\ra(0,\iy)$ as in \eq{cm2eq3}. Suppose
$X$ is a compact SL\/ $m$-fold in $M$ with conical singularities
at $x_1,\ldots,x_n$. Set\/ $X'=X\sm\{x_1,\ldots,x_n\}$, and let\/
$\rho$ be a radius function on $X'$, in the sense of Definition
\ref{cm4def1}. Define
\e
\begin{split}
Y_\sXp=\bigl\{\al\in C^\iy(T^*X'):\,&\d\al=0,\quad
\d^*(\psi^m\al)=0,\\
&\text{$\md{\na^k\al}=O(\rho^{-1-k})$ for $k\ge 0$}\bigr\}.
\end{split}
\label{cm4eq7}
\e
Then the map $\pi:Y_\sXp\ra H^1(X',\R)$ taking
$\pi:\al\mapsto[\al]$ is an isomorphism.
\label{cm4thm3}
\end{thm}

\subsection{Lagrangian Neighbourhood Theorems}
\label{cm43}

In \cite[\S 4]{Joyc7} we extend the {\it Lagrangian
Neighbourhood Theorem}, Theorem \ref{cm2thm2}, to situations
involving conical singularities, first to {\it SL
cones},~\cite[Th.~4.3]{Joyc7}.

\begin{thm} Let\/ $C$ be an SL cone in $\C^m$ with isolated
singularity at\/ $0$, and set\/ $\Si=C\cap{\cal S}^{2m-1}$.
Define $\iota:\Si\t(0,\iy)\ra\C^m$ by $\iota(\si,r)=r\si$,
with image $C\sm\{0\}$. For $\si\in\Si$, $\tau\in T_\si^*\Si$,
$r\in(0,\iy)$ and\/ $u\in\R$, let\/ $(\si,r,\tau,u)$ represent
the point\/ $\tau+u\,\d r$ in $T^*_{\smash{(\si,r)}}\bigl(\Si\!
\t\!(0,\iy)\bigr)$. Identify $\Si\t(0,\iy)$ with the zero section
$\tau\!=\!u\!=\!0$ in $T^*\bigl(\Si\t(0,\iy)\bigr)$. Define an
action of\/ $(0,\iy)$ on $T^*\bigl(\Si\!\t\!(0,\iy)\bigr)$ by
\e
t:(\si,r,\tau,u)\longmapsto (\si,tr,t^2\tau,tu)
\quad\text{for $t\in(0,\iy)$,}
\label{cm4eq8}
\e
so that\/ $t^*(\hat\om)\!=\!t^2\hat\om$, for $\hat\om$ the
canonical symplectic structure on~$T^*\bigl(\Si\!\t\!(0,\iy)\bigr)$.

Then there exists an open neighbourhood\/ $U_{\sst C}$ of\/
$\Si\t(0,\iy)$ in $T^*\bigl(\Si\t(0,\iy)\bigr)$ invariant under
\eq{cm4eq8} given by
\e
U_{\sst C}=\bigl\{(\si,r,\tau,u)\in T^*\bigl(\Si\t(0,\iy)\bigr):
\bmd{(\tau,u)}<2\ze r\bigr\}\quad\text{for some $\ze>0$,}
\label{cm4eq9}
\e
where $\md{\,.\,}$ is calculated using the cone metric $\iota^*(g')$
on $\Si\t(0,\iy)$, and an embedding $\Phi_{\sst C}:U_{\sst C}\ra\C^m$
with\/ $\Phi_{\sst C}\vert_{\Si\t(0,\iy)}=\iota$, $\Phi_{\sst
C}^*(\om')=\hat\om$ and\/ $\Phi_{\sst C}\circ t=t\,\Phi_{\sst C}$
for all\/ $t>0$, where $t$ acts on $U_{\sst C}$ as in \eq{cm4eq8}
and on $\C^m$ by multiplication.
\label{cm4thm4}
\end{thm}

In \cite[Th.~4.4]{Joyc7} we construct a particular choice of
$\phi_i$ in Definition~\ref{cm3def4}.

\begin{thm} Let\/ $(M,J,\om,\Om)$, $\psi,X,n,x_i,\up_i,C_i,\Si_i,
\mu_i,R,\Up_i$ and\/ $\iota_i$ be as in Definition \ref{cm3def4}.
Theorem \ref{cm4thm4} gives $\ze>0$, neighbourhoods $U_{\sst C_i}$
of\/ $\Si_i\t(0,\iy)$ in $T^*\bigl(\Si_i\t(0,\iy)\bigr)$ and
embeddings $\Phi_{\sst C_i}:U_{\sst C_i}\ra\C^m$ for~$i=1,\ldots,n$.

Then for sufficiently small\/ $R'\in(0,R]$ there exist unique
closed\/ $1$-forms $\eta_i$ on $\Si_i\t(0,R')$ for $i=1,\ldots,n$
written $\eta_i(\si,r)=\eta_i^1(\si,r)+\eta_i^2(\si,r)\d r$ for
$\eta_i^1(\si,r)\in T_\si^*\Si_i$ and\/ $\eta_i^2(\si,r)\in\R$,
and satisfying $\md{\eta_i(\si,r)}<\ze r$ and
\e
\bmd{\na^k\eta_i}=O(r^{\mu_i-1-k})
\quad\text{as $r\ra 0$ for $k=0,1$,}
\label{cm4eq10}
\e
computing $\na,\md{\,.\,}$ using the cone metric $\iota_i^*(g')$,
such that the following holds.

Define $\phi_i:\Si_i\t(0,R')\ra B_R$ by $\phi_i(\si,r)=\Phi_{\sst
C_i}\bigl(\si,r,\eta_i^1(\si,r),\eta_i^2(\si,r)\bigr)$. Then
$\Up_i\circ\phi_i:\Si_i\t(0,R')\ra M$ is a diffeomorphism
$\Si_i\t(0,R')\ra S_i$ for open sets $S_1,\ldots,S_n$ in $X'$
with\/ $\bar S_1,\ldots,\bar S_n$ disjoint, and\/ $K=X'\sm(S_1
\cup\cdots\cup S_n)$ is compact. Also $\phi_i$ satisfies
\eq{cm3eq11}, so that\/ $R',\phi_i,S_i,K$ satisfy
Definition~\ref{cm3def4}.
\label{cm4thm5}
\end{thm}

Next we extend Theorem \ref{cm2thm2} to SL $m$-folds with
conical singularities \cite[Th.~4.6]{Joyc7}, in a way
compatible with Theorems \ref{cm4thm4} and~\ref{cm4thm5}.

\begin{thm} Suppose $(M,J,\om,\Om)$ is an almost Calabi--Yau
$m$-fold and\/ $X$ a compact SL\/ $m$-fold in $M$ with conical
singularities at\/ $x_1,\ldots,x_n$. Let the notation $\psi,\up_i,
C_i,\Si_i,\mu_i,R,\Up_i$ and\/ $\iota_i$ be as in Definition
\ref{cm3def4}, and let\/ $\ze,U_{\sst C_i},\allowbreak
\Phi_{\sst C_i},\allowbreak R',\allowbreak \eta_i,\allowbreak
\eta_i^1,\eta_i^2,\phi_i,S_i$ and\/ $K$ be as in Theorem~\ref{cm4thm5}.

Then making $R'$ smaller if necessary, there exists an open tubular
neighbourhood\/ $U_\sXp\subset T^*X'$ of the zero section
$X'$ in $T^*X'$, such that under $\d(\Up_i\circ\phi_i):T^*\bigl(
\Si_i\t(0,R')\bigr)\ra T^*X'$ for $i=1,\ldots,n$ we have
\e
\bigl(\d(\Up_i\circ\phi_i)\bigr)^*(U_\sXp)=\bigl\{(\si,r,\tau,u)
\in T^*\bigl(\Si_i\t(0,R')\bigr):\bmd{(\tau,u)}<\ze r\bigr\},
\label{cm4eq11}
\e
and there exists an embedding $\Phi_\sXp:U_\sXp\ra M$ with\/
$\Phi_\sXp\vert_{X'}=\id:X'\ra X'$ and\/ $\Phi_\sXp^*(\om)
=\hat\om$, where $\hat\om$ is the canonical symplectic structure on
$T^*X'$, such that 
\e
\Phi_\sXp\circ\d(\Up_i\circ\phi_i)(\si,r,\tau,u)\equiv\Up_i\circ
\Phi_{\sst C_i}\bigl(\si,r,\tau+\eta_i^1(\si,r),u+\eta_i^2(\si,r)\bigr)
\label{cm4eq12}
\e
for all\/ $i=1,\ldots,n$ and\/ $(\si,r,\tau,u)\in T^*\bigl(\Si_i\t(0,R')
\bigr)$ with\/ $\bmd{(\tau,u)}<\ze r$. Here $\md{(\tau,u)}$ is computed
using the cone metric $\iota_i^*(g')$ on~$\Si_i\t(0,R')$.
\label{cm4thm6}
\end{thm}

Here is an extension of Theorem \ref{cm4thm6} to {\it families}
of almost Calabi--Yau $m$-folds $(M,J^s,\om^s,\Om^s)$ for $s\in\F$,
deduced from \cite[Th.~4.8 \& Th.~4.9]{Joyc7}.

\begin{thm} Let\/ $(M,J,\om,\Om)$ be an almost Calabi--Yau
$m$-fold and\/ $X$ a compact SL\/ $m$-fold in $M$ with conical
singularities at $x_1,\ldots,x_n$, with identifications $\up_i$
and cones $C_i$. Let the notation $R,\Up_i,\ze,\Phi_{\sst C_i},
R',\eta_i,\eta_i^1,\eta_i^2,\phi_i,S_i,K$ be as in Theorem
\ref{cm4thm5}, and let\/ $U_\sXp,\Phi_\sXp$ be as in
Theorem~\ref{cm4thm6}.

Suppose $\bigl\{(M,J^s,\om^s,\Om^s):s\in\F\bigr\}$ is
a smooth family of deformations of $(M,J,\om,\Om)$, in the
sense of Definition \ref{cm2def6}, such that\/ $\iota_*(\ga)
\cdot[\om^s]=0$ for all\/ $\ga\in H_2(X,\R)$ and\/ $s\in\F$,
where $\iota:X\ra M$ is the inclusion and\/ $\iota_*:H_2(X,\R)\ra
H_2(M,\R)$ the induced homomorphism. Define $\psi^s:M\ra(0,\iy)$
for $s\in\F$ as in \eq{cm2eq3}, but using~$\om^s,\Om^s$.

Then making $R,R'$ and\/ $U_\sXp$ smaller if necessary,
for some connected open $\F'\subseteq\F$ with\/
$0\in\F'$ and all\/ $s\in \F'$ there exist
\begin{itemize}
\item[{\rm(a)}] isomorphisms $\up_i^s:\C^m\ra T_{x_i}M$ for
$i=1,\ldots,n$ with\/ $\up_i^0=\up_i$, $(\up_i^s)^*(\om^s)=\om'$
and\/~$(\up_i^s)^*(\Om)=\psi^s(x_i)^m\Om'$,
\item[{\rm(b)}] embeddings $\Up_i^s:B_R\ra M$ for $i=1,\ldots,n$
with\/ $\Up_i^0=\Up_i$, $\Up_i^s(0)=x_i$, $\d\Up_i^s\vert_0=\up_i^s$,
$(\Up_i^s)^*(\om^s)=\om'$, and
\item[{\rm(c)}] an embedding $\Phi^s_\sXp:U_\sXp\ra M$ with\/
$\Phi_\sXp^0=\Phi_\sXp$ and\/~$(\Phi^s_\sXp)^*(\om^s)=\hat\om$,
\end{itemize}
all depending smoothly on $s\in\F'$ with
\e
\Phi_\sXp^s\circ\d(\Up_i\circ\phi_i)(\si,r,\tau,u)\equiv\Up_i^s\circ
\Phi_{\sst C_i}\bigl(\si,r,\tau+\eta_i^1(\si,r),u+\eta_i^2(\si,r)\bigr)
\label{cm4eq13}
\e
for all\/ $s\in \F'$, $i=1,\ldots,n$ and\/ $(\si,r,\tau,u)\in
T^*\bigl(\Si_i\t(0,R')\bigr)$ with\/~$\bmd{(\tau,u)}<\ze r$.
\label{cm4thm7}
\end{thm}

The condition that $\iota_*(\ga)\cdot[\om^s]=0$ for all $\ga\in
H_2(X,\R)$ essentially says that $\iota^*\bigl([\om^s]\bigr)=0$
in $H^2(X,\R)$. However, we have not put it like this as we have
not defined de Rham cohomology on the singular manifold $X$. We
could make sense of this by, for instance, interpreting $[\om^s]$
as a \v Cech cohomology class on $M$ using the equivalence of
de Rham and \v Cech cohomology, and pulling back to the \v Cech
cohomology of~$X$.

\subsection{Regularity of $X$ near $x_i$}
\label{cm44}

In \cite[\S 5]{Joyc7} we study the asymptotic behaviour of the maps
$\phi_i$ of Theorem \ref{cm4thm5}, using the elliptic regularity of
the special Lagrangian condition. Combining \cite[Th.~5.1]{Joyc7},
\cite[Lem.~4.5]{Joyc7} and \cite[Th.~5.5]{Joyc7} proves:

\begin{thm} In the situation of Theorem \ref{cm4thm5} we have
$\eta_i=\d A_i$ for $i=1,\ldots,n$, where $A_i:\Si_i\t(0,R')\ra\R$
is given by $A_i(\si,r)=\int_0^r\eta_i^2(\si,s)\d s$. Suppose
$\mu_i'\in(2,3)$ with\/ $(2,\mu_i']\cap\D_{\sst\Si_i}=
\emptyset$ for $i=1,\ldots,n$. Then
\e
\begin{gathered}
\bmd{\na^k(\phi_i-\iota_i)}=O(r^{\mu_i'-1-k}),\quad
\bmd{\na^k\eta_i}=O(r^{\mu_i'-1-k})\quad\text{and}\\
\bmd{\na^kA_i}=O(r^{\mu_i'-k})
\quad\text{as $r\ra 0$ for all\/ $k\ge 0$ and\/ $i=1,\ldots,n$.}
\end{gathered}
\label{cm4eq14}
\e

Hence $X$ has conical singularities at $x_i$ with cone $C_i$
and rate $\mu_i'$, for all possible rates $\mu_i'$ allowed by
Definition \ref{cm3def4}. Therefore, the definition of
conical singularities is essentially independent of the
choice of rate~$\mu_i$.
\label{cm4thm8}
\end{thm}

Theorem \ref{cm4thm8} in effect {\it strengthens} the definition of
SL $m$-folds with conical singularities, Definition \ref{cm3def4},
as it shows that \eq{cm3eq11} actually implies the much stronger
condition \eq{cm4eq14} on all derivatives. In \cite[Th.~6.8]{Joyc7}
we use {\it Geometric Measure Theory} to prove a {\it weakening} of
Definition \ref{cm3def4} for {\it rigid\/} cones~$C$.

\begin{thm} Let\/ $(M,J,\om,\Om)$ be an almost Calabi--Yau
$m$-fold and define $\psi:M\ra(0,\iy)$ as in \eq{cm2eq3}.
Let\/ $x\in M$ and fix an isomorphism $\up:\C^m\ra T_xM$
with\/ $\up^*(\om)=\om'$ and\/ $\up^*(\Om)=\psi(x)^m\Om'$,
where $\om',\Om'$ are as in~\eq{cm2eq1}.

Suppose that\/ $T$ is a special Lagrangian integral current
in $M$ with\/ $x\in T^\circ$, and that\/ $\up_*(C)$ is a
multiplicity $1$ tangent cone to $T$ at\/ $x$, where $C$ is
a rigid special Lagrangian cone in $\C^m$ in the sense of
Definition \ref{cm3def3}. Then $T$ has a conical singularity
at\/ $x$, in the sense of Definition~\ref{cm3def4}.
\label{cm4thm9}
\end{thm}

Here {\it integral currents}, {\it tangent cones} and
{\it multiplicity} are technical terms from Geometric
Measure Theory which are explained in \cite[\S 6]{Joyc7}.
In fact \cite[Th.~6.8]{Joyc7} applies to the larger class
of {\it Jacobi integrable} SL cones $C$, for which all
special Lagrangian Jacobi fields are integrable.

Basically, Theorem \ref{cm4thm9} shows that if a singular
SL $m$-fold $T$ in $M$ is locally modelled on a rigid SL cone
$C$ in only a very weak sense, then it necessarily satisfies
Definition \ref{cm3def4}. One moral of Theorems \ref{cm4thm8}
and \ref{cm4thm9} is that, at least for rigid SL cones $C$,
more-or-less any sensible definition of SL $m$-folds with
conical singularities is equivalent to Definition~\ref{cm3def4}.

\section{Moduli of SL $m$-folds with conical singularities}
\label{cm5}

The rest of the paper studies {\it moduli spaces} $\M_\sX$
of compact SL $m$-folds $X$ with conical singularities in an almost
Calabi--Yau manifold $M$. This section sets up the notation needed
to do this, and defines the moduli space $\M_\sX$ as a
topological space, paying particular attention to the r\^ole of
asymptotic conditions at the singular points in defining the
topology on $\M_\sX$. We continue to suppose~$m>2$.

\subsection{Notation to vary the $x_i,\up_i$}
\label{cm51}

We are interested in deformations of $X$ in $M$ that are
allowed to move the singular points $x_1,\ldots,x_n$ and
the identifications $\up_i:\C^m\ra T_{x_i}M$. We begin by
setting up some notation to allow us to do this.

\begin{dfn} Let $(M,J,\om,\Om)$ be an almost Calabi--Yau
$m$-fold and $X$ a compact SL $m$-fold in $M$ with conical
singularities at $x_1,\ldots,x_n$ with identifications
$\up_i:\C^m\ra T_{x_i}M$ and cones $C_1,\ldots,C_n$, and
use the notation of \S\ref{cm33}. Define
\e
\begin{split}
P=\bigl\{(x,\up):\,&\text{$x\in M$, $\up:\C^m\ra T_xM$ is
a real isomorphism,}\\
&\up^*(\om)=\om',\quad \up^*(\Om)=\psi(x)^m\Om'\bigr\},
\end{split}
\label{cm5eq1}
\e
where $\om',\Om'$ are as in \eq{cm2eq1}. Then $(x_i,\up_i)\in P$
for $i=1,\ldots,n$, and $P$ is the family of all possible alternative
choices of $x_i,\up_i$, by Definition~\ref{cm3def4}.

Regard each matrix $B\in\SU(m)$ as a map $\C^m\ra\C^m$. Then
if $(x,\up)\in P$ and $B\in\SU(m)$ then $(x,\up\circ B)\in P$
as $\om',\Om'$ are $\SU(m)$-invariant. Define a smooth, free
action of $\SU(m)$ on $P$ by $B:(x,\up)\mapsto(x,\up\circ B^{-1})$.
If $(x,\up),(x,\hat\up)\in P$ then $B=\hat\up^{-1}\circ\up\in\SU(m)$
and $B(x,\up)=(x,\hat\up)$. Hence the $\SU(m)$-orbits in $P$
correspond to points $x\in M$, and $P$ is a principal $\SU(m)$-bundle
over $M$. Thus~$\dim P=m^2+2m-1$.

Let $G_i$ be the Lie subgroup of $\SU(m)$ preserving the cone
$C_i$ in $\C^m$ for $i=1,\ldots,n$. Then $G_i$ acts on $P$. If
$(x,\up)$ and $(x,\hat\up)$ lie in the same $G_i$-orbit then they
define {\it equivalent\/} alternative choices for $(x_i,\up_i)$,
since $\up(C_i)$ and $\hat\up(C_i)$ are the same SL cone in $T_xM$.
Therefore if we use $P$ to parametrize alternative choices for
$(x_i,\up_i)$ we will have redundant parameters when $\dim G_i>0$,
since each cone $\up(C_i)$ in $T_xM$ is represented not by a point
in $P$ but by a submanifold isomorphic to~$G_i$.

To avoid this, let $\E_i$ be a small open ball of
dimension $\dim P-\dim G_i$ in $P$ containing $(x_i,\up_i)$
and transverse to the orbits of $G_i$ for $i=1,\ldots,n$.
Then $G_i\cdot\E_i$ is a small open neighbourhood of
the $G_i$-orbit of $(x_i,\up_i)$ in $P$. Define $\E=
\E_1\t\cdots\t\E_n$ and $e=(x_1,\up_1,\ldots,
x_n,\up_n)\in\E$. Write a general element of $\E$
as $\hat e=(\hat x_1,\hat\up_1,\ldots,\hat x_n,\hat\up_n)$.
Then $\E$ is a family of alternative choices $\hat x_i,
\hat\up_i$ of the $x_i,\up_i$, which represent all nearby
alternative choices exactly once up to equivalence, and
\e
\begin{gathered}
\dim\E_i=m^2+2m-1-\dim G_i\\
\text{and}\quad
\dim\E=n(m^2+2m-1)-\ts\sum_{i=1}^n\dim G_i.
\end{gathered}
\label{cm5eq2}
\e
The metric $g$ on $M$ induces a Riemannian metric on $P$
which restricts to $\E_i$. Let $d_\sE$ be the metric
induced on $\E=\E_1\t\cdots\t\E_n$ by the product
Riemannian metric, so that $(\E,d_\sE)$ is a metric space.
\label{cm5def1}
\end{dfn}

The following result, modelled loosely on Theorem \ref{cm4thm7},
extends $X$ to a family of Lagrangian $m$-folds $\hat X$ with
conical singularities at $\hat x_i$ and identifications $\hat\up_i$
for $\hat e=(\hat x_1,\hat\up_1,\ldots,\hat x_n,\hat\up_n)$
in an open neighbourhood $\ti\E$ of $e$ in $\E$, and
also defines Lagrangian neighbourhoods $\Phi_\sXp^{\hat e}$
for~$\hat X$.

\begin{thm} Suppose $(M,J,\om,\Om)$ is an almost Calabi--Yau
$m$-fold and\/ $X$ a compact SL\/ $m$-fold in $M$ with conical
singularities at $x_1,\ldots,x_n$. Use the notation of Theorem
\ref{cm4thm5}, let\/ $U_\sXp,\Phi_\sXp$ be as in
Theorem \ref{cm4thm6}, and\/ $e,\E$ as in Definition
\ref{cm5def1}. Then for some connected open $\ti\E\subseteq
\E$ with\/ $e\in\ti\E$ and all\/ $\hat e=(\hat x_1,
\hat\up_1,\ldots,\hat x_n,\hat\up_n)$ in $\ti\E$ there exist
\begin{itemize}
\item[{\rm(a)}] embeddings $\Up_i^{\hat e}:B_R\ra M$ for
$i=1,\ldots,n$ with
\e
\Up_i^e=\Up_i,\quad
(\Up_i^{\hat e})^*(\om)=\om',\quad
\Up_i^{\hat e}(0)=\hat x_i\quad\text{and}\quad
\d\Up_i^{\hat e}\vert_0=\hat\up_i,
\label{cm5eq3}
\e
\item[{\rm(b)}] an embedding $\Phi^{\hat e}_\sXp:U_\sXp
\ra M$ with\/ $\Phi_\sXp^e=\Phi_\sXp$ and\/
$(\Phi^{\hat e}_\sXp)^*(\om)=\hat\om$, such that\/
$\Phi_\sXp^{\hat e}\equiv\Phi_\sXp$ on~$\pi^*(K)\subset
U_\sXp$,
\end{itemize}
all depending smoothly on $\hat e\in\ti\E$, with
\e
\Phi_\sXp^{\hat e}\circ\d(\Up_i\circ\phi_i)(\si,r,\tau,u)\equiv
\Up_i^{\hat e}\circ\Phi_{\sst C_i}\bigl(\si,r,\tau+\eta_i^1(\si,r),
u+\eta_i^2(\si,r)\bigr)
\label{cm5eq4}
\e
for all\/ $\hat e\in\ti\E$, $i=1,\ldots,n$ and\/ $(\si,r,\tau,u)
\in T^*\bigl(\Si_i\t(0,R')\bigr)$ with\/~$\bmd{(\tau,u)}<\ze r$.
\label{cm5thm1}
\end{thm}

\begin{proof} We shall define $\Up_i^{\hat e}$ and $\Phi^{\hat e}_{
\sst X'}$ by modifying $\Up_i,\Phi_\sXp$ near $x_i\in M$ using
a symplectomorphism of $B_R\subset\C^m$. Let $R''\in(0,\ha R)$
satisfy conditions we will specify at the end of the proof, and
let $B_{R''},B_{2R''}\subset B_R$ be the open balls of radius
$R'',2R''$ about 0 in $\C^m$. Choose a connected open
neighbourhood $\ti\E$ of $e$ in $\E$ such that for all
$\hat e=(\hat x_1,\hat\up_1,\ldots,\hat x_n,\hat\up_n)$ in
$\ti\E$ we have $\hat x_i\in\Up_i(B_{R''})$ for
$i=1,\ldots,n$. Clearly this is possible.

Next, choose diffeomorphisms $\Xi_i^{\hat e}:B_R\ra B_R$ for
$i=1,\ldots,n$ and $\hat e\in\ti\E$ depending smoothly on
$\hat e$, such that
\begin{itemize}
\item[(i)] $\Xi_i^e$ is the identity on $B_R$ for~$i=1,\ldots,n$,
\item[(ii)] $(\Xi_i^{\hat e})^*(\om')=\om'$ for
$\hat e\in\ti\E$ and~$i=1,\ldots,n$,
\item[(iii)] $\Up_i\circ\Xi_i^{\hat e}(0)=\hat x_i$ and
$\d(\Up_i\circ\Xi_i^{\hat e})\vert_0=\hat\up_i$ for
$\hat e=(\hat x_1,\hat\up_1,\ldots,\hat x_n,\hat\up_n)
\in\ti\E$ and $i=1,\ldots,n$, and
\item[(iv)] $\Xi_i^{\hat e}$ is the identity outside $B_{2R''}
\subset B_R$ for $\hat e\in\ti\E$ and~$i=1,\ldots,n$.
\end{itemize}
Making $\ti\E$ smaller if necessary, one can do this
explicitly using standard but messy symplectic geometry
techniques, and we leave it as an exercise.

Now define an embedding $\Up_i^{\hat e}=\Up_i\circ\Xi_i^{\hat e}
:B_R\ra M$ for $i=1,\ldots,n$ and $\hat e\in\ti\E$. Then
$\Up_i^{\hat e}$ depends smoothly on $\hat e$ as $\Xi_i^{\hat e}$
does, and \eq{cm5eq3} follows immediately from $\Up_i^*(\om)=\om'$
and parts (i)--(iii) above. Regard \eq{cm5eq4} as a {\it definition}
of $\Phi_\sXp^{\hat e}$ on $\pi^*(S_i)\subset U_\sXp$ for
$i=1,\ldots,n$, and define $\Phi_\sXp^{\hat e}\equiv\Phi_\sXp$ on
$\pi^*(K)\subset U_\sXp$. Then $\Phi_\sXp^{\hat e}:U_\sXp\ra M$
is well-defined, and satisfies~\eq{cm5eq4}.

To see that $\Phi_\sXp^{\hat e}$ is smooth, we need to
show that its definitions on $\pi^*(S_i)$ and $\pi^*(K)$ join
together smoothly on $\pi^*(\pd K)$. This follows from part
(iv) above provided $\Phi_\sXp\bigl(\pi^*(\pd K)\bigr)$
does not intersect $\Up_i(B_{2R''})$, since
then when $r$ is close to $R'$ in \eq{cm5eq4} we have
$\Up_i^{\hat e}=\Up_i$, and thus $\Phi^{\hat e}_\sXp=
\Phi_\sXp$ near the boundary of $\pi^*(S_i)$ where it
joins onto~$\pi^*(K)$.

Hence, choosing $R''\in(0,\ha R)$ such that $\Phi_\sXp
\bigl(\pi^*(\pd K)\bigr)$ does not intersect $\Up_i(B_{2R''})$
for $i=1,\ldots,n$ ensures that $\Phi_\sXp^{\hat e}$ is
smooth for all $\hat e\in\ti\E$, and making $\ti\E$
smaller if necessary we can assume it is an embedding. As
$\Phi_\sXp^*(\om)=\hat\om$ we see that $(\Phi^{\hat
e}_\sXp)^*(\om)=\hat\om$ on $\pi^*(K)$, and $(\Phi^{\hat
e}_\sXp)^*(\om)=\hat\om$ on $\pi^*(S_i)$ follows from
\eq{cm5eq4} since $(\Up_i^{\hat e})^*(\om)=\om'$. Finally,
$\Phi_\sXp^e=\Phi_\sXp$ as $\Up_i^e=\Up_i$
for~$i=1,\ldots,n$.
\end{proof}

In the situation of the theorem, fix $\hat e\in\ti\E$
and define $\hat X'=\Phi^{\hat e}_\sXp(X')$, where
$X'\subset U_\sXp\subset T^*X'$ is the zero section,
and set $\hat X=\hat X'\cup\{\hat x_1,\ldots,\hat x_n\}$.
As $(\Phi^{\hat e}_\sXp)^*(\om)=\hat\om$ it follows
that $\hat X'$ is a Lagrangian submanifold of $M$, and thus
$\hat X$ is a compact {\it Lagrangian} $m$-fold in $M$
with conical singularities at $\hat x_1,\ldots,\hat x_n$,
identifications $\hat\up_1,\ldots,\hat\up_n$ and cones
$C_1,\ldots,C_n$, generalizing Definition \ref{cm3def4}
in the obvious way.

Thus we have extended $X$ to a smooth family of Lagrangian
$m$-folds $\hat X$ with conical singularities, which realize
all nearby alternative choices of $x_i,\up_i$ exactly once
up to equivalence. When $\hat e$ is close to $e$, $\hat X$
will be approximately special Lagrangian, and so we can try
to deform it to an exactly special Lagrangian $m$-fold with
the same~$\hat x_i,\hat\up_i$.

\subsection{Small deformations of $X$ and moduli spaces}
\label{cm52}

Suppose that $(M,J,\om,\Om)$ is an almost Calabi--Yau $m$-fold and
that $X,\hat X$ are compact SL $m$-folds in $M$ which both have $n$
conical singular points $x_1,\ldots,x_n$ and $\hat x_1,\ldots,\hat
x_n$ respectively, with the same cones $C_1,\ldots,C_n$ and rates
$\mu_1,\ldots,\mu_n$. When $X,\hat X$ are `sufficiently close' in
a $C^1$ sense we shall write $\hat X$ in terms of a small closed
1-form $\al$ on $X'$ with prescribed decay, using the Lagrangian
neighbourhood $\Phi^{\hat e}_\sXp$ of Theorem \ref{cm5thm1}.
Thus we shall define a topology on the set of compact SL $m$-folds
in $M$ with conical singularities.

\begin{thm} Suppose $(M,J,\om,\Om)$ is an almost Calabi--Yau $m$-fold
and\/ $X$ a compact SL\/ $m$-fold in $M$ with conical singularities
at\/ $x_1,\ldots,x_n$, with identifications $\up_i$, cones $C_i$ and
rates $\mu_i$. Let\/ $e,\E$ be as in Definition \ref{cm5def1}, and\/
$U_\sXp,\Phi_\sXp,\ti\E,\Up_i^{\hat e}$ and\/ $\Phi_\sXp^{\hat e}$
be as in Theorem~\ref{cm5thm1}.

Let\/ $\hat e=(\hat x_1,\hat\up_1,\ldots,\hat x_n,\hat\up_n)\in
\ti\E$, and suppose $\hat X$ is a compact SL\/ $m$-fold in
$M$ with conical singularities at\/ $\hat x_1,\ldots,\hat x_n$,
with identifications $\hat\up_i$, cones $C_i$ and rates $\mu_i$.
Then if\/ $\hat e,e$ are sufficiently close in $\ti\E$ and\/
$X',\hat X'$ are sufficiently close as submanifolds in a $C^1$
sense away from $x_1,\ldots,x_n$, there exists a closed\/
$1$-form $\al$ on $X'$ such that the graph\/ $\Ga(\al)$ lies in
$U_\sXp\subset T^*X'$, and\/ $\hat X'=\Phi^{\hat e}_\sXp
\bigl(\Ga(\al)\bigr)$. Furthermore we may write $\al=\be+\d f$,
where $\be$ is a closed\/ $1$-form supported in $K$ and\/~$f\in
C^\iy_{\bs\mu}(X')$.
\label{cm5thm2}
\end{thm}

\begin{proof} Apply Theorem \ref{cm4thm5} to $X$ and $\hat X$,
using $\Up_i^e=\Up_i$ for $X$ and $\Up_i^{\hat e}$ for $\hat X$,
and the same $R,\ze,U_{\sst C_i}$ and $\Phi_{\sst C_i}$ for
both. Theorem \ref{cm4thm5} then gives $R',\hat R'\in(0,R]$
and closed 1-forms $\eta_i$ on $\Si_i\t(0,R')$ and $\hat\eta_i$
on $\Si_i\t(0,\hat R')$ for $i=1,\ldots,n$ such that $X',\hat X'$
are parametrized on $S_i,\hat S_i$ using maps $\phi_i:\Si_i\t(0,R')
\ra B_R$ and $\hat\phi_i:\Si_i\t(0,\hat R')\ra B_R$ defined using
$\eta_i,\hat\eta_i$ in the usual way.

Theorem \ref{cm4thm8} defines real functions $A_i$ on $\Si_i\t(0,R')$
and $\hat A_i$ on $\Si_i\t(0,\hat R')$ with $\eta_i=\d A_i$ and
$\hat\eta_i=\d\hat A_i$, and proves results on the decay of
$\phi_i,\eta_i,A_i$ and $\hat\phi_i,\hat\eta_i,\hat A_i$ and
their derivatives. Using \eq{cm4eq10} and $\mu_i>2$ we see that
$\eta_i,\hat\eta_i=o(r)$ for small $r$. Therefore we may choose
$R''\in\bigl(0,\min(R',\hat R')\bigr]$ such that $\md{\hat\eta_i
-\eta_i}<\ze r$ on $\Si_i\t(0,R'')$ for all~$i=1,\ldots,n$.

Let $S_i'=\Up_i\circ\phi_i\bigl(\Si_i\t(0,R'')\bigr)$ and
$\hat S_i'=\Up^{\hat e}_i\circ\hat\phi_i\bigl(\Si_i\t(0,R'')\bigr)$
for $i=1,\ldots,n$, so that $S_i'\subseteq S_i\subset X'$ and
$\hat S_i'\subseteq\hat S_i\subset\hat X'$. Define a 1-form $\al$
on $S_i'$ by $\al=(\Up_i\circ\phi_i)_*(\hat\eta_i-\eta_i)$ for
$i=1,\ldots,n$. Now as $\hat\phi_i(\si,r)=\Phi_{\sst C_i}\bigl(
\si,r,\hat\eta_i^1(\si,r),\hat\eta_i^2(\si,r)\bigr)$ by Theorem
\ref{cm4thm5}, we see from \eq{cm5eq4} that if $(\si,r)\in\Si_i
\t(0,R'')$ and $(\tau,u)=(\hat\eta_i^1-\eta_i^1,\hat\eta_i^2-
\eta_i^2)(\si,r)$ then
\begin{equation*}
\Phi^{\hat e}_\sXp\bigl[\al\bigl(\Up_i\circ\phi_i(\si,r)\bigr)
\bigr]=\Phi^{\hat e}_\sXp\circ\d(\Up_i\circ\phi_i)(\si,r,\tau,u)
=\Up_i^{\hat e}\circ\hat\phi_i(\si,r)\in\hat S_i'\subset\hat X'.
\end{equation*}

Thus the subsets $\hat S_i'$ in $\hat X'$ coincide with
$\Phi^{\hat e}_\sXp\bigl(\Ga(\al)\bigr)$ on the subsets
$S_i'$ in $X'$ where $\al$ is defined so far. To show that
$\hat X'=\Phi^{\hat e}_\sXp\bigl(\Ga(\al)\bigr)$ for
some 1-form $\al$ defined on the whole of $X'$, we need that
\begin{itemize}
\item[(a)] $\hat X'$ should lie in $\Phi^{\hat e}_\sXp
(U_\sXp)$, and
\item[(b)] $\hat X'$ should intersect the image under
$\Phi^{\hat e}_\sXp$ of each fibre of $\pi:U_\sXp
\ra X'$ transversely exactly once.
\end{itemize}
We have already shown that (a) and (b) hold on the subsets~$\hat S_i'$.

Under the assumptions of the theorem $\hat e,e$ are close in $\ti\E$
and $\Phi^{\hat e}_\sXp$ and $\Phi^e_\sXp=\Phi_\sXp$ are close on
the complement of the $S_i'$. Also $X',\hat X'$ are close as
submanifolds in a $C^1$ sense away from $x_1,\ldots,x_n$, and
thus on the complement of the subsets $S_i'$ in $X'$ and
$\hat S_i'$ in $\hat X'$ for $i=1,\ldots,n$. Therefore $\hat X'$
satisfies (a) and (b) on the complement of the $\hat S_i'$, and $\al$
exists. Since $\hat X'$ is Lagrangian and $(\Phi^{\hat e}_\sXp)^*(\om)
=\hat\om$, the usual argument shows that $\al$ is closed.

Define a smooth real function $f$ on $S_i'$ by
$f=(\Up_i\circ\phi_i)_*(\hat A_i-A_i)$ for $i=1,\ldots,n$.
Then $\al=\d f$ on $S_i'$, as $\eta_i=\d A_i$ and $\hat\eta_i=
\d\hat A_i$. As $\al$ is closed and $S_i'\subseteq S_i$ are
homotopy equivalent we can extend $f$ uniquely to $S_i$ with
$\al=\d f$. Then extend $f$ smoothly over $K$. This defines a
smooth function $f$ on $X'$ with $\al=\d f$ on $S_i$ for
$i=1,\ldots,n$. Let $\be=\al-\d f$. Then $\al=\be+\d f$ and $\be$
is a closed 1-form supported in $K=X'\sm(S_1\cup\cdots\cup S_n)$,
as we have to prove. Finally, \eq{cm4eq14} for $A_i,\hat A_i$ with
$\mu_i'=\mu_i$ gives~$f\in C^\iy_{\bs\mu}(X')$.
\end{proof}

We define the {\it moduli space} $\M_\sX$ of SL $m$-folds
$\hat X$ with conical singularities in $M$, which are isotopic
to $X$ in $M$ and have the same cones~$C_1,\ldots,C_n$.

\begin{dfn} Let $(M,J,\om,\Om)$ be an almost Calabi--Yau
$m$-fold and $X$ a compact SL $m$-fold in $M$ with conical
singularities at $x_1,\ldots,x_n$ with identifications
$\up_i:\C^m\ra T_{x_i}M$ and cones $C_1,\ldots,C_n$. Define
the {\it moduli space} $\M_\sX$ {\it of deformations of\/}
$X$ to be the set of $\hat X$ such that
\begin{itemize}
\setlength{\parsep}{0pt}
\setlength{\itemsep}{0pt}
\item[(i)] $\hat X$ is a compact SL $m$-fold in $M$ with
conical singularities at $\hat x_1,\ldots,\hat x_n$ with
cones $C_1,\ldots,C_n$, for some $\hat x_i$ and
identifications~$\hat\up_i:\C^m\ra T_{\smash{\hat x_i}}M$.
\item[(ii)] There exists a homeomorphism $\hat\iota:X\ra\hat X$
with $\hat\iota(x_i)=\hat x_i$ for $i=1,\ldots,n$ such that
$\hat\iota\vert_{X'}:X'\ra\hat X'$ is a diffeomorphism and
$\hat\iota$ and $\iota$ are isotopic as continuous maps
$X\ra M$, where $\iota:X\ra M$ is the inclusion.
\end{itemize}
Note that by Theorem \ref{cm4thm8} the definition of $\hat X$
is independent of choice of rates $\mu_i$, so there is no
need to include the $\mu_i$ in~(i).

Let $\V_\sX$ be the subset of $\hat X\in\M_\sX$ such that for some
$\hat e=(\hat x_1,\hat\up_1,\ldots,\hat x_n,\hat\up_n)$ in $\ti\E$
and some 1-form $\al$ on $X'$ whose graph $\Ga(\al)$ lies in
$U_\sXp\subset T^*X'$ we have $\hat X'=\Phi^{\hat e}_\sXp\bigl(
\Ga(\al)\bigr)$, as in Theorem \ref{cm5thm2}. Note that if
$\ti X\in\M_\sX$ then $\M_{\sst\ti X}=\M_\sX$. Thus, for each
$\ti X\in\M_\sX$ we have~$\ti X\in\V_{\sst\ti X}\subset\M_\sX$.
\label{cm5def2}
\end{dfn}

The construction of $\V_\sX$ above gives a 1-1
correspondence between $\V_\sX\subseteq\M_\sX$ and
a set of pairs $(\hat e,\al)$ for $\hat e\in\ti\E$
and $\al$ a smooth 1-form on $X'$ with prescribed decay.
Using the given topology on $\ti\E$ and a suitable
choice of topology on the 1-forms $\al$, this 1-1
correspondence induces a {\it topology} on~$\V_\sX$.

To define the $\al$ topology, choose some $\bs\mu$ as in
Definition \ref{cm3def4}, and let the $C^k_{{\bs\mu}-1}$
topology on $\al$ be induced by the norm
\begin{equation*}
\nm{\al}_{C^k_{{\bs\mu}-1}}=\sum_{j=0}^k
\sup_{X'}\bmd{\rho^{-\bs{\mu}+1-j}\na^j\al},
\end{equation*}
and the $C^\iy_{{\bs\mu}-1}$ topology on $\al$ be induced by
the $C^k_{{\bs\mu}-1}$ topologies for all~$k\ge 0$.

\begin{prop} The $C^1_{{\bs\mu}-1}$ and\/ $C^\iy_{{\bs\mu}-1}$
topologies on $\al$ induce the same topology on $\V_\sX$, which
is also independent of the choice of rates~$\bs\mu$.
\label{cm5prop}
\end{prop}

\begin{proof} This is implicit in the proofs of Theorems
\ref{cm4thm8} and \ref{cm5thm2}. In particular, Theorem
\ref{cm4thm8} in effect shows that an a priori estimate
for the $C^1_{{\bs\mu}-1}$ norm of $\al$ implies a priori
estimates for the $C^k_{{\bs\mu}-1}$ norms for all $k\ge 1$,
and so the $C^1_{{\bs\mu}-1}$ and $C^\iy_{{\bs\mu}-1}$
topologies on $\al$ induce the same topology on
$\V_\sX$. It also proves independence of
the choice of~$\bs\mu$.
\end{proof}

We can now define a {\it topology} on~$\M_\sX$.

\begin{dfn} For each $\ti X\in\M_\sX$, use the 1-1
correspondence between $\V_{\sst\ti X}$ and pairs
$(\hat e,\al)$ to define a topology on $\V_{\sst\ti X}$
as in Proposition \ref{cm5prop}. We get the same topology using
the $C^1_{{\bs\mu}-1}$ or $C^\iy_{{\bs\mu}-1}$ topologies on
$\al$ for any choice of $\bs\mu$, so there is no ambiguity.
One can show that overlaps $\V_{\sst X_1}\cap\V_{\sst X_2}$
are open in $\V_{\sst X_j}$ and the $\V_{\sst X_j}$ topologies
agree on the overlaps. Piecing the topologies together therefore
defines a unique topology on~$\M_\sX$.
\label{cm5def3}
\end{dfn}

\noindent{\it Remarks.} Basically, $\M_\sX$ is the
family of compact SL $m$-folds $\hat X$ in $M$ with
conical singularities which are deformation equivalent
to $X$ in a loose sense. Note that $\M_\sX$ may not
be {\it connected}, as the isotopies in part (ii)
of Definition \ref{cm5def2} need not be through special
Lagrangian embeddings.

In Theorem \ref{cm5thm2} we assumed only that $\hat e,e$ are
close in $\ti\E$ and that $X',\hat X'$ are `sufficiently
close as submanifolds in a $C^1$ sense away from $x_1,\ldots,x_n$'.
These closeness assumptions are actually {\it very weak}, in that
we have imposed {\it no asymptotic conditions} on how $X',\hat X'$
converge to $x_i$ and $\hat x_i$, but instead required only $C^1$
closeness on large compact subsets of~$X',\hat X'$.

Because of this, we can be confident that the topology defined
on $\M_\sX$ above is a sensible choice. In particular, Theorem
\ref{cm5thm2} effectively shows that if $X,\hat X$ are close in
a very weak sense, then they are close in the $\M_\sX$ topology.
Theorem \ref{cm6thm3} below gives another way of seeing the
naturality of the topology on~$\M_\sX$.

Definitions \ref{cm5def2} and \ref{cm5def3} don't actually need
$X$ to be special Lagrangian in $(M,J,\om,\Om)$, except to ensure
that $X\in\M_\sX$. We are simply using $X$ to fix the topological
type, isotopy class and singular cones $C_i$ of $\hat X\in\M_\sX$.
In particular, given a family $\bigl\{(M,J^s,\om^s,\Om^s):s\in\F
\bigr\}$ of almost Calabi--Yau structures on $M$ with $X$ special
Lagrangian in $(M,J^0,\om^0,\Om^0)$, we can define a moduli space
$\M_\sX^s$ of special Lagrangian deformations of $X$ in
$(M,J^s,\om^s,\Om^s)$, for each~$s\in\F$.

\section{Deformations, obstructions, and smoothness}
\label{cm6}

We can now prove the first main result of the paper, Theorem
\ref{cm6thm2} below, which is an analogue of McLean's Theorem,
Theorem \ref{cm2thm1}, for compact SL $m$-folds $X$ with conical
singularities $x_1,\ldots,x_n$ in a single almost Calabi--Yau
$m$-fold $(M,J,\om,\Om)$. An important difference with the
nonsingular case is that there may be {\it obstructions} to
deforming $X$, which means that the moduli space
$\M_\sX$ may be singular.

Instead, $\M_\sX$ is locally homeomorphic by a map
$\Xi$ to the zeroes of a smooth map $\Phi:\I_\sXp\ra
\O_\sXp$ between finite-dimensional vector spaces
$\I_\sXp$, the {\it infinitesimal deformation space},
and $\O_\sXp$, the {\it obstruction space}. Here
$\I_\sXp$ is isomorphic to the image of $H^1_{\rm
cs}(X',\R)$ in $H^1(X',\R)$, and $\O_\sXp$ is a
direct sum of subspaces depending on the SL cones $C_1,\ldots,C_n$
of $X$ at~$x_1,\ldots,x_n$.

We set up the problem in \S\ref{cm61}, and define $\O_\sXp$
in \S\ref{cm62}. The main theorem is proved in \S\ref{cm63}, with some
corollaries on cases when $\M_\sX$ is smooth. Section
\ref{cm64} discusses the naturality (independence of choices) of
$\I_\sXp,\O_\sXp,\Phi$ and $\Xi$, and
\S\ref{cm65} another way to define $\I_\sXp$
and~$\O_\sXp$.

\subsection{Setting up the deformation problem}
\label{cm61}

We shall parametrize the moduli space $\M_\sX$ locally
in terms of the zeroes of a map $F$ between Banach spaces.

\begin{dfn} Let $(M,J,\om,\Om)$ be an almost Calabi--Yau
$m$-fold and $X$ a compact SL $m$-fold in $M$ with conical
singularities at $x_1,\ldots,x_n$ with identifications
$\up_i:\C^m\ra T_{x_i}M$ and cones $C_1,\ldots,C_n$. Let
$U_\sXp,\Phi_\sXp$ be as in Theorem \ref{cm4thm6},
and $e,\E$ as in Definition \ref{cm5def1}, and $\ti\E$,
$\Up_i^{\hat e}$ and $\Phi^{\hat e}_\sXp$ as in
Theorem~\ref{cm5thm1}.

Choose a vector space $\H_\sXp$ of closed 1-forms on
$X'$ supported in $K$, such that the map $\H_\sXp\ra
H^1_{\rm cs}(X',\R)$ given by $\be\mapsto[\be]$ is an isomorphism.
Since $X'$ retracts onto $K$, this is clearly possible. Now the
subspace of $\H_\sXp$ corresponding to the kernel of
the map $H^1_{\rm cs}(X',\R)\ra H^1(X',\R)$ in \eq{cm4eq4}
consists of {\it exact\/} 1-forms on $X'$, so each such 1-form
may be written $\d v$ for some~$v\in C^\iy(X')$.

Let the connected components of $S_i\cong\Si_i\t(0,R')$ be $S_i^j$
for $j=1,\ldots,b^0(\Si_i)$. As $\d v=0$ on $S_i$ we see that
$v=a_i^j$ on $S_i^j$ for some constants $a_i^j$. Since $v$ is
defined up to addition of a constant, we specify $v$ uniquely
by requiring that~$\sum_{i,j}a_i^j=0$.

Define $\K_\sXp$ to be the vector space of all such functions
$v$. Then $\d\K_\sXp=\{\d v:v\in\K_\sXp\}$ is a subspace
of $\H_\sXp$, and $\d:\K_\sXp\ra\H_\sXp$
is injective. Also $\K_\sXp$ is isomorphic to the kernel of
$H^1_{\rm cs}(X',\R)\ra H^1(X',\R)$ in \eq{cm4eq4}. Thus by
\eq{cm4eq4} we have an exact sequence
\begin{equation*}
0\ra H^0(X',\R)\ra\bigoplus_{i=1}^nH^0(\Si_i,\R)\ra\K_\sXp\ra 0,
\end{equation*}
so as $X'$ is connected we see that
\e
\dim\K_\sXp=\sum_{i=1}^nb^0(\Si_i)-1.
\label{cm6eq1}
\e

Let the {\it infinitesimal deformation space} $\I_\sXp$
be a vector subspace of $\H_\sXp$ with 
\e
\H_\sXp=\I_\sXp\op\d\K_\sXp.
\label{cm6eq2}
\e
As $\d\K_\sXp$ corresponds to the kernel of $H^1_{\rm cs}(X',\R)
\ra H^1(X',\R)$ in \eq{cm4eq4} and $\I_\sXp\cong{\cal
H}_\sXp/\d\K_\sXp$, we see that the map $\I_\sXp
\ra H^1(X',\R)$ given by $\be\mapsto[\be]$ is an isomorphism between
$\I_\sXp$ and the image of $H^1_{\rm cs}(X',\R)$ in~$H^1(X',\R)$.

Let $k>2$, $p>m$, and $\bs\mu$ be as in Definition \ref{cm3def4}.
Then $L^p_{k,\bs{\mu}}(X')$ is continuously included in
$C^2_{\bs\mu}(X')$ by Theorem \ref{cm4thm1}. Define
\e
\D_\sXp=\bigl\{(\be,f)\in\H_\sXp\t
L^p_{k,{\bs\mu}}(X'): \text{the graph of $\be+\d f$ lies
in $U_\sXp$}\bigr\}.
\label{cm6eq3}
\e
Then $\D_\sXp$ is an open subset of $\H_\sXp
\t L^p_{k,{\bs\mu}}(X')$ containing $(0,0)$. Here we use the fact
that $f$ is $C^1$ to make sense of the graph of~$\be+\d f$.

Define a map $F:\ti\E\t\D_\sXp\ra C^0(X')$ by
\e
\pi_*\bigl((\Phi^{\hat e}_\sXp)^*(\Im\Om)\vert_{\Ga(\be+\d f)}\bigr)
=F(\hat e,\be,f)\,\d V_g,
\label{cm6eq4}
\e
where $\Ga(\be+\d f)$ is the graph of $\be+\d f$ in $U_\sXp$,
and $\pi:\Ga(\be+\d f)\ra X'$ the natural projection, and
$\d V_g$ the volume form of the metric $g$ on $X'$. Since
$f$ is $C^2$ we see that $\Ga(\be+\d f)$ is a $C^1$-submanifold
of $U_\sXp$, and so $(\Phi^{\hat e}_\sXp)^*(\Im\Om)
\vert_{\Ga(\be+\d f)}$ makes sense and its image under $\pi$
is continuous. Hence $F(\hat e,\be,f)$ lies in $C^0(X')$,
the vector space of continuous functions on~$X'$.
\label{cm6def1}
\end{dfn}

The point of the definition is given in the following proposition.

\begin{prop} In the situation of Definition \ref{cm6def1},
suppose $(\hat e,\be,f)\in\ti\E\t\D_\sXp$ with\/
$F(\hat e,\be,f)=0$. Set\/ $\hat X'=\smash{\Phi^{\hat e}_\sXp}
\bigl(\Ga(\be+\d f)\bigr)$ and\/ $\hat X=\hat X'\cup\{\hat x_1,
\ldots,\hat x_n\}$, where $\hat e=(\hat x_1,\ldots,\hat\up_n)$.
Then $f\in C^\iy_{\bs\mu}(X')$ and\/ $\hat X$ is a compact SL\/
$m$-fold in $M$ with conical singularities at\/ $\hat x_i$ with
identifications $\hat\up_i$, cones $C_i$ and rates~$\mu_i$.

Thus $\hat X$ lies in $\V_\sX\subseteq\M_\sX$ in Definition
\ref{cm5def2}. Conversely, each\/ $\hat X$ in $\V_\sX$ comes
from a unique $(\hat e,\be,f)\in\ti\E\t\D_\sXp$ with\/
$F(\hat e,\be,f)=0$. Write $\Psi(\hat e,\be,f)=\hat X$.
Then $\Psi:F^{-1}(0)\ra\V_\sX$ is a homeomorphism,
with\/~$\Psi(e,0,0)=X$.
\label{cm6prop1}
\end{prop}

\begin{proof} Suppose $F(\hat e,\be,f)=0$. Then $f\in
C^2_{\bs\mu}(X')$ from above, so $f$ is locally $C^2$ and
$\hat X'$ is a $C^1$ submanifold of $M$. As $(\Phi^{\hat
e}_\sXp)^*(\om)=\hat\om$ and $\be+\d f$ is a closed
$C^1$ 1-form, we see that $\om\vert_{\hat X'}\equiv 0$ by
the usual argument. Also \eq{cm6eq4} implies that
$\Im\Om\vert_{\hat X'}\equiv 0$. Therefore, if we can
prove that $\hat X'$ is a $C^\iy$ submanifold of $M$ then
$\hat X'$ is special Lagrangian, by Definition~\ref{cm2def4}.

With $\hat e,\be$ fixed $F(\hat e,\be,f)$ depends pointwise
on $\d f,\na^2f$ by \eq{cm6eq4}, so
\e
F(\hat e,\be,f)[x]=F'\bigl(x,\d f(x),\na^2f(x)\bigr)=0,
\label{cm6eq5}
\e
where $F'$ is a smooth, nonlinear function of its arguments
defined on some domain. Now \eq{cm6eq5} is a second-order
nonlinear p.d.e., and using the ideas of \S\ref{cm23} one can
show that it is {\it elliptic}. Aubin \cite[Th.~3.56]{Aubi}
gives an elliptic regularity result for such equations which
shows that if $f$ is locally $C^2$ then $f$ is locally $C^\iy$.
Thus $f$ is smooth, so $\hat X'$ is $C^\iy$ and thus special
Lagrangian.

Recall that $A_i$ is a function and $\eta_i=\d A_i$ a 1-form
on $\Si_i\t(0,R')$ for $i=1,\ldots,n$, defined in Theorems
\ref{cm4thm5} and \ref{cm4thm8}, and that $\Up_i\circ\phi_i:
\Si_i\t(0,R')\ra S_i\subset X'$ is a diffeomorphism. Define
$\hat A_i$ and $\hat\eta_i$ on $\Si_i\t(0,R')$ by
\e
\hat A_i=f\circ\Up_i\circ\phi_i+A_i
\quad\text{and}\quad
\hat\eta_i=\d\hat A_i=\d(f\circ\Up_i\circ\phi_i)+\eta_i.
\label{cm6eq6}
\e
Let $\hat\eta_i^1,\hat\eta_i^2$ be the components of $\hat\eta_i$
as in Theorem \ref{cm4thm5}, and define
\e
\hat\phi_i:\Si_i\t(0,R')\ra B_R
\quad\text{by}\quad
\hat\phi_i(\si,r)=\Phi_{\sst C_i}\bigl(\si,r,\hat\eta_i^1
(\si,r),\hat\eta_i^2(\si,r)\bigr).
\label{cm6eq7}
\e

Combining \eq{cm4eq10}, \eq{cm4eq14}, \eq{cm6eq7} and
$f\in C^2_{\bs\mu}(X')$ from above, we prove that
\e
\begin{gathered}
\bmd{\na^k(\hat\phi_i-\iota_i)}=O(r^{\mu_i-1-k})
\quad\text{and}\quad
\bmd{\na^k\hat\eta_i}=O(r^{\mu_i-1-k})
\quad\text{for $k=0,1$}\\
\text{and}\quad
\bmd{\na^k\hat A_i}=O(r^{\mu_i-k})
\quad\text{for $k=0,1,2$, as $r\ra 0$ for $i=1,\ldots,n$.}
\end{gathered}
\label{cm6eq8}
\e
Using \eq{cm5eq4} and the facts that $\hat X'=\Phi^{\hat
e}_\sXp\bigl(\Ga(\be+\d f)\bigr)$ and $\be=0$ in $S_i$, we
find that $\Up_i^{\hat e}\circ\hat\phi_i:\Si_i\t(0,R')\ra M$
maps into $\hat X'$, and defines a diffeomorphism $\Si_i\t(0,R')
\ra\hat S_i$ with an open subset $\hat S_i$ of $\hat X'$. Also
the natural diffeomorphism $X'\ra\hat X'$ identifies $S_i$ and
$\hat S_i$, and thus $\hat K=\hat X'\sm(\hat S_1\cup\cdots\cup
\hat S_n)$ is compact.

Therefore all the conditions of Definition \ref{cm3def4} are
satisfied, and so $\hat X$ is a compact SL $m$-fold in $M$ with
conical singularities at $\hat x_1,\ldots,\hat x_n$, with
identifications $\hat\up_i$, cones $C_i$ and rates $\mu_i$, as
we have to prove. Applying Theorem \ref{cm4thm8} to $X$ and
$\hat X$ then shows that $\md{\na^kA_i}=O(r^{\mu_i-k})$ and
$\md{\na^k\hat A_i}=O(r^{\mu_i-k})$ for all $k\ge 0$. Thus
\eq{cm6eq6} gives $\md{\na^kf}=O(\rho^{\mu_i-k})$ on $S_i$
for all $k\ge 0$ and $i=1,\ldots,n$. Since $f$ is smooth this
implies that $f\in C^\iy_{\bs\mu}(X')$, as we have to prove.

Definition \ref{cm5def2} now shows that $\hat X\in\V_\sX$.
Conversely, if $\hat X\in\V_\sX$ then Definition \ref{cm5def2}
gives $\hat X'=\Phi^{\hat e}_\sXp\bigl(\Ga(\al)\bigr)$ for
some $\hat e\in\ti\E$ and 1-form $\al$ on $X'$ whose graph
$\Ga(\al)$ lies in $U_\sXp$. The proof of Theorem \ref{cm5thm2}
then shows that $\al=\ti\be+\d\ti f$, where $\ti\be$ is a closed
$1$-form supported in $K$ and~$\ti f\in C^\iy_{\bs\mu}(X')$.

Let $\be$ be the unique element of $\H_\sXp$ with
$[\be]=[\ti\be]$ in $H^1_{\rm cs}(X',\R)$, where $\H_{\sst
X'}$ is as in Definition \ref{cm6def1}. Then $\ti\be-\be=\d\ga$,
for $\ga\in C^\iy_{\rm cs}(X')$. Set $f=\ti f+\ga$. Then $f\in
C^\iy_{\bs\mu}(X')$ with $\al=\be+\d f$. Theorem \ref{cm4thm8}
shows that we can improve the rates $\mu_i$ of the singularities
$\hat x_i$ of $\hat X$ to some rates $\mu_i'>\mu_i$ for
$i=1,\ldots,n$. It follows that $f\in C^\iy_{{\bs\mu}'}(X')$, and
therefore $f\in L^p_{k,{\bs\mu}}(X')$ as $C^\iy_{{\bs\mu}'}(X')
\subset L^p_{k,{\bs\mu}}(X')$. Therefore $(\be,f)\in\D_\sXp$
by~\eq{cm6eq3}.

As $\hat X'$ is special Lagrangian $\Im\Om\vert_{\hat X'}\equiv 0$,
and it follows from \eq{cm6eq4} that $F(\hat e,\be,f)=0$. Thus each
$\hat X$ in $\V_\sX$ comes from some $(\hat e,\be,f)\in\ti\E\t
\D_\sXp$ with $F(\hat e,\be,f)=0$. Since there are no nontrivial
$G_1\t\cdots\t G_n$ equivalences in $\ti\E$ by construction,
$\hat X$ determines $\hat e$ uniquely, and $\hat X,\hat e$ then
determine $\al$ and so $\be,f$ uniquely. Thus $(\hat e,\be,f)$ is
unique.

Thus writing $\Psi(\hat e,\be,f)=\hat X$ defines a bijection
$\Psi:F^{-1}(0)\ra\V_\sX$ with $\Psi(e,0,0)=X$. We
must show that $\Psi$ is a {\it homeomorphism}. The topology
on $\V_\sX$ is defined using pairs $(\hat e,\al)$,
where $\hat e$ has the $\ti\E$ topology and $\al$ either
the $C^1_{{\bs\mu}-1}$ or the $C^\iy_{{\bs\mu}-1}$ topology
on 1-forms for any choice of $\bs\mu$, and $\Psi$
takes~$(\hat e,\be,f)\mapsto(\hat e,\be+\d f)$.

Now $f$ has the $L^p_{k,{\bs\mu}}$ topology, so $\d f$ has the
$L^p_{k-1,{\bs\mu}-1}$ topology. This is intermediate between
the $C^1_{{\bs\mu}-1}$ and $C^\iy_{{\bs\mu}'-1}$ topologies on
$\al$ for $\mu_i'>\mu_i$ as above, as $C^\iy_{{\bs\mu}'-1}(T^*X')
\subset L^p_{k-1,{\bs\mu}-1}(T^*X')\subset C^1_{{\bs\mu}-1}(T^*X')$
by Theorem \ref{cm4thm1}. But the $C^1_{{\bs\mu}-1}$ and
$C^\iy_{{\bs\mu}'-1}$ topologies on $\al$ induce the same
topology on $\V_\sX$ by Proposition \ref{cm5prop}. Thus the
$L^p_{k-1,{\bs\mu}-1}$ topology on $\d f$ also induces
the same topology on $\V_\sX$, and it follows quickly
that $\Psi$ is a homeomorphism.
\end{proof}

Here is an analogue of Proposition \ref{cm2prop2} for~$F$.

\begin{prop} In the situation above, for $x\in X'$ we may write
\e
F(\hat e,\be,f)[x]\!=\!-\d^*\bigl(\psi^m(\be\!+\!\d f)\bigr)[x]
\!+\!Q\bigl(\hat e,x,(\be\!+\!\d f)(x),(\na\be\!+\!\na^2f)(x)\bigr),
\label{cm6eq9}
\e
where $Q:\bigl\{(\hat e,x,y,z):\hat e\in\ti\E$, $x\in X'$,
$y\in T^*_xX'\cap U_\sXp$, $z\in\ot^2T^*_xX'\bigr\}\ra\R$ is
smooth, and for $\rho(x)^{-1}\md{y}$, $\md{z}$ and\/ $d_{\sst\cal
E}(\hat e,e)$ small we have
\e
Q(\hat e,x,y,z)=O\bigl(\rho(x)^{-2}\ms{y}+\ms{z}+
\rho(x)d_\sE(\hat e,e)\bigr),
\label{cm6eq10}
\e
and more generally for $\rho(x)^{-1}\md{y}$, $\md{z}$ and\/
$d_\sE(\hat e,e)$ small and\/ $a,b,c\ge 0$ we have
\e
\begin{split}
(\na_x)^a(\pd_y)^b(\pd_z)^cQ(\hat e,x,y,z)&=
O\bigl(\rho(x)^{-a-\max(2,b)}\md{y}^{\max(0,2-b)}\\
&+\rho(x)^{-a}\md{z}^{\max(0,2-c)}
+\rho(x)^{1-a-b}d_\sE(\hat e,e)\bigr),
\end{split}
\label{cm6eq11}
\e
where $\na_x,\pd_y,\pd_z$ are the partial derivatives of\/ $Q$
in the $x,y,z$ variables, using the Levi-Civita connection
$\na$ of\/ $g$ to form~$\na_x$.
\label{cm6prop2}
\end{prop}

\begin{proof} The value of $F(\hat e,\be,f)$ at $x\in X'$ depends
on $\hat e$, via $(\Phi^{\hat e}_\sXp)^*(\Im\Om)$, and on the
tangent space to $\Ga(\be+\d f)$ at $x'$, where $x'\in\Ga(\al)$
with $\pi(x')=x$. But $T_{x'}\Ga(\be+\d f)$ depends on both
$(\be+\d f)\vert_x$ and $(\na\be+\na^2 f)\vert_x$. Therefore
$F(\hat e,\be,f)$ depends pointwise on the arguments of
$Q$ in~\eq{cm6eq9}.

As in the proof of Proposition \ref{cm2prop2} we may take
\eq{cm6eq9} as a {\it definition} of $Q$, and $Q$ is then
well-defined on the given domain, which is the set of all
$\hat e,x,y,z$ realized by $\hat e,\be,f$ in the domain of
$F$. As $\pi,\psi,\Im\Om,\d V_g$ are smooth and $\Phi^{\hat
e}_\sXp$ is smooth and depends smoothly on $\hat e$,
we see that $Q$ is a smooth function of its arguments.

Since $\Phi^e_\sXp=\Phi_\sXp$ and $\Phi_\sXp$ is
the identity on $X'=\Ga(0)\subset U_\sXp$ we see that
$F(e,0,0)\,\d V_g=\Im\Om\vert_{X'}=0$ as $X'$ is special
Lagrangian. Thus $F(e,0,0)=0$, and so $Q(e,x,0,0)=0$. 
Following the proof of Proposition \ref{cm2prop2} we
can also show that~$\pd_yQ(e,x,0,0)=\pd_zQ(e,x,0,0)=0$.

Therefore by Taylor expansion of $Q(\hat e,x,y,z)$ about
$\hat e=e$, $y=z=0$ we see that {\it for fixed\/} $x$ {\it in}
$X'$ and small $\md{y},\md{z},d_\sE(\hat e,e)$, we have
\e
Q(\hat e,x,y,z)=O\bigl(\ms{y}+\ms{z}+d_\sE(\hat e,e)\bigr),
\label{cm6eq12}
\e
and more generally for fixed $x$, small $\md{y},\md{z},
d_\sE(\hat e,e)$, and $a,b,c\ge 0$ we have
\e
(\na_x)^a(\pd_y)^b(\pd_z)^cQ(\hat e,x,y,z)=
O\bigl(\md{y}^{\max(0,2-b)}+\md{z}^{\max(0,2-c)}
+d_\sE(\hat e,e)\bigr).
\label{cm6eq13}
\e

To prove \eq{cm6eq10} and \eq{cm6eq11} we have to extend
\eq{cm6eq12} and \eq{cm6eq13} to hold uniformly for $x\in X'$ by
inserting appropriate functions of $x$ as multipliers. Careful
consideration of the asymptotic behaviour of $F$ and $Q$ and
their derivatives near $x_i$ for $i=1,\ldots,n$ shows that
the powers of $\rho$ given in \eq{cm6eq10} and \eq{cm6eq11}
suffice. These powers are independent of $\bs\mu$ as the
inequalities $\mu_i>2$ imply that the terms given dominate
other error terms involving the~$\mu_i$.
\end{proof}

We can also refine the image of $F$ in~$C^0(X')$.

\begin{prop} In the situation above, $F$ maps
\e
F:\ti\E\t \D_\sXp\ra\bigl\{u\in L^p_{k-2,{\bs\mu}-2}
(X'):\ts\int_{X'}u\,\d V_g=0\bigr\},
\label{cm6eq14}
\e
and this is a smooth map of Banach manifolds.
\label{cm6prop3}
\end{prop}

\begin{proof} If $(\hat e,\be,f)\in\ti\E\t\D_{\sst
X'}$ then $\be$ is smooth and compactly-supported and $f\in
L^p_{k,{\bs\mu}}(X')$, so $-\d^*\bigl(\psi^m(\be+\d f)\bigr)\in
L^p_{k-2,{\bs\mu}-2}(X')$. Hence we must show that the $Q$ term
in \eq{cm6eq9} also lies in $L^p_{k-2,{\bs\mu}-2}(X')$. For
$x\in X'$, write
\begin{equation*}
y(x)=(\be\!+\!\d f)(x),\;\>
z(x)=(\na\be\!+\!\na^2f)(x)
\;\>\text{and}\;\>
v(x)=Q\bigl(\hat e,x,y(x),z(x)\bigr).
\end{equation*}
Then we must show that~$v\in L^p_{k-2,{\bs\mu}-2}(X')$.

As $L^p_{k,\bs{\mu}}(X')\subset C^2_{\bs\mu}(X')$ by Theorem
\ref{cm4thm1}, we have $\md{y}=O(\rho^{{\bs\mu}-1})$ and
$\md{z}=O(\rho^{{\bs\mu}-2})$. Equation \eq{cm6eq10} then gives
\begin{equation*}
v=Q\bigl(\hat e,x,y(x),z(x)\bigr)=O(\rho^{2{\bs\mu}-4})+
O(\rho^{2{\bs\mu}-4})+O\bigl(\rho(x)d_\sE(\hat e,e)\bigr).
\end{equation*}
Now $2\mu_i-4>\mu_i-2$ and $1>\mu_i-2$ as $2<\mu_i<3$, so
$v$ decays faster than $\rho^{{\bs\mu}-2}$ near $x_i$, and
it follows that~$v\in L^p_{0,{\bs\mu}-2}(X')$.

For the derivatives of $v$, by the chain rule we have
\begin{align*}
\bmd{\na^jv}\le\, &j!\sum_{\substack{a,b,c\ge 0\\a+b+c\le j}}
\Big\vert(\na_x)^a(\pd_y)^b(\pd_z)^cQ\bigl(\hat e,x,y(x),z(x)\bigr)
\Big\vert\\
&\t\sum_{\substack{m_1,\ldots,m_b,n_1,\ldots,n_c\ge 1\\
a+m_1+\cdots+m_b+n_1+\cdots+n_c=j}}
\prod_{i=1}^b\bmd{\na^{m_i}y(x)}\cdot
\prod_{i=1}^c\bmd{\na^{n_i}z(x)}\,.
\end{align*}
Using \eq{cm6eq13} to estimate $\bmd{(\na_x)^a(\pd_y)^b(\pd_z)^c
Q\bigl(\hat e,x,y(x),z(x)\bigr)}$ and noting that $y\in
L^p_{k-1,{\bs\mu}-1}(X')$ and $z\in L^p_{k-2,{\bs\mu}-2}(X')$,
after some calculations using Theorem \ref{cm4thm1} and
H\"older's inequality we can show that $\md{\na^jv}\in
L^p_{0,{\bs\mu}-2-j}(X')$ for $j=0,\ldots,k-2$, so
that~$v\in L^p_{k-2,{\bs\mu}-2}(X')$.

Therefore $F$ maps $\ti\E\t\D_\sXp\ra
L^p_{k-2,{\bs\mu}-2}(X')$. As in  Proposition \ref{cm6prop1},
each $(\hat e,\be,f)\in\smash{\ti\E}\t\D_\sXp$ defines a
compact $C^1$ Lagrangian $m$-fold $\hat X$ in $M$ with conical
singularities. Regard $\hat X,X$ as $m$-chains in homology. Then
$[\hat X]=[X]\in H_m(M,\Z)$ as $\hat X,X$ are isotopic. So using
\eq{cm6eq4} we see that
\begin{equation*}
\int_{X'}F(\hat e,\be,f)\,\d V_g=\int_{\hat X'}\Im\Om
=[\hat X]\cdot[\Im\Om]=[X]\cdot[\Im\Om]=
\int_{\hat X'}\Im\Om=0,
\end{equation*}
as $\Im\Om$ is closed and $X'$ is special Lagrangian. Thus $F$
maps to the r.h.s.\ of \eq{cm6eq14}, as we have to prove. The
smoothness of $F$ as a map between Banach manifolds easily
follows from the smoothness of $Q$ and general limiting arguments.
\end{proof}

\subsection{The obstruction space}
\label{cm62}

We shall determine the derivative $\d F\vert_{(e,0,0)}$
of $F$ at~$(e,0,0)$.

\begin{prop} There exists a unique linear map $\chi:T_e\ti\E
\ra C^\iy_{\bs 0}(X')$, where ${\bs 0}=(0,\ldots,0)\in\R^n$ and\/
$\chi(y)\equiv 0$ on $K$ for all\/ $y\in T_e\ti\E$, such
that\/ $\d F\vert_{(e,0,0)}:T_e\ti\E\t\H_\sXp\t
L^p_{k,{\bs\mu}}(X')\ra L^p_{k-2,{\bs\mu}-2}(X')$ is given by
\e
\d F\vert_{(e,0,0)}:(y,\be,f)\mapsto
\d^*\bigl(\psi^m(\d[\chi(y)]-\be-\d f)\bigr).
\label{cm6eq15}
\e
\label{cm6prop4}
\end{prop}

\begin{proof} As $F$ is smooth by Proposition \ref{cm6prop3},
$\d F\vert_{(e,0,0)}$ is well-defined. Equation \eq{cm6eq9}
then shows that $\d F\vert_{(e,0,0)}$ maps $(0,\be,f)\mapsto
-\d^*(\psi^m(\be+\d f))$, since \eq{cm6eq10} implies that the
$Q$ term in \eq{cm6eq9} can only have derivative 0 in $\be,f$
at $(0,0)$. This gives the final two terms in~\eq{cm6eq15}.

Let $y\in T_e\ti\E$, and differentiate $\Phi^{\hat e}_\sXp$
w.r.t.\ $\hat e$ in the direction of $y$ at $\hat e=e$. This gives
$\pd_y\Phi^{\hat e}_\sXp\vert_{\hat e=e}$, which is a section
of the vector bundle $(\Phi^e_\sXp)^*(TM)$ over $U_\sXp$.
Now $\d\Phi^e_\sXp$ induces an isomorphism of $TU_\sXp$
and $(\Phi^e_\sXp)^*(TM)$ as vector bundles over $U_\sXp$.
Therefore $v=(\d\Phi^e_\sXp)^*\bigl(\pd_y\Phi^{\hat e}_{\sst
X'}\vert_{\hat e=e}\bigr)$ is a section of $TU_\sXp$, that is,
a vector field on $U_\sXp$, which depends linearly on~$y$.

Differentiating $(\Phi^{\hat e}_\sXp)^*(\Im\Om)$ w.r.t.\
$\hat e$ in the direction of $y$, we find that
\begin{equation*}
\pd_y(\Phi^{\hat e}_\sXp)^*(\Im\Om)\vert_{\hat e=e}=
{\cal L}_v(\Phi^e_\sXp)^*(\Im\Om),
\end{equation*}
where ${\cal L}_v$ is the Lie derivative. But restricting to
$X'\subset U_\sXp$ we have
\begin{equation*}
(\Phi^{\hat e}_\sXp)^*(\Im\Om)\vert_{X'}=F(\hat e,0,0)\,\d V_g,
\end{equation*}
by \eq{cm6eq4}. Combining the last two equations gives
\e
\pd_yF(e,0,0)\,\d V_g=\bigl({\cal L}_v(\Phi^e_\sXp)^*
(\Im\Om)\bigr)\vert_\sXp.
\label{cm6eq16}
\e

Define a 1-form $\al$ on $U_\sXp$ by $\al=v\cdot\hat\om$.
Then from \eq{cm6eq16} and the proof of Proposition \ref{cm2prop2}
we find that
\e
\pd_yF(e,0,0)=\d^*(\psi^m\al\vert_\sXp).
\label{cm6eq17}
\e
Since $(\Phi^{\hat e}_\sXp)^*(\om)=\hat\om$ for all $\hat e
\in\ti\E$, it follows that ${\cal L}_v\hat\om\equiv 0$, and
hence $\al$ is a {\it closed\/} 1-form on $U_\sXp$. Also
$v=\al=0$ on $\pi^*(K)$ as $\Phi_\sXp^{\hat e}\equiv
\Phi_\sXp$ on $\pi^*(K)\subset U_\sXp$, by
Theorem~\ref{cm5thm1}.

Thus $\al\vert_\sXp$ is a closed 1-form on $X'$ which is zero
on $K$. Since $X'$ retracts onto $K$ there exists a unique smooth
function $\chi(y):X'\ra\R$ with $\al\vert_{X'}=\d[\chi(y)]$ and
$\chi(y)\equiv 0$ on $K$. Clearly $\chi(y)$ is linear in $y$, and
\eq{cm6eq17} gives
\begin{equation*}
\d F\vert_{(e,0,0)}\bigl((y,0,0)\bigr)=
\pd_yF(e,0,0)=\d^*\bigl(\psi^m\d[\chi(y)]\bigr).
\end{equation*}
This completes the proof of~\eq{cm6eq15}.

It remains to show that $\chi$ maps $T_e\ti\E\ra C^\iy_0(X')$.
As $\Phi^{\hat e}_\sXp$ satisfies \eq{cm5eq4}, one can show
that $v$ and $\al$ on $\pi^*(S_i)\subset U_\sXp$ are the
pull-backs under $\Phi^e_\sXp$ of a smooth vector field $v'$
and a smooth closed 1-form $\al'$ on $\Up^e_i(B_R)$, where
$\Up_i^*(v')=\pd_y\Up^{\hat e}_i\vert_{\hat e=e}$ and $\al'=v'
\cdot\om$. This implies estimates on the decay of $\al$ and
its derivatives on $S_i$ for $i=1,\ldots,n$, which imply that
$\chi(y)\in C^\iy_{\bs 0}(X')$, as we want.
\end{proof}

To apply the Implicit Mapping Theorem to $F$ in \S\ref{cm63},
we will need to know how close $\d F\vert_{(e,0,0)}$ is to being
injective and surjective. First we show that $\d F\vert_{(e,0,0)}$
is injective on a large subspace of its domain.

\begin{prop} The restriction of\/ $\d F\vert_{(e,0,0)}$ to
$T_e\ti\E\t\d\K_\sXp\t L^p_{k,{\bs\mu}}(X')$ is
injective, where $\d\K_\sXp\le\H_\sXp$
as in Definition~\ref{cm6def1}.
\label{cm6prop5}
\end{prop}

\begin{proof} Let $(y,\d v,f)\in T_e\ti\E\t\d\K_\sXp\t
L^p_{k,{\bs\mu}}(X')$ with $\d F\vert_{(e,0,0)}(y,\d v,f)=0$. Then
\begin{equation*}
\d^*\bigl(\psi^m\d[\chi(y)-v-f]\bigr)=0
\end{equation*}
by \eq{cm6eq15}. Multiplying this equation by $\chi(y)-v-f$ and
integrating over $X'$ by parts, we find
\begin{equation*}
\int_{X'}\psi^m\bms{\d[\chi(y)-v-f]}\d V_g=0.
\end{equation*}
This holds even though $X'$ is noncompact, because of the
asymptotic behaviour of $\chi(y)-v-f$ and its derivatives
near $x_i$, and may be proved rigorously using
\cite[Lem.~2.13]{Joyc7}. Thus~$\d[\chi(y)-v-f]=0$.

Now $(y,\d v,f)$ corresponds to an infinitesimal deformation of
$X$ as a Lagrangian $m$-fold in $M$ with conical singularities,
locally the graph of $\d[\chi(y)-v-f]=0$. As $\d[\chi(y)-v-f]=0$
this infinitesimal deformation is trivial, and so cannot change
the singular points $x_i$ or identifications $\up_i$. Therefore
$y=0$, as $\ti\E$ parametrizes nonequivalent choices of
$x_i,\up_i$ by definition.

Hence $\d(v+f)=0$, so $v+f\equiv c\in\R$. As $f\in C^0_{\bs\mu}(X')$
by Theorem \ref{cm4thm1} we have $f(x)\ra 0$ as $x\ra x_i$ in $X'$.
But $v\equiv a_i^j$ on $S_i^j$ and $\sum_{i,j}a_i^j=0$, by Definition
\ref{cm6def1}. Taking $x\ra x_i$ shows that $a_i^j=c$ for all $i,j$,
and thus $c=0$ as $\sum_{i,j}a_i^j=0$. Hence $v=0$ on $S_i$ for all
$i$, and $v$ is {\it compactly-supported}, so that $[\d v]=0$ in
$H^1_{\rm cs}(X',\R)$. Since the map $\K_\sXp\ra H^1_{\rm cs}
(X',\R)$ given by $v\mapsto[\d v]$ is injective, by Definition
\ref{cm6def1}, we see that $v=0$, and hence $f=0$. Therefore
$\d F\vert_{(e,0,0)}$ is injective on~$T_e\ti\E\t\d\K_\sXp\t
L^p_{k,{\bs\mu}}(X')$.
\end{proof}

Next we in effect measure how close $\d F\vert_{(e,0,0)}$ is
to being surjective.

\begin{prop} In the situation above, the map
$L^p_{k,{\bs\mu}}(X')\ra L^p_{k-2,{\bs\mu}-2}(X')$
given by $f\mapsto\d^*(\psi^m\d f)$ is Fredholm with
cokernel of dimension~$\sum_{i=1}^nN_{\sst\Si_i}(2)$.
\label{cm6prop6}
\end{prop}

\begin{proof} This is just the map $P_{\bs\mu}:L^p_{k,{\bs\mu}}(X')
\ra L^p_{k-2,{\bs\mu}-2}(X')$ of Theorem \ref{cm4thm2}. Thus part
(b) of Theorem \ref{cm4thm2} shows that $P_{\bs\mu}$ is injective,
and then part (a) proves that $P_{\bs\mu}$ is Fredholm with cokernel
of dimension $\sum_{i=1}^nN_{\sst\Si_i}(\mu_i)$. But $N_{\sst\Si_i}
(\mu_i)=N_{\sst\Si_i}(2)$ by \eq{cm3eq10}, as $N_{\sst\Si_i}$ is
upper semicontinuous and discontinuous exactly on ${\cal
D}_{\sst\Si_i}$ by Definition~\ref{cm3def2}.
\end{proof}

Now we can define the {\it obstruction space} in our problem.

\begin{dfn} Proposition \ref{cm6prop3} shows that
\begin{equation*}
\d F\vert_{(e,0,0)}\bigl(T_e\ti\E\t\d\K_\sXp\t
L^p_{k,{\bs\mu}}(X')\bigr)\subseteq
\bigl\{u\in L^p_{k-2,{\bs\mu}-2}(X'):\ts\int_{X'}u\,\d V_g=0\bigr\},
\end{equation*}
and Propositions \ref{cm6prop4} and \ref{cm6prop6} show that this
inclusion is of finite codimension. Choose a finite-dimensional
vector subspace $\O_\sXp$ of smooth, compactly-supported functions
$v$ on $X'$ with $\int_{X'}v\,\d V_g=0$, such that
\e
\begin{split}
\bigl\{u\in L^p_{k-2,{\bs\mu}-2}(X'):\,&\ts\int_{X'}u\,\d V_g=0\bigr\}=\\
&\O_\sXp\op\d F\vert_{(e,0,0)}\bigl(T_e\ti\E\t
\d\K_\sXp\t L^p_{k,{\bs\mu}}(X')\bigr).
\end{split}
\label{cm6eq18}
\e
This is possible as such functions $v$ are dense in the l.h.s.\
of \eq{cm6eq18}. We call $\O_\sXp$ the {\it obstruction space}.
Propositions \ref{cm6prop4}--\ref{cm6prop6} imply that
\e
\begin{split}
\dim\O_\sXp&=
\sum_{i=1}^nN_{\sst\Si_i}(2)-\dim\ti\E-\dim\K_\sXp-1\\
&=\sum_{i=1}^nN_{\sst\Si_i}(2)-n(m^2+2m-1)+\sum_{i=1}^n\dim G_i
-\sum_{i=1}^nb^0(\Si_i)\\
&=\sum_{i=1}^n\bigl(N_{\sst\Si_i}(2)-b^0(\Si_i)
-m^2-2m+1+\dim G_i\bigr)\\
&=\sum_{i=1}^n\sind(C_i),
\end{split}
\label{cm6eq19}
\e
where $\dim\ti\E=\dim\E$ is given in \eq{cm5eq2} and
$\dim\K_\sXp$ in \eq{cm6eq1}, we use \eq{cm3eq5} in the
last line, and $\sind(C_i)\ge 0$ is the {\it stability index}
of Definition~\ref{cm3def3}.
\label{cm6def2}
\end{dfn}

We may interpret \eq{cm6eq19} by saying that each singular point
$x_i$ contributes an obstruction space of dimension $\sind(C_i)$
to deforming $X$ as an SL $m$-fold with conical singularities,
and $\O_\sXp$ is the sum of these obstruction spaces.

\subsection{The main result}
\label{cm63}

We are now ready to prove our main results on the moduli
space $\M_\sX$ of compact SL $m$-folds with
conical singularities. The key tool is the {\it Implicit
Mapping Theorem}. The following version may be proved
from Lang~\cite[Th.~2.1, p.~131]{Lang}.

\begin{thm} Let\/ $Y,Z$ and\/ $T$ be Banach spaces, and\/ $W$
an open neighbourhood of\/ $(0,0)$ in $Y\t Z$. Suppose that the
function $G:W\ra T$ is a smooth map of Banach manifolds with\/
$G(0,0)=0$, and that\/ $\d G_{(0,0)}\vert_Z:Z\ra T$ is an
isomorphism of\/ $Z,T$ as vector and topological spaces. Then
there exist open neighbourhoods $U,V$ of\/ $0$ in $Y$ and\/ $Z$
with\/ $U\t V\subseteq W$ and a smooth map $H:U\ra V$ with\/
$H(0)=0$ such that if\/ $(u,v)\in U\t V$ then $G(u,v)=0$ if
and only if\/~$v=H(u)$.
\label{cm6thm1}
\end{thm}

Here is our first main result, describing $\M_\sX$
near~$X$.

\begin{thm} Suppose $(M,J,\om,\Om)$ is an almost Calabi--Yau
$m$-fold and\/ $X$ a compact SL\/ $m$-fold in $M$ with conical
singularities at\/ $x_1,\ldots,x_n$ and cones $C_1,\ldots,C_n$.
Let\/ $\M_\sX$ be the moduli space of deformations
of\/ $X$ as an SL\/ $m$-fold with conical singularities in $M$,
as in Definition \ref{cm5def2}. Set\/~$X'=X\sm\{x_1,\ldots,x_n\}$.

Then there exist natural finite-dimensional vector spaces
$\I_\sXp$, $\O_\sXp$ such that\/ $\I_\sXp$ is isomorphic to
the image of\/ $H^1_{\rm cs}(X',\R)$ in $H^1(X',\R)$ and\/ 
$\dim\O_\sXp=\sum_{i=1}^n\sind(C_i)$, where $\sind(C_i)$ is
the stability index of Definition \ref{cm3def3}. There exists
an open neighbourhood\/ $U$ of\/ $0$ in $\I_\sXp$, a smooth
map $\Phi:U\ra\O_\sXp$ with\/ $\Phi(0)=0$, and a map
$\Xi:\{u\in U:\Phi(u)=0\}\ra\M_\sX$ with\/ $\Xi(0)=X$ which is
a homeomorphism with an open neighbourhood of\/ $X$ in~$\M_\sX$.
\label{cm6thm2}
\end{thm}

\begin{proof} As $\ti\E$ is an open neighbourhood of $e$
in $\E$, which is an open ball, we can choose a smooth
identification of $\ti\E$ with an open neighbourhood of
0 in $T_e\ti\E$ which identifies $e$ with 0 and induces
the identity map on $T_e\ti\E$. Define
\ea
\begin{split}
Y&=\I_\sXp, \quad
Z=\O_\sXp\t T_e\ti\E\t\K_\sXp\t
L^p_{k,{\bs\mu}}(X'),\\
T&=\{u\in L^p_{k-2,{\bs\mu}-2}(X'):\ts\int_{X'}u\,\d V_g=0\bigr\}
\quad\text{and}
\end{split}
\label{cm6eq20}\\
W&=\bigl\{(\be,\ga,\hat e,v,f)\in Y\t Z:{\hat e}\in\ti\E
\subset T_e\ti\E,\quad (\be+\d v,f)\in\D_\sXp\bigr\}.
\label{cm6eq21}
\ea
Then $0\in Z$ is $(0,e,0,0)$. Choose any norms on the
finite-dimensional spaces $\I_\sXp,\O_\sXp,
T_e\ti\E,\K_\sXp$, and use the usual norms on
$L^p_{k,{\bs\mu}}(X')$ and $T$. Then $Y,Z,T$ are Banach spaces,
and $W$ is an open neighbourhood of $(0,0)$ in $Y\t Z$, as in
Theorem~\ref{cm6thm1}.

Define a map $G:W\ra T$ by $G(\be,\ga,\hat e,v,f)=\ga+F(\hat e,
\be+\d v,f)$. This is a smooth map of Banach manifolds, by
Proposition \ref{cm6prop3}, and $G(0,0)=G(0,0,e,0,0)=0$
as $F(e,0,0)=0$. The map $\d G_{(0,0)}\vert_Z$ is given by
\e
\d G_{(0,0)}\vert_Z:(\ga,y,v,f)\mapsto\ga+\d F_{(e,0,0)}(y,\d v,f).
\label{cm6eq22}
\e
Now Proposition \ref{cm6prop5} proves that $(y,v,f)\mapsto\d
F_{(e,0,0)}(y,\d v,f)$ is an injective map on $T_e\ti\E
\t\K_\sXp\t L^p_{k,{\bs\mu}}(X')$. Also \eq{cm6eq18} implies
that $\O_\sXp$ intersects the image of $\d F_{(e,0,0)}$
only in 0. Therefore $\d G_{(0,0)}:Z\ra T$ is {\it injective}.

But \eq{cm6eq18} shows that $\d G_{(0,0)}$ is {\it surjective}.
Thus $\d G_{(0,0)}$ is an isomorphism of $Z,T$ as vector spaces.
Since $\d G_{(0,0)}$ is continuous, it is an isomorphism of
$Z,T$ as topological spaces by the Open Mapping Theorem. Hence
the hypotheses of Theorem \ref{cm6thm1} hold, and the theorem
gives open neighbourhoods $U$ of 0 in $\I_\sXp$ and
$V$ of 0 in $Z$ and a smooth map $H:U\ra V\subset Z$ with~$H(0)=0$.

Since $(\be,v)\mapsto\be+\d v$ is a homeomorphism
$\I_\sXp\t\K_\sXp\ra\H_\sXp$ by \eq{cm6eq2}, we see from
\eq{cm6eq21} that the map
\begin{equation*}
\bigl\{(\be,0,\hat e,v,f)\in W\bigr\}\ra\ti\E\t\D_\sXp
\quad\text{given by}\quad
(\be,0,\hat e,v,f)\mapsto(\hat e,\be+\d v,f)
\end{equation*}
is a homeomorphism. Applying Proposition \ref{cm6prop1} we see that
\begin{itemize}
\item[(a)] The map $\bigl\{(\be,0,\hat e,v,f)\in G^{-1}(0)\subset W
\bigr\}\ra\V_\sX$ given by $(\be,0,\hat e,v,f)\mapsto\Psi(
\hat e,\be+\d v,f)$ is a homeomorphism taking~$(0,0,e,0,0)\mapsto X$.
\end{itemize}

Define $\Phi:U\ra\O_\sXp$, $H_1:U\ra T_e\ti\E$,
$H_2:U\ra\K_\sXp$ and $H_3:U\ra L^p_{k,{\bs\mu}}(X')$ by
$H(u)=\bigl(\Phi(u),H_1(u),H_2(u),H_3(u)\bigr)\in V\subset Z$.
Then $\Phi,H_1,H_2,H_3$ are smooth as $H$ is smooth, and
$\Phi(0)=0$, $H_j(0)=0$ as $H(0)=0$. By Theorem \ref{cm6thm1},
if $(u,v)\in U\t V$ then $G(u,v)=0$ if and only if $v=H(u)$.
That is:
\begin{itemize}
\item[(b)] if $(\be,\ga,\hat e,v,f)\in U\t V\subseteq W$ then
$G(\be,\ga,\hat e,v,f)=0$ if and only if $\ga=\Phi(\be)$,
$\hat e=H_1(\be)$, $v=H_2(\be)$ and~$f=H_3(\be)$.
\end{itemize}

Combining (a), (b) proves that $\Xi:\{u\in U:\Phi(u)=0\}\ra
\V_\sX$ given by $\Xi(u)=\Psi\bigl(H_1(u),u+\d H_2(u),
H_3(u)\bigr)$ is a homeomorphism from $U$ to an open
neighbourhood of $X$ in $\V_\sX$ with $\Xi(0)=X$.
This completes the proof.
\end{proof}

Here are two simple corollaries of Theorem \ref{cm6thm2}.
Firstly, if $X$ has {\it stable singularities} in the sense
of Definition \ref{cm3def4} then $\sind(C_i)=0$, so
$\dim\O_\sXp=0$, and $\M_\sX$ is
locally homeomorphic to $\I_\sXp$. Thus
$\M_\sX$ is a {\it manifold\/} near~$X$.

But all SL $m$-folds $\hat X\in\M_\sX$ have the same cones
$C_i$, so all $\hat X\in\M_\sX$ have stable singularities,
and $\M_\sX$ is a manifold everywhere. The maps $\Xi$ of
Theorem \ref{cm6thm2} provide coordinate charts on $\M_\sX$.
It is easy to see that the transition maps are smooth (this
follows for instance from Theorem \ref{cm6thm3} below), so
$\M_\sX$ is a {\it smooth manifold}. This gives:

\begin{cor} Suppose $(M,J,\om,\Om)$ is an almost Calabi--Yau
$m$-fold and\/ $X$ a compact SL\/ $m$-fold in $M$ with stable
conical singularities, and let\/ $\M_\sX$ and\/
$\I_\sXp$ be as in Theorem \ref{cm6thm2}. Then
$\M_\sX$ is a smooth manifold of
dimension~$\dim\I_\sXp$.
\label{cm6cor1}
\end{cor}

Here is another simple condition for $\M_\sX$ to be
a manifold near~$X$.

\begin{dfn} Let $(M,J,\om,\Om)$ be an almost Calabi--Yau
$m$-fold and $X$ a compact SL $m$-fold in $M$ with conical
singularities, and let $\I_\sXp,\O_\sXp,U$
and $\Phi$ be as in Theorem \ref{cm6thm2}. We call $X$ {\it
transverse} if the linear map $\d\Phi\vert_0:\I_\sXp
\ra\O_\sXp$ is surjective. It is not difficult to
see that this definition is independent of the choices made in
defining $\I_\sXp,\O_\sXp,U$ and~$\Phi$.
\label{cm6def3}
\end{dfn}

If $X$ is transverse then $\{u\in U:\Phi(u)=0\}$ is a manifold
near 0, so we prove:

\begin{cor} Suppose $(M,J,\om,\Om)$ is an almost Calabi--Yau $m$-fold
and\/ $X$ a transverse compact SL\/ $m$-fold in $M$ with conical
singularities, and let\/ $\M_\sX,\I_\sXp$ and\/ $\O_\sXp$ be as
in Theorem \ref{cm6thm2}. Then $\M_\sX$ is near $X$ a smooth manifold
of dimension~$\dim\I_\sXp-\dim\O_\sXp$.
\label{cm6cor2}
\end{cor}

\subsection{Naturality of $\I_\sXp,\O_\sXp,\Phi$ and $\Xi$}
\label{cm64}

In the course of proving Theorem \ref{cm6thm2} we made a
considerable number of {\it arbitrary choices} in \S\ref{cm43},
\S\ref{cm5} and \S\ref{cm6}, including $\Up_i,\ze,U_\sXp,
\Phi_\sXp,\E,\ti\E,\Up_i^{\hat e},\allowbreak
\Phi_\sXp^{\hat e},\allowbreak
\H_\sXp,\allowbreak
\I_\sXp,\allowbreak
\O_\sXp$ and $U$. We now consider to what extent the
final result depends on these choices, in particular the vector
spaces $\I_\sXp,\O_\sXp$ and maps~$\Phi,\Xi$.

Now $\I_\sXp$ is naturally isomorphic to the image of
$H^1_{\rm cs}(X',\R)$ in $H^1(X',\R)$ by \S\ref{cm61}. Thus {\it as
a vector space} $\I_\sXp$ depends only on $X'$, though
{\it as a vector space of\/ $1$-forms} it depends on an arbitrary
choice. Let us {\it identify} $\I_\sXp$ with the image of
$H^1_{\rm cs}(X',\R)$ in $H^1(X',\R)$, so that $\I_\sXp$
is independent of choices. 

Then $\Xi$ maps $\Phi^{-1}(0)\subset\I_\sXp\subseteq
H^1(X',\R)$ to $\M_\sX$, as a local homeomorphism.
In the next theorem we shall construct an {\it inverse} $\Th$
for $\Xi$, defined near $X$ in $\M_\sX$ and mapping
into $H^1(X',\R)$, which is independent of all arbitrary choices.
This proves that both $\Xi$ and its domain $\{u\in U:\Phi(u)=0\}
\subset\I_\sXp$ are independent of arbitrary choices
near 0 in~$\I_\sXp$.

In \S\ref{cm65} we will explain an alternative construction of
$\O_\sXp$ {\it as a vector space} which is independent of choices.
The author does not know to what extent $\Phi$ is natural where
it is nonzero, but this does not seem a very important question.
The theorem is based on the construction of {\it natural
coordinates} on moduli spaces $\M_\sX$ of compact, nonsingular
SL $m$-folds, which is described by Hitchin \cite[\S 4]{Hitc}
and the author~\cite[\S 9.4]{Joyc6}.

\begin{thm} Let\/ $(M,J,\om,\Om),X,X',\M_\sX,U,\Xi$
and\/ $\Phi$ be as in Theorem \ref{cm6thm2}, and let\/ $V$ be
a path-connected, simply-connected open neighbourhood of\/
$X$ in $\M_\sX$. Then there exists a natural, continuous map
$\Th:V\ra H^1(X',\R)$ depending only on $M,\om,X$ and\/ $V$,
such that\/ $\Th,\Xi$ are inverse maps on the connected
component of\/ $V\cap\Xi(U)$ containing~$X$.
\label{cm6thm3}
\end{thm}

\begin{proof} Let $\hat X\in V$. As $V$ is path-connected and
simply-connected there is a unique isotopy class of continuous
paths $\ga:[0,1]\ra V$ with $\ga(0)=X$ and $\ga(1)=\hat X$.
This determines a unique isotopy class of continuous maps
$\Pi:[0,1]\t X'\ra M$ with $\Pi\bigl(\{0\}\t X'\bigr)=X'$
and $\Pi\bigl(\{1\}\t X'\bigr)=\hat X'$. Let $\Pi$ be
a smooth map in this isotopy class. Then $\Pi^*(\om)$ is
a closed 2-form on $[0,1]\t X'$ vanishing on $\{0,1\}\t X'$,
since $X',\hat X'$ are Lagrangian.

Thus $[\Pi^*(\om)]$ defines a class in $H^2\bigl([0,1]\t X';
\{0,1\}\t X',\R\bigr)$, the {\it relative de Rham cohomology
group}, which depends only on $M,\om,V,X$ and $\hat X$. Define
$\Th(\hat X)$ to be the class in $H^1(X',\R)$ corresponding to
$[\Pi^*(\om)]$ under the natural isomorphism $H^1(X',\R)\cong
H^2\bigl([0,1]\t X';\{0,1\}\t X',\R\bigr)$. Then $\Th:V\ra
H^1(X',\R)$ depends only on $M,\om,X$ and $V$, and is clearly
continuous.

We must show that $\Th,\Xi$ are inverse near $X$. Let $\hat X$
lie in the connected component of $V\cap\Xi(U)$ containing $X$.
From \S\ref{cm61}--\S\ref{cm63} we find that $\hat X'=\Phi_{\sst
X'}^{\hat e}\bigl(\Ga(\be+\d f)\bigr)$ for some $(\hat e,\be,f)
\in\ti\E\t\D_\sXp$, and that $[\be]\in H^1(X',\R)$
lies in $U\subset\I_\sXp\subset H^1(X',\R)$ with
$\Phi\bigl([\be]\bigr)=0$ and~$\Xi\bigl([\be]\bigr)=\hat X$.

Now $\Phi_\sXp^{\hat e}\equiv\Phi_\sXp$ on $\pi^*(K)$. Assuming
the fibres of $\pi:U_\sXp\ra X'$ are convex for simplicity, we
may take $\Pi\vert_{[0,1]\t K}$ above to be $\Pi(t,x)=\Phi_\sXp
\bigl(t(\be+\d f)\vert_x\bigr)$. This has the correct isotopy
class as $X,\hat X$ lie in the same component of $V\cap\Xi(U)$.
Since $\Phi_\sXp^*(\om)=\hat\om$, a short calculation then
shows that $\Pi^*(\om)=(\be+\d f)\w\d t$ on $[0,1]\t K$. As $X'$
retracts onto $K$, we find that $\Th(\hat X)$ is $[\be+\d f]=
[\be]\in H^1(X',\R)$. But $\Xi\bigl([\be]\bigr)=\hat X$, so
$\Th,\Xi$ are inverse.
\end{proof}

The theorem implies that the topology on $\M_\sX$ is locally
induced from the Euclidean topology on $H^1(X',\R)$ via $\Th$.
This gives another way of seeing the naturality of the topology
on~$\M_\sX$.

\subsection{Another way of thinking about $\I_\sXp,\O_\sXp$}
\label{cm65}

In \S\ref{cm23} we saw that for a compact, nonsingular SL $m$-fold
$N$ in an almost Calabi--Yau $m$-fold $M$, the infinitesimal
deformations correspond to 1-forms $\al$ on $N$ with $\d\al=
\d^*(\psi^m\al)=0$, which form a vector space naturally isomorphic
to $H^1(N,\R)$. To extend this to SL $m$-folds $X$ with conical
singularities $x_1,\ldots,x_n$ with rates $\mu_1,\ldots,\mu_n$,
we need to regard $\al$ as a 1-form on $X'$ with {\it asymptotic
conditions} on $\al$ and its derivatives.

We saw in Theorem \ref{cm4thm3} that the most natural asymptotic
condition on $\al$ from the point of view of {\it Hodge theory}
is $\md{\na^k\al}=O(\rho^{-1-k})$ for all $k\ge 0$. The vector
space $Y_\sXp$ of such $\al$ is isomorphic to $H^1(X',\R)$.
Consider for the moment only deformations of $X$ that fix the
$x_i$ and $\up_i$. Then the most natural asymptotic condition
on $\al$ for the {\it deformation theory} of $X$ is
$\md{\na^k\al}=O(\rho^{{\bs\mu}-1-k})$ for all~$k\ge 0$.

Clearly if $\md{\na^k\al}=O(\rho^{{\bs\mu}-1-k})$ then
$\md{\na^k\al}=O(\rho^{-1-k})$. So define
\begin{equation*}
Z_\sXp=\bigl\{\al\in Y_\sXp:\text{$\md{\na^k\al}
=O(\rho^{{\bs\mu}-1-k})$ for all $k\ge 0$}\bigr\}.
\end{equation*}
This is an obvious candidate for the infinitesimal deformations
of $X$ which fix the $x_i,\up_i$. Therefore we ask: {\it how big
a subspace of\/ $Y_\sXp\cong H^1(X',\R)$ is}~$Z_\sXp$?

First note that if the image of $[\al]\in H^1(X',\R)$ under
the map $H^1(X',\R)\ra\bigoplus_{i=1}^nH^1(\Si_i,\R)$ of
\eq{cm4eq4} is nonzero, then one can easily see from the
proof of Theorem \ref{cm4thm3} in \cite[\S 2.5]{Joyc7} that
$\al$ decays exactly at rate $O(\rho^{-1})$ near some $x_i$,
and thus $\al\notin Z_\sXp$. Hence $Z_\sXp$
corresponds to a subspace of the kernel of $H^1(X',\R)\ra
\bigoplus_{i=1}^nH^1(\Si_i,\R)$, that is, to a subspace
of the image of $H^1_{\rm cs}(X',\R)\ra H^1(X',\R)$ in
\eq{cm4eq4}, which is isomorphic to~$\I_\sXp$.

Define $G_i$ to be the space of {\it germs of smooth\/ $1$-forms
on $X'$ near} $x_i$, that is, smooth 1-forms $\xi$ defined on
$U_i\sm\{x_i\}$ for some small open neighbourhood $U_i$ of $x_i$ in
$X$, where two such 1-forms are equivalent if they agree on the
intersection of their domains. For $i=1,\ldots,n$ define
\begin{equation*}
\O_i=\frac{\bigl\{\xi\in G_i:\text{$\xi$ is exact,
$\d^*(\psi^m\xi)=0$, $\md{\na^k\xi}=O(\rho^{-1-k})$
for all $k\ge 0$}\bigr\}}{
\bigl\{\xi\in G_i:\text{$\xi$ is exact, $\d^*(\psi^m\xi)=0$,
$\md{\na^k\xi}=O(\rho^{\mu_i-1-k})$ for all $k\ge 0$}\bigr\}}\,.
\end{equation*}

Then one can show that $\O_i$ is a vector space of
dimension $N_{\sst\Si_i}(2)-b^0(\Si_i)$, an {\it obstruction
space}. Each $\xi$ in the subspace of $Y_\sXp$ corresponding
to $\I_\sXp$ has a natural projection to $\O_i$
for $i=1,\ldots,n$, and $\xi\in Z_\sXp$ if and only if all
of these projections are zero. Thus the infinitesimal deformation
space $Z_\sXp$ is the kernel of a linear map $\I_{\sst
X'}\ra\bigoplus_{i=1}^n\O_i$, and each obstruction space
$\O_i$ depends only on the germ of $X$ at $x_i$, and
essentially only on the cone~$C_i$.

In fact $\bigoplus_{i=1}^n\O_i$ does not correspond exactly
to the obstruction space $\O_\sXp$ of \S\ref{cm62}, as
$\O_\sXp$ is the obstructions to deformations which can
vary $x_i,\up_i$. Each $\O_i$ contains a vector subspace
${\cal P}_i$ isomorphic to $T_{(x_i,\up_i)}\E_i$, corresponding
to infinitesimal deformations $\xi$ which vary $x_i,\up_i$. It can be
shown that there is a natural isomorphism $\O_\sXp\cong
\bigoplus_{i=1}^n\O_i/{\cal P}_i$. The corresponding linear
map $\I_\sXp\ra\O_\sXp$ is $\d\Phi\vert_0$,
in the notation of~\S\ref{cm63}.

This way of thinking about the infinitesimal deformation and
obstruction spaces $\I_\sXp,\O_\sXp$ has
the advantages of being closer to McLean's method, and of
presenting $\O_\sXp$ as a direct sum of contributions
from each singular point $x_i$, in a way that was implicit in
\eq{cm6eq19} but was not brought out in \S\ref{cm62}. However, the
author did not find it helpful in actually writing down a proof.

\section{Extension to families $\bigl\{(M,J^s,\om^s,\Om^s):s\in\F\bigr\}$}
\label{cm7}

We now extend the material of \S\ref{cm5} and \S\ref{cm6}
from a single almost Calabi--Yau $m$-fold $(M,J,\om,\Om)$
to a {\it smooth family of deformations} $\bigl\{(M,J^s,
\om^s,\Om^s):s\in\F\bigr\}$ of $(M,J,\om,\Om)$, as in
Definition \ref{cm2def6}. The basic idea is that we consider
deformations $\hat X$ of a compact SL $m$-fold $X$ in
$(M,J,\om,\Om)$ with conical singularities, not just in
$(M,J,\om,\Om)$ but in $(M,J^s,\om^s,\Om^s)$ for~$s\in\F$.

We collect these deformations $(s,\hat X)$ into a big moduli
space $\M_\sX^\sF$ with a natural topology and a continuous
projection $\pi^\sF:\M_\sX^\sF\ra\F$, generalizing \S\ref{cm5}.
Then we show that $\M_\sX^\sF$ is homeomorphic near $(0,X)$ to
the zeroes of a smooth map $\Phi^\sF:\F\t\I_\sXp\ra\O_\sXp$
between finite-dimensional spaces, generalizing~\S\ref{cm6}.

\subsection{Moduli spaces of SL $m$-folds in families
$(M,J^s,\om^s,\Om^s)$}
\label{cm71}

We first explain how to extend \S\ref{cm5} to {\it families}
$\bigl\{(M,J^s,\om^s,\Om^s):s\in\F\bigr\}$ of almost Calabi--Yau
$m$-folds, as in Definition \ref{cm2def6}. In fact this is not
very much work, as we are {\it already} dealing with families
$\E$ of choices of $x_i,\up_i$, so we simply have to enlarge
these families to include $\F$, and make appropriate changes.
Consider the following situation.

\begin{dfn} Let $(M,J,\om,\Om)$ be an almost Calabi--Yau
$m$-fold and $X$ a compact SL $m$-fold in $M$ with conical
singularities at $x_1,\ldots,x_n$ with identifications
$\up_i:\C^m\ra T_{x_i}M$, cones $C_1,\ldots,C_n$ and rates
$\mu_i$. Suppose $\bigl\{(M,J^s,\om^s,\Om^s):s\in\F\bigr\}$
is a smooth family of deformations of $(M,J,\om,\Om)$, where
$\F\subset\R^d$ is the {\it base space}, such that
$\iota_*(\ga)\cdot[\om^s]=0$ for all $\ga\in H_2(X,\R)$
and $s\in{\cal F}$ and $[X]\cdot[\Im\Om^s]=0$ for all
$s\in{\cal F}$. Here $\iota:X\ra M$ is the inclusion,
$\iota_*:H_2(X,\R)\ra H_2(M,\R)$ the induced map,
$[\om^s]\in H^2(M,\R)$, $[X]\in H_m(M,\R)$
and~$[\Im\Om^s]\in H^m(M,\R)$.
\label{cm7def1}
\end{dfn}

The point of this definition is that $\iota_*(\ga)\cdot[\om^s]=0$
for all $\ga$ and $[X]\cdot[\Im\Om^s]=0$ are {\it necessary
conditions} for there to exist an SL $m$-fold $\hat X$ in
$(M,J^s,\om^s,\Om^s)$ with conical singularities, isotopic
to $X$ in $M$. For if $\hat\iota:\hat X\ra\R$ is the inclusion
then by isotopy $\hat\iota_*(\ga)=\iota_*(\ga)$ under the
natural isomorphism $H_2(\hat X,\R)\cong H_2(X,\R)$ and
$[\hat X]=[X]$. But clearly $\hat\iota_*(\ga)\cdot[\om^s]=0$
for all $\ga\in H_2(\hat X,\R)$ and $[\hat X]\cdot[\Im\Om^s]=0$,
since~$\om^s\vert_{\hat X'}\equiv\Im\Om^s\vert_{\hat X'}\equiv 0$.

We have written these conditions in an odd way. In effect
$\iota_*(\ga)\cdot[\om^s]=0$ for all $\ga$ and $[X]\cdot[\Im\Om^s]
=0$ simply mean that $[\om^s\vert_X]=[\Im\Om^s\vert_X]=0$ in
$H^*(X,\R)$. However, we have not defined the de Rham cohomology
$H^*(X,\R)$ of the singular manifold $X$, so this does not make
sense. The conditions $[\om^s\vert_{X'}]=[\Im\Om^s\vert_{X'}]=0$
in $H^*(X',\R)$ do make sense, but are not strong enough.

Here are the analogues of Definition \ref{cm5def1} and
Theorems \ref{cm5thm1} and~\ref{cm5thm2}.

\begin{dfn} In the situation of Definition \ref{cm7def1}, for
$s\in\F$ define $\psi^s:M\ra(0,\iy)$ as in \eq{cm2eq3}, but
using $\om^s,\Om^s$. Extending \eq{cm5eq1}, define
\e
\begin{split}
P^\sF=\bigl\{(s,x,\up):\,&\text{$s\in\F$, $x\in M$,
$\up:\C^m\ra T_xM$ is a real isomorphism,}\\
&\up^*(\om^s)=\om',\quad \up^*(\Om^s)=\psi^s(x)^m\Om'\bigr\},
\end{split}
\label{cm7eq1}
\e
where $\om',\Om'$ are as in \eq{cm2eq1}. Define
$\pi^\sF:P^\sF\ra\F$ by $\pi^\sF:(s,x,\up)\mapsto s$.
Define a free $\SU(m)$-action on $P^\sF$ by $B:(s,x,\up)
\mapsto(s,x,\up\circ B^{-1})$. Then $P^\sF$ is a principal
$\SU(m)$-bundle over~$\F\t M$.

Let $G_i$ be the Lie subgroup of $\SU(m)$ preserving $C_i$.
Let $0\in\F'\subseteq\F$ and $\up_i^s:\C^m\ra T_{x_i}M$ for
$i=1,\ldots,n$ and $s\in\F'$ be as in Theorem \ref{cm4thm7}.
Then $(s,x_i,\up_i^s)\in P^\sF$ for $i=1,\ldots,n$ and
$s\in\F'$. Let $\E_i,\E$ be as in Definition~\ref{cm5def1}.

For $i=1,\ldots,n$ let $\E_i^\sFp$ be a submanifold of dimension
$\dim P^\sF-\dim G_i$ in $(\pi^\sF)^*(\F')\subseteq P^\sF$ such
that $\pi^\sF:\E_i^\sFp\ra\F'$ is a submersion, $(\pi^\sF)^{-1}(s)$
is a small ball containing $(s,x_i,\up_i^s)$ for $s\in\F'$ which
is transverse to the orbits of $G_i$, and $(\pi^\sF)^{-1}(0)=\{0\}
\t\E_i$. Making $\F'$ smaller if necessary, such $\E_i^\sFp$ exist.
Define
\e
\E^\sFp=\bigl\{(s,\hat x_1,\hat\up_1,\ldots,\hat x_n,\hat\up_n):
\text{$(s,\hat x_i,\hat\up_i)\in\E_i^\sFp$ for $i=1,\ldots,n$}\bigr\}.
\label{cm7eq2}
\e
Write a general element of $\E^\sFp$ as $(s,\hat e)$ for $s\in\F'$
and $\hat e=(\hat x_1,\hat\up_1,\ldots,\hat x_n,\hat\up_n)$ as in
\S\ref{cm51}, and let $e^s=(x_1,\up_1^s,\ldots,x_n,\up_n^s)$,
so that $(s,e^s)\in\E^\sFp$ for all $s\in\F'$. Define
$\pi^\sF:\E^\sFp\ra\F'$ by $\pi^\sF:(s,\hat x_1,\ldots,
\hat\up_n)\mapsto s$. Then~$(\pi^\sF)^{-1}(0)=\{0\}\t\E$.

This $\E^\sFp$ is a family of $(s,\hat x_i,\hat\up_i)$ such
that $\hat x_i,\hat\up_i$ are close to $x_i,\up_i$, and are
valid alternative choices of $x_i,\up_i$ in $(M,J^s,\om^s,\Om^s)$,
noting that $\hat\up_i:\C^m\ra T_{\hat x_i}M$ has to be compatible
with $\om^s,\Om^s$ as in \S\ref{cm33}. Each $G_1\t\cdots\t G_n$
equivalence class of choices of $s,\hat x_i,\hat\up_i$ close to
$0,x_i,\up_i$ is represented {\it exactly once} in~$\E^\sFp$.
\label{cm7def2}
\end{dfn}

\begin{thm} In the situation above, use the notation of
Theorem \ref{cm4thm5}, let\/ $U_\sXp,\Phi_\sXp$ be as in
Theorem \ref{cm4thm6}, let\/ $0\in\F'\subseteq\F$ and\/
$\up_i^s,\Up_i^s,\Phi_\sXp^s$ for $s\in\F'$ be as in
Theorem \ref{cm4thm7}, and let\/ $\ti\E,\Up_i^{\hat e}$
and\/ $\Phi_\sXp^{\hat e}$ be as in Theorem~\ref{cm5thm1}.

Then making $\F'$ smaller if necessary, there exists a
connected open subset\/ $\ti\E^\sFp\subseteq\E^\sFp$
with\/ $(s,e^s)\in\ti\E^\sFp$ for all\/ $s\in\F'$ and\/
$(\pi^\sF)^{-1}(0)\cap\ti\E^\sFp=\{0\}\t\ti\E$, and for all\/
$(s,\hat e)=(s,\hat x_1,\hat\up_1,\ldots,\hat x_n,\hat\up_n)$
in $\ti\E^\sFp$ there exist
\begin{itemize}
\item[{\rm(a)}] embeddings $\Up_i^{s,\hat e}:B_R\ra M$ for
$i=1,\ldots,n$ with
\end{itemize}
\e
\Up_i^{s,e^s}=\Up_i^s,\quad
\Up_i^{s,\hat e}(0)=\hat x_i,\quad
\d\Up_i^{s,\hat e}\vert_0=\hat\up_i\quad\text{and}\quad
(\Up_i^{s,\hat e})^*(\om^s)=\om',
\label{cm7eq3}
\e
\begin{itemize}
\item[{\rm(b)}] an embedding $\Phi^{s,\hat e}_\sXp:U_\sXp
\ra M$ with\/ $\Phi_\sXp^{s,e^s}=\Phi_\sXp^s$ and\/
$(\Phi^{s,\hat e}_\sXp)^*(\om^s)=\hat\om$, such that\/
$\Phi_\sXp^{s,\hat e}\equiv\Phi_\sXp^s$ on $\pi^*(K)\subset U_\sXp$,
\end{itemize}
all depending smoothly on $(s,\hat e)\in\ti\E^\sFp$, with\/
$\Up_i^{0,\hat e}=\Up_i^{\hat e}$ and\/ $\Phi_\sXp^{0,\hat e}
=\Phi_\sXp^{\hat e}$ for all\/ $\hat e\in\ti\E$ and
\e
\Phi_\sXp^{s,\hat e}\circ\d(\Up_i\circ\phi_i)(\si,r,\tau,u)\equiv
\Up_i^{s,\hat e}\circ\Phi_{\sst C_i}\bigl(\si,r,\tau+\eta_i^1(\si,r),
u+\eta_i^2(\si,r)\bigr)
\label{cm7eq4}
\e
for all\/ $(s,\hat e)\in\ti\E^\sFp$, $1\!\le\!i\!\le\!n$
and\/ $(\si,r,\tau,u)\!\in\!T^*\bigl(\Si_i\t(0,R')\bigr)$
with\/~$\bmd{(\tau,u)}\!<\!\ze r$.
\label{cm7thm1}
\end{thm}

\begin{thm} In the situation above, let\/ $(s,\hat e)=(s,\hat x_1,
\hat\up_1,\ldots,\hat x_n,\hat\up_n)\in\ti\E^\sFp$, and suppose
$\hat X$ is a compact SL\/ $m$-fold in $(M,J^s,\om^s,\Om^s)$ with
conical singularities at\/ $\hat x_1,\ldots,\hat x_n$, with
identifications $\hat\up_i$, cones $C_i$ and rates $\mu_i$. Then
if\/ $(s,\hat e),(0,e)$ are sufficiently close in $\ti\E^\sFp$
and\/ $X',\hat X'$ are sufficiently close as submanifolds
in a $C^1$ sense away from $x_1,\ldots,x_n$, there exists a
closed\/ $1$-form $\al$ on $X'$ such that the graph\/ $\Ga(\al)$
lies in $U_\sXp\subset T^*X'$, and\/ $\hat X'=\Phi^{s,\hat e}_\sXp
\bigl(\Ga(\al)\bigr)$. Furthermore we may write $\al=\be+\d f$,
where $\be$ is a closed\/ $1$-form supported in $K$ and\/~$f\in
C^\iy_{\bs\mu}(X')$.
\label{cm7thm2}
\end{thm}

The proofs are straightforward modifications of \S\ref{cm51},
replacing Theorem \ref{cm4thm6} by Theorem \ref{cm4thm7}. Here
is the analogue of Definitions \ref{cm5def2} and~\ref{cm5def3}.

\begin{dfn} In the situation above, define the {\it moduli space
$\M_\sX^\sF$ of deformations of\/ $X$ in the family} $\F$ to be
the set of pairs $(s,\hat X)$ such that
\begin{itemize}
\setlength{\parsep}{0pt}
\setlength{\itemsep}{0pt}
\item[(i)] $s\in\F$ and $\hat X$ is a compact SL $m$-fold in
$(M,J^s,\om^s,\Om^s)$ with conical singularities at $\hat x_1,
\ldots,\hat x_n$ with cones $C_1,\ldots,C_n$, for some~$\hat x_i$.
\item[(ii)] There exists a homeomorphism $\hat\iota:X\ra\hat X$
with $\hat\iota(x_i)=\hat x_i$ for $i=1,\ldots,n$ such that
$\hat\iota\vert_{X'}:X'\ra\hat X'$ is a diffeomorphism and
$\hat\iota$ and $\iota$ are isotopic as continuous maps
$X\ra M$, where $\iota:X\ra M$ is the inclusion.
\end{itemize}
Define $\pi^\sF:\M_\sX^\sF\ra\F$ by $\pi^\sF:(s,\hat X)\mapsto s$.
Let $\V_{\sst 0,X}^\sF$ be the subset of $(s,\hat X)\in\M_\sX^\sF$
such that for some $(s,\hat e)\in\ti\E^\sFp$ and some 1-form $\al$
on $X'$ whose graph $\Ga(\al)$ lies in $U_\sXp\subset T^*X'$ we
have $\hat X'=\Phi^{s,\hat e}_\sXp\bigl(\Ga(\al)\bigr)$, as in
Theorem~\ref{cm7thm2}.

This gives a 1-1 correspondence between $\V_{\smash{\sst 0,X}}^\sF$
and a set of triples $(s,\hat e,\al)$ for $(s,\hat e)\in\ti\E^\sFp$
and $\al$ a smooth 1-form on $X'$ with prescribed decay. Also
$(\pi^\sF)^{-1}(0)\cap\V_{\sst 0,X}^\sF=\{0\}\t\V_\sX$, where
$\V_\sX$ is as in Definition \ref{cm5def2}, and the triples
$(0,\hat e,\al)$ for $(\pi^\sF)^{-1}(0)\cap\V_{\sst 0,X}^\sF$
agree with the pairs $(\hat e,\al)$ for $\V_\sX$ in
Definition~\ref{cm5def2}.

Use this 1-1 correspondence to define a topology on $\V_{\sst
0,X}^\sF$, using the natural topology on $\smash{\ti\E^\sFp}$
and either the $C^1_{{\bs\mu}-1}$ or the $C^\iy_{{\bs\mu}-1}$
topology on $\al$, defined as in \S\ref{cm52}. The analogue
of Proposition \ref{cm5prop} shows that these yield the same
topology on $\V_{\sst 0,X}^\sF$, which is also independent
of choice of rates~$\mu_i$.

For each $(\ti s,\ti X)\in\M_\sX^\sF$ we can regard
$\bigl\{(M,J^s,\om^s,\Om^s):s\in\F\bigr\}$ as a family of
deformations of $(M,J^{\ti s},\om^{\ti s},\Om^{\ti s})$
rather than of $(M,J^0,\om^0,\Om^0)$, and we can redo the
whole of this section replacing $0\in\F$ by $\ti s\in F$
and $X$ by $\ti X$. In this way we define a subset $\V_{
\smash{\sst\ti s,\ti X}}^\sF$ of $\M_\sX^\sF$ containing
$(\ti s,\ti X)$ with a 1-1 correspondence between $\V_{
\smash{\sst\ti s,\ti X}}^\sF$ and a set of triples $(s,\hat
e,\al)$, and a topology on~$\V_{\smash{\sst\ti s,\ti X}}^\sF$.

One can show that the topologies on different neighbourhoods
$\V_{\smash{\sst\ti s,\ti X}}^\sF$ agree on the overlaps, and
that the overlaps are open in each. Piecing the topologies
together therefore defines a unique topology on $\M_\sX^\sF$.
In this topology $\pi^\sF:\M_\sX^\sF\ra\F$ is continuous,
and $\V_{\sst 0,X}^\sF$ is an open neighbourhood of~$(0,X)$.

Note that $(\pi^\sF)^{-1}(0)\subset\M_\sX^\sF$ is just
$\{0\}\t\M_\sX$ in the notation of \S\ref{cm52}, and the
subspace topology on $(\pi^\sF)^{-1}(0)$ agrees with the
topology on $\M_\sX$ in Definition \ref{cm5def3}. More
generally, if $(s,\hat X)\in\M_\sX^\sF$ then $(\pi^\sF)^{-1}
(s)\subset\M_\sX^\sF$ is $\{s\}\t\M_{\sst\hat X}$ as a
topological space, where $\M_{\sst\hat X}$ is the moduli
space of deformations of $\hat X$ in~$(M,J^s,\om^s,\Om^s)$.
\label{cm7def3}
\end{dfn}

\noindent{\it Remarks.} Basically, $\M_\sX^\sF$ is the
family of pairs $(s,\hat X)$ where $s\in\F$ and $\hat X$ is
a compact SL $m$-fold in $M$ with conical singularities,
which is deformation equivalent to $X$ in a loose sense.
Note that $\M_\sX^\sF$ may not be {\it connected}. The fibres
$(\pi^\sF)^{-1}(s)$ of $\pi^\sF:\M_\sX^\sF\ra\F$ are (as
topological spaces) moduli spaces of compact SL $m$-folds in
$(M,J^s,\om^s,\Om^s)$ with conical singularities, deformation
equivalent to $X$, and with~$(\pi^\sF)^{-1}(0)=\M_\sX$.

The whole point of constructing $\M_\sX^\sF$, and its topology, is
that we can now make sense of the idea of a continuous family of
compact SL $m$-folds $\hat X$ in $M$ with conical singularities,
in which the underlying almost Calabi--Yau structure is allowed
to vary. That is, we can continuously deform $X$ not just in
$(M,J,\om,\Om)$ but also in $(M,J^s,\om^s,\Om^s)$ for~$s\in\F$.

\subsection{The main result for families $(M,J^s,\om^s,\Om^s)$}
\label{cm72}

Next we extend \S\ref{cm6} to the families case. Here are the
analogues of Definition \ref{cm6def1} and Proposition~\ref{cm6prop1}.

\begin{dfn} Let $\H_\sXp,\K_\sXp,\I_\sXp,k,p,{\bs\mu},\D_\sXp$
and $F$ be as in Definition \ref{cm6def1}. Define a map
$F^\sF:\ti\E^\sFp\t\D_\sXp\ra C^0(X')$ by
\e
\pi_*\bigl((\Phi^{s,\hat e}_\sXp)^*(\Im\Om)\vert_{\Ga(\be+\d f)}\bigr)
=F^\sF(s,\hat e,\be,f)\,\d V_g
\label{cm7eq6}
\e
for $(s,\hat e)\in\ti\E^\sFp$ and $(\be,f)\in\D_\sXp$. Then
$F^\sF(0,\hat e,\be,f)\equiv F(\hat e,\be,f)$ on~$\ti\E\t\D_\sXp$.
\label{cm7def4}
\end{dfn}

\begin{prop} In the situation above, suppose $(s,\hat e,\be,f)\in
\ti\E^\sFp\t\D_\sXp$ with\/ $F^\sF(s,\hat e,\be,f)=0$. Set\/
$\hat X'=\smash{\Phi^{s,\hat e}_\sXp}\bigl(\Ga(\be+\d f)\bigr)$
and\/ $\hat X=\hat X'\cup\{\hat x_1,\ldots,\hat x_n\}$, where
$\hat e=(\hat x_1,\ldots,\hat\up_n)$. Then $f\in C^\iy_{\bs\mu}(X')$
and\/ $\hat X$ is a compact SL\/ $m$-fold in $(M,J^s,\om^s,\Om^s)$
with conical singularities at\/ $\hat x_i$ with identifications
$\hat\up_i$, cones $C_i$ and rates $\mu_i$. Thus $(s,\hat X)$ lies
in $\V_\sX^\sF\subset\M_\sX^\sF$ in Definition \ref{cm7def3}.
Conversely, each\/ $\hat X$ in $\V_\sX^\sF$ comes from a unique
$(s,\hat e,\be,f)\in\ti\E^\sFp\t\D_\sXp$ with\/ $F^\sF(s,\hat e,
\be,f)=0$. Write $\Psi^\sF(s,\hat e,\be,f)=(s,\hat X)$. Then
$\Psi^\sF:(F^\sF)^{-1}(0)\ra\V_\sX^\sF$ is a homeomorphism,
with\/~$\Psi^\sF(0,e,0,0)=(0,X)$.
\label{cm7prop1}
\end{prop}

The modifications to the proof of Proposition \ref{cm6prop1}
are just trivial notational ones. We shall use Proposition
\ref{cm6prop2} as it is. The analogue of Proposition
\ref{cm6prop3} is

\begin{prop} In the situation above, $F^\sF$ maps
\e
F^\sF:\ti\E^\sFp\t\D_\sXp\ra\bigl\{u\in L^p_{k-2,{\bs\mu}-2}
(X'):\ts\int_{X'}u\,\d V_g=0\bigr\},
\label{cm7eq7}
\e
and this is a smooth map of Banach manifolds.
\label{cm7prop2}
\end{prop}

Again, the modifications to the proof are just trivial changes
in notation. We shall use all of \S\ref{cm62} as it is. The
point is that $F^\sF\vert_{s=0}\equiv F$, so the calculations
in \S\ref{cm62} about $\d F\vert_{(e,0,0)}$ immediately tell
us about the restriction of $\d F^\sF\vert_{(0,e,0,0)}$ to
the vector subspace with~$s=0$.

We can now prove the main result of this section, the analogue of
Theorem \ref{cm6thm2} for families, which describes $\M_\sX^\sF$
near~$(0,X)$.

\begin{thm} Suppose $(M,J,\om,\Om)$ is an almost Calabi--Yau
$m$-fold and\/ $X$ a compact SL\/ $m$-fold in $M$ with conical
singularities at\/ $x_1,\ldots,x_n$. Let\/ $\M_\sX,
X',\allowbreak
\I_\sXp,\allowbreak
\O_\sXp,\allowbreak
U,\allowbreak
\Phi$ and $\Xi$ be as in Theorem~\ref{cm6thm2}.

Suppose $\bigl\{(M,J^s,\om^s,\Om^s):s\in\F\bigr\}$ is
a smooth family of deformations of\/ $(M,J,\om,\Om)$, in the
sense of Definition \ref{cm2def6}, such that\/ $\iota_*(\ga)
\cdot[\om^s]=0$ for all\/ $\ga\in H_2(X,\R)$ and\/ $s\in\F$,
where $\iota:X\ra M$ is the inclusion, and\/ $[X]\cdot[\Im\Om^s]=0$
for all\/ $s\in\F$, where $[X]\in H_m(M,\R)$ and\/
$[\Im\Om^s]\in H^m(M,\R)$. Let\/ $\M_\sX^\sF$ and\/
$\pi^\sF:\M_\sX^\sF\ra\F$ be as in Definition~\ref{cm7def3}.

Then there exists an open neighbourhood\/ $U^\sF$ of\/ $(0,0)$
in $\F\t U$, a smooth map $\Phi^\sF:U^\sF\ra\O_\sXp$ with\/
$\Phi^\sF(0,u)\equiv\Phi(u)$, and a map $\Xi^\sF:\{(s,u)\in
U^\sF:\Phi^\sF(s,u)=0\}\ra\M_\sX^\sF$ with\/ $\Xi^\sF(0,u)\equiv
\bigl(0,\Xi(u)\bigr)$ and\/ $\pi^\sF\circ\Xi^\sF(s,u)\equiv s$,
which is a homeomorphism with an open neighbourhood of\/
$(0,X)$ in~$\M_\sX^\sF$.
\label{cm7thm3}
\end{thm}

\begin{proof} Recall that $0\in\F'\subseteq\F\subset\R^d$
and $\pi^\sF:\ti\E^\sFp\ra\F'$ is a submersion with fibres
open balls, and $\ti\E^\sFp\supset(\pi^\sF)^{-1}(0)
=\{0\}\t\ti\E$. Thus we can choose a smooth identification
of $\ti\E^\sFp$ with an open neighbourhood of $(0,0)$ in
$\F'\t T_e\ti\E\subset\R^d\t T_e\ti\E$ which identifies the
projections $\pi^\sF:\ti\E^\sFp\ra\F'$ and $\pi^\sF:\F'\t
T_e\ti\E\ra\F'$, and on $(\pi^\sF)^{-1}(0)=\{0\}\t\ti\E$
and $\{0\}\t T_e\ti\E$ agrees with the identification
between $\ti E$ and a subset of $T_e\ti\E$ chosen in the
proof of Theorem \ref{cm6thm2}. Define
\begin{align*}
Y^\sF&=\R^d\t\I_\sXp, \quad
Z=\O_\sXp\t T_e\ti\E\t\K_\sXp\t
L^p_{k,{\bs\mu}}(X'),\\
T&=\{u\in L^p_{k-2,{\bs\mu}-2}(X'):\ts\int_{X'}u\,\d V_g=0\bigr\}
\quad\text{and}\\
W^\sF&=\bigl\{(s,\be,\ga,\hat e,v,f)\!\in\!Y^\sF\t Z:
(s,\hat e)\!\in\!\ti\E^\sFp\!\subset\!\R^d\t T_e\ti\E,\;
(\be+\d v,f)\!\in\!\D_\sXp\bigr\}.
\end{align*}
Then $0\in Z$ is $(0,e,0,0)$. Choose any norms on the
finite-dimensional spaces $\R^d,\I_\sXp,\O_\sXp,
T_e\ti\E,\K_\sXp$, and use the usual norms on
$L^p_{k,{\bs\mu}}(X')$ and $T$. Then $Y^\sF,Z,T$ are Banach spaces,
and $W^\sF$ is an open neighbourhood of $(0,0)$ in $Y^\sF\t Z$, as
in Theorem~\ref{cm6thm1}.

Define a map $G^\sF:W^\sF\ra T$ by $G(s,\be,\ga,\hat e,v,f)=
\ga+F^\sF(s,\hat e,\be+\d v,f)$. This is a smooth map of Banach
manifolds, by Proposition \ref{cm7prop2}, and $G^\sF(0,0)=
G^\sF(0,0,0,e,0,0)=0$ as $F^\sF(0,e,0,0)=0$. The map
$\d G^\sF_{(0,0)}\vert_Z$ is given by
\e
\d G^\sF_{(0,0)}\vert_Z:(\ga,y,v,f)\!\mapsto\!\ga\!+\!
\d F^\sF_{(0,e,0,0)}(0,y,\d v,f)\!=\!\ga\!+\!\d F_{(e,0,0)}(y,\d v,f),
\label{cm7eq8}
\e
since $F^\sF\vert_{s=0}\equiv F$, as in Definition~\ref{cm7def4}.

Comparing \eq{cm7eq8} with \eq{cm6eq22} we see that $\d G^\sF_{
(0,0)}\vert_Z:Z\ra T$ agrees with $\d G_{(0,0)}\vert_Z:Z\ra T$
in the proof of Theorem \ref{cm6thm2}. Therefore $\d G^\sF_{
(0,0)}\vert_Z$ is an isomorphism of topological vector spaces
as in the proof of Theorem \ref{cm6thm2}, and we can apply
Theorem \ref{cm6thm1} to $Y^\sF,Z,T,W^\sF$ and $G^\sF$.
The rest of the proof is a straightforward modification of
that of Theorem~\ref{cm6thm2}.
\end{proof}

Here is the analogue of Corollary \ref{cm6cor1}. Note
the similarity to Theorem~\ref{cm2thm3}.

\begin{cor} Let\/ $(M,J,\om,\Om)$ be an almost Calabi--Yau
$m$-fold, $X$ a compact SL\/ $m$-fold in $M$ with stable
conical singularities, let\/ $\bigl\{(M,J^s,\om^s,\Om^s):s\in
\F\bigr\}$ be a smooth family of deformations of\/ $(M,J,\om,\Om)$
for $\F\subset\R^d$ with\/ $\iota_*(\ga)\cdot[\om^s]=0$ and\/
$[X]\cdot[\Im\Om^s]=0$ for all\/ $\ga\in H_2(X,\R)$ and\/ $s\in\F$,
and let\/ $\M_\sX,\M_\sX^\sF,\pi^\sF$, and\/ $\I_\sXp$ be as in
Theorem~\ref{cm7thm3}.

Then $\M_\sX^\sF$ is a smooth manifold of dimension
$d+\dim\I_\sXp$ and\/ $\pi^\sF:\M_\sX^\sF\ra\F$ a smooth
submersion. For all\/ $s\in\F$ sufficiently close to $0$ the
fibre $(\pi^\sF)^{-1}(s)$ is a nonempty smooth manifold of
dimension $\dim\I_\sXp$, with\/~$(\pi^\sF)^{-1}(0)=\M_\sX$.
\label{cm7cor1}
\end{cor}

Here $\pi^\sF:\M_\sX^\sF\ra\F$ is a {\it submersion\/} means that
$\pi^\sF:T_{\smash{\sst(s,\hat X)}}\M_\sX^\sF\!\ra\!T_s\F\!=\!\R^d$
is surjective for all $(s,\hat X)\in\M_\sX^\sF$. This follows near
$(0,X)\in\M_\sX^\sF$ as $\Xi^\sF$ is a diffeomorphism from
$U^\sF\subset\F\t U$ to a neighbourhood of $(0,X)\in\M_\sX^\sF$
which identifies the projections $\pi^\sF:\M_\sX^\sF\ra\F$ and
$\pi^\sF:\F\t U\ra\F$. Thus it holds near every $(s,\hat X)\in
\M_\sX^\sF$, by applying Theorem \ref{cm7thm3} with $(M,J,\om,
\Om)$ replaced by $(M,J^s,\om^s,\Om^s)$ and $X$ by~$\hat X$.

Corollary \ref{cm7cor1} implies the analogue of Theorem
\ref{cm2thm3} for compact SL $m$-folds $X$ in $M$ with
stable conical singularities. That is, it shows that
there are no local obstructions to deforming $X$ to
nearby almost Calabi--Yau structures $(J^s,\om^s,\Om^s)$
on $M$, except the obvious cohomological ones.

Here are the analogues of Definition \ref{cm6def3} and
Corollary~\ref{cm6cor2}.

\begin{dfn} Let $(M,J,\om,\Om)$ be an almost Calabi--Yau
$m$-fold, $X$ a compact SL $m$-fold in $M$ with conical
singularities, $\bigl\{(M,J^s,\om^s,\Om^s):s\in\F\bigr\}$
a smooth family of deformations of $(M,J,\om,\Om)$ for
$\F\subset\R^d$ with $\iota_*(\ga)\cdot[\om^s]=0$ and
$[X]\cdot[\Im\Om^s]=0$ for all $\ga\in H_2(X,\R)$ and
$s\in\F$, and let $\I_\sXp,\O_\sXp,U^\sF,\Phi$ and
$\Phi^\sF$ be as in Theorem~\ref{cm7thm3}.

We call $X$ {\it transverse in} $\F$ if the linear map
$\d\Phi^\sF\vert_{(0,0)}:\R^d\t\I_\sXp\ra\O_\sXp$ is
surjective. This definition is independent of the choices
made in defining $\I_\sXp,\O_\sXp,U^\sF$ and $\Phi^\sF$.
Since the restriction of $\d\Phi^\sF\vert_{(0,0)}$ to
$\I_\sXp\subset\R^d\t\I_\sXp$ is $\d\Phi\vert_0$, we
see that if $X$ is {\it transverse} in the sense of
Definition \ref{cm6def3} then it is also transverse in
$\F$, for any family~$\F$.
\label{cm7def5}
\end{dfn}

\begin{cor} Let\/ $(M,J,\om,\Om)$ be an almost Calabi--Yau
$m$-fold, $X$ a compact SL\/ $m$-fold in $M$ with conical
singularities, let\/ $\bigl\{(M,J^s,\om^s,\Om^s):s\in\F\bigr\}$
be a smooth family of deformations of\/ $(M,J,\om,\Om)$ for
$\F\subset\R^d$ with\/ $\iota_*(\ga)\cdot[\om^s]=0$ and\/
$[X]\cdot[\Im\Om^s]=0$ for all\/ $\ga\in H_2(X,\R)$ and\/
$s\in\F$, and let\/ $\M_\sX^\sF,\I_\sXp$ and\/ $\O_\sXp$
be as in Theorem \ref{cm7thm3}. Suppose $X$ is transverse
in $\F$. Then $\M_\sX^\sF$ is near $(0,X)$ a smooth manifold
of dimension $d+\dim\I_\sXp-\dim\O_\sXp$, and\/ $\pi^\sF:
\M_\sX^\sF\ra\F$ is a smooth map near~$(0,X)$.
\label{cm7cor2}
\end{cor}

Here Theorem \ref{cm7thm3} implies that near $(0,X)$ we can
identify $\M_\sX^\sF$ with a submanifold of $\F\t U$, and
$\pi^\sF$ then coincides with the projection $\pi^\sF:\F\t
U\ra\F$, so $\pi^\sF$ is smooth near $(0,X)$. Corollary
\ref{cm7cor2} will be important in \S\ref{cm9}, as we will show
that for any compact SL $m$-fold $X$ in $(M,J,\om,\Om)$ with
conical singularities, there exists a family of deformations
$\bigl\{(M,J^s,\om^s,\Om^s):s\in\F\bigr\}$ of $(M,J,\om,\Om)$
such that $X$ is transverse in~$\F$.

\section{Other extensions of Theorems \ref{cm6thm2} and \ref{cm7thm3}}
\label{cm8}

Section \ref{cm7} discussed the extension of the deformation
theory of \S\ref{cm5}--\S\ref{cm6} to {\it families} of almost
Calabi--Yau $m$-folds. We now briefly consider other possible
extensions of the theory, first to {\it immersed\/} rather
than {\it embedded\/} submanifolds, and secondly to ways in
which we can allow the SL cones $C_1,\ldots,C_n$ to vary over
the moduli spaces $\M_\sX,\M_\sX^\sF$, rather than being the
same at every point. Allowing the $C_i$ to vary reduces the
dimension of the obstruction space $\O_\sXp$, and so increases
the (expected) dimension of~$\M_\sX,\M_\sX^\sF$.

\subsection{Immersions}
\label{cm81}

So far, for simplicity, we have worked throughout with
{\it embedded\/} submanifolds. In fact, nearly everything
we have done can be generalized to {\it immersed\/}
submanifolds in an obvious way, with only trivial,
notational changes. Here are a few of the details
involved in doing this.

Instead of regarding compact SL $m$-folds $X$ in
$(M,J,\om,\Om)$ with conical singularities as subsets of $M$,
we instead regard $X$ as a {\it Riemannian manifold with
conical singularities}, in the sense of \cite[\S 2]{Joyc7},
together with an {\it isometric immersion} $\iota:X\ra M$,
which is locally but not necessarily globally injective.
The singular points $x_1,\ldots,x_n\in X$ are distinct,
but their images $\iota(x_1),\ldots,\iota(x_n)\in M$
may not be.

The $\Si_i$ become compact Riemannian manifolds with
isometric immersions $\Si_i\ra{\cal S}^{2m-1}$, and the
cones $C_i$ on $\Si_i$ become {\it Riemannian cones} in
the sense of \cite[\S 2.1]{Joyc7}, with isometric
immersions $C_i\ra\C^m$ which need not be locally injective
near 0. The $\Up_i$ can still be embeddings, but their
images may overlap. The $\phi_i,\iota_i,\Phi_{\sst C_i},
\Phi_{\sst X'}$, etc., should be taken to be immersions.

The only point the author is aware of where there is a
significant problem in changing from embeddings to
immersions is in the Geometric Measure Theory of
\cite[\S 6]{Joyc7}, in particular Theorem \ref{cm4thm9}
above, where the tangent cone $C$ must have $C\sm\{0\}$
a genuine embedded submanifold. However, we do not use
Theorem \ref{cm4thm9} in this paper, so this does not
affect the results of~\S\ref{cm5}--\S\ref{cm7}.

Suppose $C$ is an embedded SL cone in $\C^m$ with an isolated
singularity at 0, so that $\Si=C\cap{\cal S}^{2m-1}$ is a
compact $(m\!-\!1)$-manifold. If $\Si$ is not simply-connected
we may be able to take a {\it finite cover} $\pi:\ti\Si\ra\Si$.
Then $\ti\Si$ is an {\it immersed\/} minimal Legendrian
$(m\!-\!1)$-fold in ${\cal S}^{2m-1}$, with a corresponding
immersed SL cone $\ti C$ in~$\C^m$.

This construction considerably increases the supply of possible
SL cones available as model singularities in the immersed case.
It is particularly effective when $m=3$, as then $\Si$ is an
oriented Riemann surface of genus $g\ge 1$, and so admits
many finite covers. A similar phenomenon is described in
\cite[Th.~11.6]{Joyc4}, which constructs a large family of
immersed SL 3-folds in $\C^3$ diffeomorphic to ${\cal S}^1
\t\R^2$, which are asymptotic at infinity to the double
cover of an embedded SL $T^2$-cone in~$\C^3$.

\subsection{Cones $C_i$ with multiple ends}
\label{cm82}

The moduli spaces $\M_\sX$ and $\M_\sX^\sF$ defined in
\S\ref{cm5} and \S\ref{cm71} have the {\it same} set of
SL cones $C_1,\ldots,C_n$ (up to $\SU(m)$ equivalence)
for every $\hat X\in\M_\sX$ or $(s,\hat X)\in\M_\sX^\sF$.
There are various ways of relaxing this, and enlarging
the moduli spaces $\M_\sX,\M_\sX^\sF$ by allowing the
SL cones $C_i$ to vary. Consider the case in which
$\Si_1,\ldots,\Si_n$ are not all connected, so that
$b^0(\Si_i)>1$ for at least one $i$. We shall explain
two ways to generalize $\M_\sX$ and~$\M_\sX^\sF$.

The first way is to regard $X$ as an {\it immersed\/}
SL $m$-fold in $M$ with conical singularities, as in
\S\ref{cm81}. That is, instead of $X$ having $n$ singular
points $x_1,\ldots,x_n$, we regard it as having $\check n
=\sum_{i=1}^nb^0(\Si_i)$ distinct singular points $y_1,
\ldots,y_{\check n}$, where $\check n>n$, which happen to
be mapped to $n$ points in $M$ in groups of $b^0(\Si_i)$
for $i=1,\ldots,n$ by the immersion~$\iota:X\ra M$.

Essentially, we replace $X$ by $\check X=X'\cup\{y_1,
\ldots,y_{\check n}\}$, where each $y_i$ compactifies
one of the $\check n$ noncompact ends of $X'$. Then we
deform $\check X$ to get a moduli space $\check\M_\sX$
or $\check\M_\sX^\sF$ of immersed SL $m$-folds with
$\check n$ singular points. Note that for general
elements of $\check\M_\sX$ or $\check\M_\sX^\sF$,
there will be up to $\check n$ distinct singular
points in $M$, rather than just~$n$.

The second way is to retain the number $n$ of singular points,
but to allow the $b^0(\Si_i)$ components of $C_i'$ to move
around separately under $\SU(m)$ rotations. Let $\Si_i^j$ be
the connected components of $\Si_i$ for $j=1,\ldots,b^0(\Si_i)$,
and let $C_i^j$ be the cone on $\Si_i^j$ in $\C^m$, so
that~$C_i=\bigcup_{j=1}^{\smash{\sst b^0(\Si_i)}}C_i^j$.

Then in defining $\M_\sX,\M_\sX^\sF$ we allow the SL $m$-folds
$\hat X$ with conical singularities at $\hat x_1,\ldots,\hat x_n$
to have cones $\hat C_i=\bigcup_{j=1}^{\smash{\sst b^0(\Si_i)}}
B_i^jC_i^j$ for $B_i^j\in\SU(m)$ with $B_i^1=1$. This enlarges
the family of SL cones allowed in $\M_\sX,\M_\sX^\sF$, and so
enlarges~$\M_\sX,\M_\sX^\sF$.

In \S\ref{cm5}--\S\ref{cm7} we have to enlarge $\E$, etc.,
by including possible values of $B_i^j$ near 1 for $j>1$.
The main effect that this has on the final results is that
it {\it reduces} the dimension of the obstruction space
$\O_\sXp$, and thus {\it increases} the (expected) dimension
of $\M_\sX,\M_\sX^\sF$. The old formula \eq{cm6eq19} for
$\dim\O_\sXp$ should be replaced~by
\e
\dim\O_\sXp=\sum_{i=1}^n\Bigl(-2m+\sum_{j=1}^{b^0(\Si_i)}
\bigl(\sind(C_i^j)+2m\bigr)\Bigr).
\label{cm8eq1}
\e
If $b^0(\Si_i)>1$ one can show that this
does strictly reduce~$\dim\O_\sXp$.

The new obstruction space $\O_\sXp$ is a quotient of
the old by a vector subspace, which is the extra
obstructions we can overcome by moving the $C_i^j$
around separately under $\SU(m)$. The new infinitesimal
deformation space $\I_\sXp$ is the same as the old one.

There is one special case to be considered above. In
Definition \ref{cm3def3} and throughout we have assumed
that the SL cone $C_i$ has an isolated singularity at 0.
It could be that if $b^0(\Si_i)>1$ then some of the
$C_i^j$ above are SL planes $\R^m$ in $\C^m$, and thus
are nonsingular at 0, and so are not covered by
Definition~\ref{cm3def3}.

In this case \eq{cm3eq4} fails for $\Si_i^j={\cal S}^{m-1}$,
as $m_{\smash{\sst\Si_i^j}}(1)=m$. To compensate for this,
the appropriate value of $\sind(C_i^j)$ in \eq{cm8eq1} is
$\sind(C_i^j)=-m$. This is because the term $\sind(C_i^j)+2m$
in \eq{cm8eq1} contains a contribution $2m$ on the assumption
that $m_{\smash{\sst\Si_i^j}}(1)=2m$, and this has to be
reduced from $2m$ to~$m$.

\subsection{Families of special Lagrangian cones}
\label{cm83}

Let $X$ be a compact SL $m$-fold in $(M,J,\om,\Om)$ with
conical singularities at $x_1,\ldots,x_n$ with cones
$C_1,\ldots,C_n$. Here is a more general way of relaxing
the condition that the SL $m$-folds $\hat X$ in $\M_\sX$,
$\M_\sX^\sF$ must all have the {\it same} SL cones
$C_1,\ldots,C_n$ at their singular points.

Suppose ${\cal C}_i$ is a {\it smooth, connected family} of
{\it distinct\/} SL cones in $\C^m$ with $C_i\in{\cal C}_i$
for $i=1,\ldots,n$. Since we can always move cones through
$\SU(m)$ rotations by changing the identifications $\up_i$,
suppose for simplicity that ${\cal C}_i$ is closed under
the action of $\SU(m)$. Then in defining $\M_\sX,\M_\sX^\sF$
we allow the SL $m$-folds $\hat X$ with conical singularities
at $\hat x_1,\ldots,\hat x_n$ to have cones $\hat C_i\in
{\cal C}_i$ for~$i=1,\ldots,n$.

If ${\cal C}_i$ is the $\SU(m)$-orbit of $C_i$, then this
yields exactly the same moduli spaces $\M_\sX,\M_\sX^\sF$ as
in \S\ref{cm5}--\S\ref{cm7}. In the situation of \S\ref{cm82},
if $C_i=\smash{\bigcup_{j=1}^{\smash{\sst b^0(\Si_i)}}}C_i^j$
and we take ${\cal C}_i$ to be an open subset of the product of
the $\SU(m)$-orbits of $C_i^j$ for $j=1,\ldots,b^0(\Si_i)$, so
that ${\cal C}_i$ consists of cones $\hat C_i$ got by moving
the $C_i^j$ about independently with $\SU(m)$ rotations, then
this recovers the `second way' of~\S\ref{cm82}.

But if ${\cal C}_i$ contains {\it nontrivial\/} deformations
of $C_i$ not obtained by $\SU(m)$ rotations of the components
of $C_i'$, then this is a true generalization of the problem,
which will enlarge $\M_\sX,\M_\sX^\sF$ and their (expected)
dimension. Intuitively one might expect that special Lagrangian
cones are pretty rigid things and will not admit nontrivial
deformations in this way, so that there do not exist any
interesting families ${\cal C}_i$ to use in this construction.

However, at least when $m=3$, this is not true. There exists a
complicated theory which describes all {\it special Lagrangian
$T^2$-cones} in $\C^3$ using {\it integrable systems}, which is
described in McIntosh \cite{McIn} and the author \cite{Joyc5}.
It establishes a 1-1 correspondence between SL $T^2$-cones in
$\C^3$ up to isometry and collections of {\it spectral data},
including a Riemann surface $Y$ with even genus called the
{\it spectral curve}, and a holomorphic line bundle~$L\ra Y$.

As \cite[\S 4.2]{McIn} and \cite[\S 4.3]{Joyc5}, it turns
out that an SL $T^2$-cone with spectral curve $Y$ of genus
$2d\ge 4$ is part of a smooth $(d\!-\!2)$-dimensional family
of SL $T^2$-cones up to isometries of $\C^3$, which have the
same spectral curve $Y$ but varying line bundles $L\ra Y$.
Ian McIntosh (personal communication) and Emma Carberry have
recently announced a proof of the existence of SL $T^2$-cones
with spectral curves of every even genus. Thus there exist
smooth families ${\cal C}_i$ of SL $T^2$-cones in $\C^3$
with arbitrarily high dimension, to which we can apply this
deformation theory.

The main changes to the final results are that we
replace the definition of $\sind(C_i)$ in \eq{cm3eq5} by
$\sind_{\sst{\cal C}_i}(C_i)=N_{\sst\Si_i}(2)-b^0(\Si_i)-2m
-\dim{\cal C}_i$, the {\it stability index of\/ $C_i$ in\/}
${\cal C}_i$, and then the old formula \eq{cm6eq19} for
$\dim\O_\sXp$ should be replaced by $\dim\O_\sXp=\sum_{i=1}^n
\sind_{\sst{\cal C}_i}(C_i)$. The new infinitesimal deformation
space $\I_\sXp$ is the same as the old one.

\section{Transversality and genericity results}
\label{cm9}

Finally we discuss the question: if $(M,J,\om,\Om)$ is a
{\it generic} almost Calabi--Yau $m$-fold, are moduli spaces
$\M_\sX$ of compact SL $m$-folds $X$ in $M$ with conical
singularities necessarily smooth?

Consider what we mean by {\it generic} here. The conditions
$\iota_*(\ga)\cdot[\om]=0$ for $\ga\in H_2(X,\R)$ and $[X]
\cdot[\Im\Om]=0$ mean that when $[\om],[\Im\Om]$ are
generic there {\it will not exist\/} any such SL $m$-folds $X$
in $(M,J,\om,\Om)$. Thus, choosing $(M,J,\om,\Om)$ generically
in the family of all almost Calabi--Yau $m$-folds is too strong.
Instead, we shall require only that $\om$ {\it is generic in
its K\"ahler class}.

That is, given an almost Calabi--Yau $m$-fold $(M,J,\om,\Om)$
containing a compact SL $m$-fold $X$ with conical singularities,
we consider generic perturbations $(M,J,\check\om,\Om)$ with
$\check\om=\om+\d(J\,\d f)$ for some {\it K\"ahler potential\/}
$f\in C^\iy(M)$, so that $[\check\om]=[\om]\in H^2(M,\R)$. Then
there are no cohomological obstructions to the existence of SL
$m$-folds $\check X$ with conical singularities in $(M,J,\check
\om,\Om)$ isotopic to $X$, and we wish to know whether the
moduli space $\check\M_\sX$ of such $\check X$ is smooth.

We begin by showing that for any compact SL $m$-fold $X$ with
conical singularities, there exists a family of deformations
$\F$ with $X$ transverse in~$\F$.

\begin{thm} Let\/ $(M,J,\om,\Om)$ be an almost Calabi--Yau
$m$-fold and\/ $X$ a compact SL\/ $m$-fold in $M$ with conical
singularities. Let\/ $\I_\sXp,\O_\sXp$ be as in \S\ref{cm6}.
Then there exists a smooth family of deformations
$\bigl\{(M,J,\om^s,\Om):s\in\F\bigr\}$ of\/ $(M,J,\om,\Om)$
with\/ $[\om^s]=[\om]\in H^2(M,\R)$ for all\/ $s\in\F$, such
that\/ $X$ is transverse in $\F$, in the sense of Definition
\ref{cm7def5}, and\/ $\dim\F=\dim\O_\sXp$. Hence the moduli
space $\M_\sX^\sF$ of\/ \S\ref{cm7} is a manifold near~$(0,X)$.
\label{cm9thm1}
\end{thm}

\begin{proof} Use the notation of \S\ref{cm6}--\S\ref{cm7}.
Recall from Definition \ref{cm6def2} that $\O_\sXp$ consists of
smooth, compactly-supported functions $v$ on $X'$ with $\int_{X'}
v\,\d V_g=0$. Since $H^m_{\rm cs}(X',\R)=0$, we see that each
such $v$ may be written as $\d^*(\psi^m\al)$ for $\al$ a smooth,
compactly-supported 1-form on $X'$. Let $d=\dim\O_\sXp$, and
choose smooth, compactly-supported 1-forms $\al_1,\ldots,\al_d$
on $X'$ with
\e
\O_\sXp=\ban{\d^*(\psi^m\al_1),\ldots,\d^*(\psi^m\al_d)}.
\label{cm9eq1}
\e

Suppose $f\in C^\iy(M)$ with $f\vert_{X'}\equiv 0$. Then
$\d f\vert_{X'}\in C^\iy(\nu^*)$, where $\nu\ra X'$ is the
normal bundle to $X'$ in $M$. But the complex structure
$J$ induces an isomorphism $\nu\cong TX'$, so we can regard
$\d f\vert_{X'}$ as an element of $C^\iy(T^*X')$, that is,
a 1-form on $X'$.

Choose smooth functions $f_1,\ldots,f_d\in C^\iy(M)$ such that
$f_j\vert_{X'}\equiv 0$, and $f_j$ is supported on a small open
neighbourhood $U_j$ in $M$ of the support of $\al_j$ in $X'$
with $x_i\notin U_j$ for $i=1,\ldots,n$, and $\d f_j\vert_{X'}$
is identified with $\al_j$ under the isomorphism $C^\iy(\nu^*)
\cong C^\iy(T^*X')$ above, for $j=1,\ldots,d$. It is easy to show
that this is possible.

For $s=(s_1,\ldots,s_d)\in\R^d$, define a closed real (1,1)-form
$\om^s$ on $M$ by
\e
\om^s=\om+\ts\sum_{j=1}^ds_j\,\d\bigl(J(\d f_j)\bigr).
\label{cm9eq2}
\e
Choose an open neighbourhood $\F$ of 0 in $\R^d$ such that
$\om^s$ is the K\"ahler form of a K\"ahler metric $g^s$ on
$(M,J)$ for all $s\in\F$. This is true for small $s\in\R^d$.
Then $\bigl\{(M,J,\om^s,\Om):s\in\F\bigr\}$ is a smooth
family of deformations of $(M,J,\om,\Om)$, in the sense of
Definition~\ref{cm2def6}.

The definition of $f_j$ implies that $(J\,\d f_j)\vert_{X'}=\al_j$.
Thus \eq{cm9eq2} gives
\e
\ts\om^s\vert_{X'}=\sum_{j=1}^ds_j\,\d\al_j.
\label{cm9eq3}
\e
Applying Theorem \ref{cm4thm7} gives $0\in\F'\subseteq\F$ and
family of maps $\Phi^s_\sXp:U_\sXp\ra M$ for $s\in\F'$ with
$(\Phi_\sXp^s)^*(\om^s)=\hat\om$. Identifying $X'$ with the zero
section in $U_\sXp$, we see from \eq{cm9eq2} and \eq{cm9eq3} that
\e
(\Phi_\sXp^s)^*(\om)\vert_{X'}=-\sum_{j=1}^ds_j\,\d\al_j
+O\bigl(\ms{s}\bigr)\quad\text{for small $s\in\F'$.}
\label{cm9eq4}
\e
As the restriction of $\hat\om$ on $U_\sXp$ to the graph
$\Ga(\al)$ of a 1-form $\al$ is $-\d\al$, examining the
proof of Theorem \ref{cm4thm7} in \cite{Joyc7} we find that
we can choose $\Phi^s_\sXp$ such that
\e
\Phi_\sXp^s(x)=\Phi_\sXp\bigl(\ts\sum_{j=1}^ds_j\,\al_j\bigr)
+O\bigl(\ms{s}\bigr)\quad\text{for $x\in X'$ and small $s\in\F'$.}
\label{cm9eq5}
\e
That is, the image of the zero section under $\Phi_\sXp^s$
approximates the image of the graph of $\sum_{j=1}^ds_j\,\al_j$
under~$\Phi_\sXp$.

The proof of Proposition \ref{cm2prop2} now shows that
\e
(\Phi_\sXp^s)^*(\Im\Om)\vert_{X'}=-\sum_{j=1}^ds_j\,
\d^*(\psi^m\al_j)\d V_g+O\bigl(\ms{s}\bigr)
\quad\text{for small $s\in\F'$.}
\label{cm9eq6}
\e
But $\Phi_\sXp^{s,e^s}=\Phi_\sXp^s$ in Theorem \ref{cm7thm1}
and \eq{cm7eq6} in Definition \ref{cm7def4} imply that
\e
(\Phi^s_\sXp)^*(\Im\Om)\vert_{X'}=F^\sF(s,e^s,0,0)\d V_g.
\label{cm9eq7}
\e
Combining equations \eq{cm9eq1}, \eq{cm9eq6} and \eq{cm9eq7}
shows that the projection to $\O_\sXp$ of the derivative
$\d\F^\sF\vert_{(0,e,0,0)}$ is surjective. It easily follows
that in Theorem \ref{cm7thm3}, the map $\d\Phi^\sF\vert_{(0,0)}:
\R^d\t\I_\sXp\ra\O_\sXp$ is surjective. Hence $X$ is transverse
in $\F$ by Definition \ref{cm7def5}. The last part follows from
Corollary~\ref{cm7cor2}.
\end{proof}

Let $F:P\ra Q$ be a smooth map between finite-dimensional
manifolds. Recall that $q\in Q$ is called a {\it critical
value} of $F$ if $q=F(p)$ for some $p\in P$ for which
$\d F\vert_p:T_pP\ra T_qQ$ is not surjective. Points
$q\in Q$ which are not critical values are called
{\it regular values}. Then {\it Sard's Theorem}
(see Bredon \cite[\S II.6 \& App.~C]{Bred} for a proof)
says that the set of critical values of $F$ is of measure
zero in $Q$. Thus, almost all points in $Q$ are regular values.

This is important because if $q\in Q$ is a regular value then
$F^{-1}(q)$ is a {\it submanifold\/} of $Q$, of dimension
$\dim P-\dim Q$. Now in Theorem \ref{cm9thm1} we know that
$\M_\sX^\sF$ is a manifold and $\pi^\sF:\M_\sX^\sF\ra\F$ a
smooth map near $(0,X)$. Thus Sard's Theorem shows that
$(\pi^\sF)^{-1}(s)$ is a manifold near $(0,X)$ for small
generic $s\in\F$. So we prove:

\begin{cor} In the situation of Theorem \ref{cm9thm1}, for small
generic $s\in\F$ the moduli space $\M_\sX^s=(\pi^\sF)^{-1}(s)\subset
\M_\sX^\sF$ of deformations of\/ $X$ in $(M,J,\om^s,\Om)$ is near
$(0,X)$ a manifold of dimension~$\dim\I_\sXp-\dim\O_\sXp$.
\label{cm9cor1}
\end{cor}

If $\dim\I_\sXp-\dim\O_\sXp<0$ then $\M_\sX^s$ is empty near
$(0,X)$ for small generic $s$. We can generalize Theorem
\ref{cm9thm1} in the following way. As transversality is an
open condition, $\hat X$ is transverse to $\F$ for $\hat X$
in an open neighbourhood of $X$ in $\M_\sX$. In the same way,
for each $\ti X\in\M_\sX$ we can construct a family of
deformations $\F_{\sst\ti X}$ of $(M,J,\om,\Om)$ and an open
neighbourhood of $\ti X$ in $\M_\sX$ in which all $\hat X$
are transverse to~$\F_{\sst\ti X}$.

Let $W\subseteq\M_\sX$ be compact. Taking a finite subcover of
$W$ from this collection of open neighbourhoods in $\M_\sX$, we
get families of deformations $\F_1,\ldots,\F_l$ of $(M,J,\om,\Om)$
such that every $\hat X\in W$ is transverse in $\F_j$ for some
$j=1,\ldots,l$. Choose a family $\F$ of deformations of
$(M,J,\om,\Om)$ containing open neighbourhoods of 0 in
$\F_1,\ldots,\F_l$. This is easily done, as the $\F_j$ are open
neighbourhoods of $\om$ in affine subspaces ${\cal A}_1,\ldots,
{\cal A}_l$ of the K\"ahler class of $\om$, and we can take $\F$
to be an open neighbourhood of $\om$ in the affine subspace
spanned by ${\cal A}_1,\ldots,{\cal A}_l$. Then all $\hat X\in W$
are transverse in $\F$, giving:

\begin{thm} Let\/ $(M,J,\om,\Om)$ be an almost Calabi--Yau
$m$-fold and\/ $X$ a compact SL\/ $m$-fold in $M$ with conical
singularities. Let\/ $\M_\sX,\I_\sXp,\O_\sXp$ be as in
\S\ref{cm5}--\S\ref{cm6}, and suppose $W\subseteq\M_\sX$ is a
compact subset. Then there exists a smooth family of deformations
$\bigl\{(M,J,\om^s,\Om):s\in\F\bigr\}$ of\/ $(M,J,\om,\Om)$
with\/ $[\om^s]=[\om]\in H^2(M,\R)$ for all\/ $s\in\F$, such
that\/ $\hat X$ is transverse in $\F$ for all\/ $\hat X\in W$.
Hence the moduli space $\M_\sX^\sF$ of\/ \S\ref{cm7} is a
manifold near~$\{0\}\t W$.
\label{cm9thm2}
\end{thm}

The analogue of Corollary \ref{cm9cor1} is:

\begin{cor} In the situation of Theorem \ref{cm9thm2}, for small
generic $s\in\F$ the moduli space $\M_\sX^s=(\pi^\sF)^{-1}(s)\subset
\M_\sX^\sF$ of deformations of\/ $X$ in $(M,J,\om^s,\Om)$ is near
$\{0\}\t W$ a manifold of dimension~$\dim\I_\sXp-\dim\O_\sXp$.
\label{cm9cor2}
\end{cor}

Roughly speaking, Corollaries \ref{cm9cor1} and \ref{cm9cor2}
imply that for a small generic perturbation $(M,J,\check\om,\Om)$
of $(M,J,\om,\Om)$ in the same K\"ahler class, the perturbed
moduli space $\check\M_\sX$ is a manifold near $X$, or more
generally near a compact subset $W$ of $\M_\sX$. Of course,
$X$ and $W$ do not lie in $\check\M_\sX$, but the idea does
make sense. We conjecture that if $\check\om$ is sufficiently
generic then $\check\M_\sX$ is a manifold everywhere.

\begin{conj} Let\/ $(M,J,\om,\Om)$ be an almost Calabi--Yau
$m$-fold and\/ $X$ a compact SL\/ $m$-fold in $M$ with conical
singularities, and let\/ $\I_\sXp,\O_\sXp$ be as in \S\ref{cm6}.
Then for a second category subset of K\"ahler forms $\check\om$
in the K\"ahler class of\/ $\om$, the moduli space $\check\M_\sX$
of compact SL\/ $m$-folds $\hat X$ with conical singularities in
$(M,J,\check\om,\Om)$ isotopic to $X$ is a manifold of
dimension~$\dim\I_\sXp-\dim\O_\sXp$.
\label{cm9conj}
\end{conj}

Recall that a subset of a topological space is of {\it second
category} if it can be written as the intersection of a countable
number of open dense sets. Using the {\it Baire category theorem}
one can show that second category subsets of the K\"ahler class of
$\om$ are dense. Thus, the conjecture implies that $\check\M_\sX$
is smooth for generic~$\check\om$.

As a countable intersection of second category subsets is second
category, the conjecture also implies that by choosing $\check\om$
generically we can make a countable number of moduli spaces
$\check\M_{\sst X_1},\check\M_{\sst X_2},\ldots$ simultaneously
smooth. However, we have not extended Conjecture \ref{cm9conj} to
the tempting, much simpler statement that for generic $\check\om$,
{\it all\/} the moduli spaces $\check\M_\sX$ are smooth.

This is because, as in \S\ref{cm83}, there can exist smooth,
positive-dimensional families of SL cones in $\C^m$ which are 
distinct under $\SU(m)$ transformations. Now with the definitions
of \S\ref{cm5}, every $\hat X\in\M_\sX$ has the same cones
$C_1,\ldots,C_n$. If these cones $C_i$ are allowed to vary in
positive-dimensional families, we would get corresponding
uncountable families of moduli spaces $\check\M_\sX^t$, and it
is too much to expect all of these to be simultaneously smooth.

Results similar to Conjecture \ref{cm9conj} are proved by
Donaldson and Kronheimer \cite[\S 4.3]{DoKr} for moduli
spaces of instantons on 4-manifolds w.r.t.\ a generic $C^l$
metric, and by McDuff and Salamon \cite[\S 3]{McSa1} for
smoothness of moduli spaces of pseudo-holomorphic curves
on a symplectic manifold w.r.t.\ a generic $C^l$ or smooth
almost complex structure.

Following these proofs, the author has a sketch proof of a
version of Conjecture \ref{cm9conj} using $C^l$ K\"ahler forms
$\check\om$ rather than smooth K\"ahler forms, for large
$l\ge 3$. It involves messy issues in infinite-dimensional
analysis, so we will not give it. The reason for using $C^l$
K\"ahler forms is to be able to apply the Sard--Smale Theorem,
a version of Sard's Theorem for Banach manifolds. The author
cannot yet see how to extend this to smooth K\"ahler forms.

\end{document}